\documentclass[12pt]{article}
\usepackage{latexsym,amsmath,amscd,amssymb,graphics}
\usepackage{enumerate}
     
\newcommand{\bfi}{\bfseries\itshape}

\makeatletter

\@addtoreset{figure}{section}
\def\thefigure{\thesection.\@arabic\c@figure}
\def\fps@figure{h, t}
\@addtoreset{table}{bsection}
\def\thetable{\thesection.\@arabic\c@table}
\def\fps@table{h, t}
\@addtoreset{equation}{section}

\makeatother

\newcommand{\one}{\ensuremath{L^1(\mathcal{M})}}
\newcommand{\two}{\ensuremath{L^2(\mathcal{M})}}
\newcommand{\bounded}{\ensuremath{L^\infty(\mathcal{M})}}
\newcommand{\compact}{\ensuremath{\mathcal{K(M)}}}

\newcommand{\tr}{\ensuremath{\operatorname{tr}}}
\newcommand{\invertible}{\ensuremath{GL^\infty(\mathcal{M})}}
\newcommand{\shermone}{\ensuremath{\mathfrak{u}^1(\mathcal{M})}}
\newcommand{\hermone}{\ensuremath{\mathfrak{h}^1(\mathcal{M})}}
\newcommand{\sherm}{\ensuremath{\mathfrak{u}^\infty(\mathcal{M})}}
\newcommand{\herm}{\ensuremath{\mathfrak{h}^\infty(\mathcal{M})}}
\newcommand{\unitary}{\ensuremath{U^\infty(\mathcal{M})}}
\newcommand{\norm}[1]{\left\Vert#1\right\Vert}
\newcommand{\im}{\ensuremath{\operatorname{im}}}
\newcommand{\ad}{\ensuremath{\operatorname{ad}}}
\renewcommand{\b}{\mathfrak{b}}

\newcommand{\C}{\mathbb C}

\newcommand{\m}{\ensuremath{\mathfrak m\;}}
\newcommand{\M}{\ensuremath{\mathcal M\;}}
\begin{document}

\newtheorem{theorem}{Theorem}[section]
\newtheorem{definition}[theorem]{Definition}
\newtheorem{lemma}[theorem]{Lemma}
\newtheorem{remark}[theorem]{Remark}
\newtheorem{proposition}[theorem]{Proposition}
\newtheorem{corollary}[theorem]{Corollary}
\newtheorem{example}[theorem]{Example}
      
\makeatletter

\title{Banach Lie-Poisson spaces and reduction}
\author{Anatol Odzijewicz$^{1}$ and Tudor S. Ratiu$^{2}$} 
\addtocounter{footnote}{1} 
\footnotetext{Institute of Physics, University of Bialystok, Lipowa 41, PL-15424
Bialystok, Poland. \texttt{aodzijew@labfiz.uwb.edu.pl}}
\addtocounter{footnote}{1} 
\footnotetext{Institut de Math\'ematiques Bernoulli,
\'Ecole Polytechnique F\'ed\'erale de Lausanne. CH--1015 Lausanne.
Switzerland. \texttt{Tudor.Ratiu@epfl.ch} }
\date{April 16, 2003}
\maketitle 

\makeatother

\maketitle



\begin{abstract}
The category of Banach Lie-Poisson spaces is introduced and studied. It is
shown that the category of $W^*$-algebras can be considered as one of
its subcategories. Examples and applications of Banach Lie-Poisson spaces to
quantization and integration of Hamiltonian systems are given. The 
relationship between classical and quantum reduction is discussed.
\end{abstract}

\tableofcontents
\section{Introduction}
\label{section: introduction}

This paper investigates the foundations of Banach Poisson differential geometry,
including such topics as Banach Lie-Poisson spaces, classical and quantum reduction,
integration and quantization of Hamiltonian systems with the aid of the momentum
map. We were inspired to study this circle of problems due to the appearance of 
formal Poisson structures in a large number of works devoted to the integration
of infinite dimensional systems and the crucial role played by the momentum map
in these approaches. 

The notion of a Lie-Poisson space is as old as the concept of a Lie algebra
and both were introduced simultaneously by Lie [1890]. A
Lie-Poisson space is a Poisson vector space with the
property that its dual is invariant under the Poisson bracket, which is
equivalent to the statement that the Poisson bracket is linear. In the
finite dimensional case the notions of Lie algebras and Lie-Poisson spaces
are equivalent in the sense that for any Lie algebra $\mathfrak{g}$ its dual
$\mathfrak{g}^\ast$ is a Lie-Poisson space and, conversely, given a
Lie-Poisson space its dual is a Lie algebra. This is so because finite dimensional
vector spaces are reflexive, the operation of taking the dual defines an 
isomorphism between these two categories. To generalize this to infinite
dimensions, it is reasonable to assume that a Lie-Poisson space is a
Banach space $\mathfrak{b}$ endowed with a Poisson bracket $\{\cdot,
\cdot\}$ such that the bracket of any two linear continuous functions is
again a linear continuous function. This implies that $(\mathfrak{b}^\ast,
\{\cdot, \cdot\})$ is a Banach Lie algebra. In order to preserve the
correspondence between Banach Lie-Poisson spaces and Banach Lie algebras
it is necessary to restrict to those Banach Lie algebras $(\mathfrak{g},
[\cdot, \cdot])$ that admit a predual $\mathfrak{g}_\ast$ and satisfy in
addition the condition that $\operatorname{ad}^\ast_{\mathfrak{g}}:
\mathfrak{g}^\ast \rightarrow \mathfrak{g}^\ast$ preserves the predual
$\mathfrak{g}_\ast$. Thus, in the infinite dimensional case, Banach
Lie-Poisson spaces form a subcategory of the category of Banach Lie
algebras. A crucial example is the Banach space $\one$ of linear trace
class operators on a separable Hilbert space $\mathcal{M}$ which is
predual to the Banach Lie algebra $\bounded$ of all linear bounded
operators on $\mathcal{M}$. As far as we know, the Lie-Poisson structure on
$\one$ was first found by Bona [2000].

Momentum maps are an efficient way to encode integrals of motion for a
Hamiltonian system. In its modern formulation due to Kostant [1966] and
Souriau [1966], [1967] a momentum map is naturally associated to an
infinitesimal Poisson action of a Lie algebra on a Poisson manifold and it
maps the phase space to the dual of the Lie algebra of symmetries. It turns
out, that in finite dimensions, a momentum map is characterized by the
property that it is Poisson, when one endows the dual of the Lie algebra
with the Lie-Poisson structure (see, e.g. Marsden and Ratiu [1994] for a
proof of this fact). 

In infinite dimensions, due to existence of non-reflexive Banach spaces,
we will define a momentum map to be a Poisson map from a Banach Poisson
manifold to a Banach Lie-Poisson space, which can always be considered as
the predual of the Banach Lie algebra of symmetries. 
It is shown that the momentum map so defined has all the usual properties,
such as being conserved along the flow of any symmetry invariant
Hamiltonian vector field (Noether's theorem).  Like in finite dimensions, also in
the infinite dimensional case the notion of momentum map is an important
tool in the study of Hamiltonian systems. For example, the
knowledge of momentum maps leads to integrals of motion of the considered
Hamiltonian system, as will be illustrated here through the example of the
infinite Toda lattice.

In the special case when one assumes that the momentum map is an
injective immersion and its range is linearly dense in the target Banach
Lie-Poisson space one discovers that it is the coherent states map in the
sense of Odzijewicz [1992]. So, it can be used to quantize the system
under consideration. This method of quantization, called Ehrenfest
quantization, is a natural unification of the Kostant-Souriau geometric
(Kostant [1970], Souriau [1966], [1967]) and $*$-product quantization; for
details see Odzijewicz [1992] and \S \ref{section: momentum maps and
reduction}. 

The structure of the paper is as follows. In \S \ref{section: banach poisson 
manifolds} the notion of a Banach Poisson manifold is introduced modeled
on the example of a strong symplectic manifold. Its elementary properties 
are presented as well as some comments on the compatibility of the Poisson 
structure with almost complex, complex, and holomorphic structures. 

Classical reduction for Banach Poisson manifolds is discussed 
in \S \ref{section: classical reduction}. The Poisson reduction theorem 
of Marsden and Ratiu [1986] and its consequences are generalized to the Banach
manifold context. 

Banach Lie-Poisson spaces and their properties are analyzed in 
\S \ref{section: banach lie-poisson spaces}. Linear
continuous Poisson maps are studied in detail. The realification and
complexification of a Banach Lie-Poisson space is also presented. The
upshot of this section is the establishment of an isomorphism between
the category of Banach Lie-Poisson spaces and a specific subcategory of
Banach Lie algebras. 

The entirety of \S \ref{section: preduals} is devoted to one crucial
example: the predual of a $W^\ast$-algebra and the dual to a
$C^\ast$-algebra are naturally a Banach Lie-Poisson spaces. As a
consequence it is  shown that various spaces related to operator
algebras (for example the space of
Hermitian trace class operators on a separable Hilbert space)
are Banach Lie-Poisson spaces. 

In \S \ref{section: quantum reduction} we show that quantum measurement
operation in the sense of von Neumann can be considered as a Poisson
projection. We shall give examples of other physically important Poisson 
projections. These examples justify the interpretation  of Poisson projection
as a quantum reduction procedure.

The internal structure of Banach Lie-Poisson spaces is presented in 
\S \ref{section: symplectic leaves}. If the Lie
algebra of a Banach Lie group admits a predual which is invariant
under the coadjoint representation it is shown that a large class of 
coadjoint orbits in the
predual, which is naturally a Banach Lie-Poisson space, are symplectic
leaves in a weak sense: they are weak symplectic manifolds and are
weakly immersed submanifolds (the inclusion is  smooth and has injective 
derivative, but no splitting condition, or even a closed range condition
on the derivative, usually imposed in the definition of an immersion,
holds). Among these orbits a subclass is determined for which the
symplectic form is strong and the orbit is injectively immersed.
The section ends with the standard example of a dual pair based on
the cotangent bundle (in this case on the precotangent bundle) that
illustrates that our definition of a Banach Poisson manifold is violated
in this important case and that once one leaves the category of
$W^\ast$-algebras a weakening of this notion will be needed.

Section \ref{section: momentum maps and reduction} introduces 
momentum maps as Poisson maps from a Banach
Poisson manifold to a Banach Lie-Poisson space. It is shown that the
coherent states map is a momentum map with certain special properties.
In this way the quantization procedure based on the coherent states map
has Banach Poisson geometrical interpretation. The relationship between 
classical and quantum reduction is explored. Using both
procedures of reduction, classical and quantum, one can construct a new 
momentum map from a given one. The description of the infinite Toda lattice
in Banach Poisson geometrical terms is presented. Among others, it is shown
that the Flaschka transformation is a densely defined momentum map of some Banach weak symplectic space
into the Banach Lie-Poisson space of the lower triangular trace class operators.

\section{Banach Poisson manifolds}
\label{section: banach poisson manifolds}

Throughout the paper, given a Banach space $\b$, the notation
$\b^\ast$ will always mean the Banach space dual to
$\mathfrak{b}$. Given $x \in \b^\ast$ and $b \in \b$, we shall
denote by $\langle x, b
\rangle$ the value of $x$ on $b$. Thus $\langle \cdot,
\cdot \rangle : \b^\ast \times \b \rightarrow \mathbb{R}$ (or
$\mathbb{C}$, depending on whether we work with real or complex
Banach spaces and functions) will  denote the natural bilinear
continuous duality pairing between $\b$ and its dual $\b^\ast$.

A real finite dimensional Poisson manifold is a pair
$(P, \{\cdot, \cdot\})$ consisting of a manifold $P$ whose space
of Fr\'echet smooth functions is endowed with a
Lie algebra structure $\{\cdot, \cdot \}$ satisfying the Leibniz
property in each factor; this bilinear operation $\{\cdot, \cdot\}$ is called a
Poisson bracket. As we shall discuss below, this
definition is not appropriate in infinite dimensions and a more
stringent condition needs to be imposed.

To see this, assume that on the space $C^\infty(P)$ of smooth functions on the
infinite dimensional smooth Banach manifold $P$
there is a Poisson bracket $\{\cdot, \cdot \}$. Due to the Leibniz property, 
the value of the Poisson bracket at a given point $p\in P$ depends only on the
differentials $d f(p)$, $d g (p) \in T_p ^\ast P$ which implies
that there is a smooth section $\varpi$ of the vector bundle
$\bigwedge^2 T^{\ast \ast}P$ satisfying
\[
\{f, g\} = \varpi(df, dg ).
\]
This means that for each $p\in P$ the map $\varpi_p : T_p^\ast P
\times T_p ^\ast P \rightarrow \mathbb{R}$  is a
continuous bilinear antisymmetric map that depends smoothly on the
base point $p$. In addition, denoting by $[\cdot, \cdot]_S$ the Schouten 
bracket on skew symmetric contravariant tensors, the equality (see e.g. Marsden
and Ratiu [1998], \S 10.6) 
\[
\{\{f, g\}, h\} +  \{\{g, h\}, f\} + \{\{h, f\}, g\} =
i_{[\varpi, \varpi]_S}(df \wedge dg \wedge dh),
\]
shows that the Jacobi identity is equivalent to $[\varpi, \varpi]_S = 0$, which
is an additional differential quadratic condition on $\varpi$. 

Let $\sharp: T^\ast P \rightarrow T^{\ast\ast} P$ be the bundle map covering
the identity defined by $\sharp_p(dh (p)) : = \varpi(\cdot, dh)(p)$, that is, 
$\sharp_p(dh (p)) (dg(p)) = \{g , h\}(p)$, for any locally defined functions 
$g$ and $h$.

Denote by $\mathfrak{b}$ the Banach space modeling the Banach
manifold $P$. Thus $T_pP \cong \mathfrak{b}$, $T_p^\ast P  \cong
\mathfrak{b}^\ast$, and $T_p^{\ast \ast} P  \cong
\mathfrak{b}^{\ast \ast}$. If $\mathfrak{b}$ is not reflexive,
that is, $\mathfrak{b} \subset \mathfrak{b}^{\ast \ast}$ and
$\mathfrak{b} \neq \mathfrak{b}^{\ast \ast}$, then
\[
X_f:= \varpi(\cdot, df) = \sharp(df), \quad \text
{or,~as~a~derivation~on~functions,} \quad X_f = \{ \cdot, f \}
\]
is a smooth section of $T^{\ast \ast} P$ and hence is not, in
general, a vector field on $P$. In analogy with the finite
dimensional case, we want $X_f$ to be the Hamiltonian vector field
defined by the function $f$. In order to achieve this, we are
forced to make the assumption that the Poisson bracket on $P$ satisfies the condition
$\sharp (T^\ast P) \subset TP \subset T^{\ast \ast}P$. Thus we give the following
definition.

\begin{definition}
\label{Poisson manifold definition}
A {\bfi Banach Poisson manifold\/} is a pair $(P, \{\cdot, \cdot\})$ consisting of a
smooth Banach manifold and a bilinear operation $\{\cdot, \cdot \}$ satisfying the
following conditions:
\begin{enumerate}
\item[{\rm (i)}] $(C^\infty(P), \{\cdot, \cdot \})$ is a Lie algebra;
\item[{\rm (ii)}] $\{\cdot, \cdot \}$ satisfies the Leibniz identity
on each factor;
\item[{\rm (iii)}] the vector bundle map
$\sharp: T^\ast P \rightarrow T^{\ast\ast} P$ covering the identity  satisfies
$\sharp (T^\ast P) \subset TP$.
\end{enumerate}
\end{definition}

\noindent Condition (iii) allows one to introduce for any  function $h \in
C^\infty(P)$  the {\bfi Hamiltonian vector field\/}  by 
\[
X_h[f] := \langle df, X_h \rangle = \{f, h\}
\]
where $f$ is an arbitrary smooth locally defined function on $P$.
\medskip

Given two Banach Poisson manifolds $(P_1, \{\,,\}_1)$ and 
$(P_2, \{\,,\}_2)$, a smooth map 
$\varphi: P_1 \rightarrow P_2$ is said to be {\bfi canonical\/}
or a {\bfi Poisson map\/} if 
\begin{equation}
\label{Poisson map definition}
\varphi^\ast \{f, g\}_2 = \{\varphi^\ast f, \varphi^\ast g\}_1
\end{equation}
for any two smooth locally defined functions $f$ and $g$ on $P_2$.
Condition (iii) in the previous definition implies, like in
the finite dimensional case, that \eqref{Poisson map definition} is equivalent 
to
\begin{equation}
\label{Poisson map condition}
X^2_f \circ  \varphi = T\varphi \circ X^1_{f \circ \varphi}
\end{equation}
for any smooth locally defined function $f$ on $P_2$ (for the proof see e.g. 
Marsden and Ratiu [1998], \S 10.3). Therefore, the  flow of a Hamiltonian 
vector field is a Poisson map and Hamilton's equations in Poisson bracket 
formulation are valid.
\medskip

For later applications we shall need the notion of the product of Banach 
Poisson manifolds. The definition we shall give is the one used in finite 
dimensions (see, e.g. Weinstein [1983] or Vaisman [1994]). However, the proof  of the
theorem characterizing the product needs some care due to the infinite 
dimensionality of the manifolds and the additional condition (iii) imposed in 
Definition \ref{Poisson manifold definition}. For this reason we shall sketch 
it below.

\begin{theorem}
\label{product theorem}
Given the Banach Poisson manifolds $(P_1, \{\,,\}_1)$ and 
$(P_2, \{\,,\}_2)$ there is a unique Banach Poisson structure 
$\{\,,\}_{12}$ on the product manifold $P_1 \times P_2$ such that:
\begin{enumerate}
\item[{\rm (i)}] the canonical projections 
$\pi_1: P_1 \times P_2 \rightarrow P_1$ 
and $\pi_2: P_1 \times P_2 \rightarrow P_2$ are Poisson maps;
\item[{\rm (ii)}] $\pi_1^\ast (C^\infty(P_1))$ and $\pi_2^\ast (C^\infty(P_2))$
are Poisson commuting subalgebras of $C^\infty(P_1 \times P_2)$.
\end{enumerate}
This unique Poisson structure on $P_1 \times P_2$ is called the 
{\bfi product\/} Poisson structure and its bracket is given by the formula
\begin{equation}
\label{product bracket}
\{f, g\}_{12}(p_1, p_2) = \{f_{p_2}, g_{p_2}\}_1(p_1) + 
\{f_{p_1}, g_{p_1}\}_2(p_2),
\end{equation}
where $f_{p_1}, g_{p_1}  \in C^\infty(P_2)$ and 
$f_{p_2}, g_{p_2} \in C^\infty(P_1)$ are the partial functions given by
     $f_{p_1}(p_2) = f_{p_2}(p_1) = f(p_1, p_2)$ and similarly for $g$. 
     \end{theorem}  

\noindent \textbf{Proof}. Recall that if $f \in C^\infty(P_1 \times P_2)$ 
then the partial exterior derivative $d_1 f(p_1, p_2)$ relative to 
$P_1$ is defined by $d_1 f(p_1, p_2) := df_{p_2} (p_1) 
= (\pi_1^\ast df_{p_2})(p_1, p_2)= d(\pi_1^\ast f_{p_2})(p_1, p_2)$ and 
similarly $d_2 f(p_1, p_2) = d(\pi_2^\ast f_{p_1})(p_1,p_2)$. Therefore,
$d f(p_1, p_2) = d_1 f(p_1, p_2) + d_2 f(p_1, p_2) 
= d(\pi_1^\ast f_{p_2})(p_1, p_2) + d(\pi_2^\ast f_{p_1})(p_1, p_2)$. Thus the 
functions $f$ and $\pi_1^\ast f_{p_2} + \pi_2^\ast f_{p_1}$ have the 
same derivatives at the point $(p_1, p_2) \in P_1 \times P_2$. Similarly,
$g$ and $\pi_1^\ast g_{p_2} + \pi_2^\ast g_{p_1}$ have the 
same derivatives at the point $(p_1, p_2) \in P_1 \times P_2$. 

Assume that there is a Poisson bracket $\{\,,\}_{12}$ on $P_1\times P_2$ 
satisfying the conditions in the theorem. Since any
Poisson bracket depends only on the first derivatives of the functions 
we necessarily have
\begin{align*}
\{f, g\}_{12}(p_1, p_2) 
&= \{\pi_1^\ast f_{p_2} + \pi_2^\ast f_{p_1}, 
\pi_1^\ast g_{p_2} + \pi_2^\ast g_{p_1} \}_{12}(p_1, p_2) \\
& = \{\pi_1^\ast f_{p_2}, \pi_1^\ast g_{p_2} \}_{12}(p_1, p_2)
+ \{\pi_1^\ast f_{p_2} , \pi_2^\ast g_{p_1} \}_{12}(p_1, p_2) \\
& \qquad \quad 
+\{ \pi_2^\ast f_{p_1}, \pi_1^\ast g_{p_2} \}_{12}(p_1, p_2)
+ \{\pi_2^\ast f_{p_1},  \pi_2^\ast g_{p_1} \}_{12}(p_1, p_2) \\
&= (\pi_1^\ast \{ f_{p_2} , g_{p_2} \}_1)(p_1, p_2)
+ (\pi_2^\ast \{ f_{p_1},  g_{p_1} \}_2)(p_1, p_2) \\
&= \{ f_{p_2} , g_{p_2} \}_1(p_1)
+ \{ f_{p_1},  g_{p_1} \}_2(p_2),
\end{align*}
where condition (ii) and (i) were used in the third equality. This shows that 
the Poisson bracket, if it exists, is 
unique and is given by \eqref{product bracket}.

Now define $\{\,,\}_{12}$ by \eqref{product bracket}. It remains to show
that the axioms in Definition \ref{Poisson manifold definition} hold. It is 
obvious that this operation satisfies the Leibniz identity, 
is bilinear, and skew symmetric. By Definition 
\ref{Poisson manifold definition} (iii), one can use Hamiltonian 
vector fields to express $\{\{f,g\},h\}$. A direct computation gives
\begin{align*}
&\{\{f,g\},h\}_{12}(p_1, p_2) 
= \{\{f_{p_2} ,g_{p_2}\}_1,h_{p_2}\}_1(p_1)
+\{\{f_{p_1} ,g_{p_1}\}_2,h_{p_1}\}_2(p_2)\\
&\qquad + d_1 d_2 f(p_1, p_2)\left( X_{h_{p_2}}^1(p_1),
X_{g_{p_1}}^2(p_2) \right) 
+ d_1 d_2 f(p_1, p_2)\left( X_{g_{p_2}}^1(p_1),
X_{h_{p_1}}^2(p_2)\right)\\
&\qquad - d_1 d_2 g(p_1, p_2)\left( X_{h_{p_2}}^1(p_1),
X_{f_{p_1}}^2(p_2) \right)
- d_1 d_2 g(p_1, p_2)\left( X_{f_{p_2}}^1(p_1),
X_{h_{p_1}}^2(p_2) \right),
\end{align*}
where $d_1d_2 f $ denotes the second mixed partial derivative 
of $f$ and where $X_{f_{p_2}}^1$ is the Hamiltonian vector 
field on $P_1$ corresponding to the function 
$f_{p_2} \in C^\infty (P_1)$ and similarly for the other ones.
Adding the other two terms obtained by circular permutation gives
zero since the  first two terms summed with their analogues vanish by the 
Jacobi identity on $P_1$ and $P_2$ respectively and the other terms cancel.

Since Hamiltonian vector fields on $P_1$ and $P_2$ exist by Definition \ref
{Poisson manifold definition}, formula 
\eqref{product bracket} shows that the Hamiltonian vector field on 
$P_1 \times P_2$ exists and is given by 
\begin{equation}
\label{product Hamiltonian vector field}
X^{12}_h (p_1,p_2) = \left(X^1_{h_{p_2}}(p_1), 
X^2_{h_{p_1}}(p_2)\right) \in T_{p_1} P_1 \times T_{p_2} P_2,
\end{equation}
where condition (iii) in Definition \ref{Poisson manifold definition} was used 
on $P_1$ and $P_2$; we have identified here $T_{(p_1, p_2)}(P_1 \times P_2)$ with
$T_{p_1} P_1 \times T_{p_2} P_2$.
Thus all conditions in Definition \ref{Poisson manifold definition} hold which
proves that $P_1 \times P_2$ is a Banach Poisson manifold. $\quad 
\blacksquare$
\medskip

We remark that \eqref{product bracket} implies that the product is functorial, 
that is, if $\varphi_1 :P_1 \rightarrow P_1'$ and $\varphi_2 :P_2 \rightarrow 
P_2'$ are Poisson maps then their product
$\varphi_1 \times \varphi_2 :P_1 \times P_2 \rightarrow 
P_1' \times P_2'$ is also a Poisson map.
\medskip

Returning to Definition \ref{Poisson manifold definition}, it should be noted 
that the condition $\sharp(T^\ast P) \subset TP$ is
automatically satisfied in certain cases: 
\begin{itemize}
\item if $P$ is a smooth manifold
modeled on a reflexive Banach space, that is $\mathfrak{b}^{\ast \ast} = 
\mathfrak{b}$, or
\item if $P$ is a smooth manifold
modeled on a Banach space whose norm is $C^\infty$ away from the origin, or
\item $P$ is a strong symplectic manifold with symplectic form $\omega$.
\end{itemize}
The second condition holds because the assignment $f \mapsto \{\cdot, f\}$ is a
derivation on the space of smooth functions on $P$ which guarantees that it is defined
by a vector field; see Abraham, Marsden, and Ratiu [1988], \S4.2, for the proof. In
particular, the first two conditions hold if $P$ is a Hilbert (and, in particular, a
finite dimensional) manifold.

Any strong symplectic manifold $(P, \omega)$ is a Poisson manifold
in the sense of Definition \ref{Poisson manifold definition}.
Recall that {\bfi strong\/} means that for each $p\in P$ the map
\begin{equation}
\label{symplectic duality}
v_p \in T_p P \mapsto \omega(p)(v_p, \cdot) \in T_p^\ast P
\end{equation}
is a bijective continuous linear map. Therefore, given a smooth function $f:P
\rightarrow \mathbb{R}$ there exists a vector field $X_f$ such
that $df = \omega(X_f, \cdot)$. The Poisson bracket is defined by
$\{f, g \} = \omega(X_f, X_g) = \langle df, X_g \rangle$, thus $\sharp df = 
X_f$, so  $\sharp( T^\ast P) \subset TP.$

On the other hand, a weak symplectic manifold is not a Poisson
manifold in the sense of Definition \ref{Poisson manifold
definition}. Recall that {\bfi weak} means that the map defined by
\eqref{symplectic duality} is an injective continuous linear map
that is, in general, not surjective. Therefore, one cannot
construct the map that associates to every differential $df$ of a
smooth function $f: P \rightarrow \mathbb{R}$ the Hamiltonian
vector field $X_f$. Since the definition of the Poisson bracket
should be $\{f, g \} = \omega(X_f, X_g)$, one cannot define this
operation on functions and hence weak symplectic manifold
structures do not define, in general, Poisson manifold structures in the
sense of Definition \ref{Poisson 
manifold definition}. There are various ways to deal
with this problem. One of them is to restrict the space of
functions on which one is working, as is often done in field theory. Another is
to deal with densely defined vector
fields and invoke the theory of (nonlinear) semigroups; see
Chernoff and Marsden [1974] for this approach. A simple example illustrating the 
importance  of the underlying topology is given by the canonical symplectic 
structure on $\mathfrak{b} \times \mathfrak{b}^\ast$, where $\mathfrak{b}$ is
a Banach space. This canonical symplectic structure is in general weak; if 
$\mathfrak{b}$ is reflexive then it is strong.

In this paper we shall not address these important questions regarding 
weak symplectic manifolds and their relation to Poisson structures 
and we shall exclusively consider Banach Poisson manifolds as given by Definition
\ref{Poisson manifold definition}. Thus, in some sense, the
Poisson manifolds considered in this paper are generalizations of
strong symplectic manifolds. However, in \S \ref{section: symplectic leaves}
and \S\ref{section: momentum maps and reduction} we
shall give examples illustrating the need for a weakening of 
Definition \ref{Poisson manifold definition}.

We shall need in the sequel various notions of Poisson structures defined 
on almost complex and complex manifolds. We briefly summarize the various
possibilities below. 
\medskip

Assume that the real Banach manifold $P$ underlying the Poisson structure
given by the tensor field $\varpi$ has also the structure of an almost
complex manifold, that is, there is a smooth vector bundle map $I:TP
\rightarrow TP$ covering the identity which satisfies $I^2 = - id$. The
question then arises what does it mean for the Poisson and almost complex
structures to be compatible. The Poisson structure $\varpi$ is said to be
{\bfi compatible with the almost complex structure\/} $I$ if the following
diagram commutes:

\unitlength=5mm
\begin{center}
\begin{picture}(9,8)
\put(1,7){\makebox(0,0){$T^\ast P$}} 
\put(9,7){\makebox(0,0){$TP$}}
\put(1,2){\makebox(0,0){$T^\ast P$}}
\put(9,2){\makebox(0,0){$TP$}}
\put(1,3){\vector(0,1){3}} 
\put(9,6){\vector(0,-1){3}}
\put(2,7){\vector(1,0){5.5}} 
\put(2,2){\vector(1,0){5.7}}
\put(0.3,4.3){\makebox(0,0){$I^\ast$}}
\put(9.4,4.3){\makebox(0,0){$I$}} 
\put(5,7.5){\makebox(0,0){$\sharp$}}
\put(5,2.5){\makebox(0,0){$\sharp$}}
\end{picture}
\end{center}
that is, 
\begin{equation}
\label{almost complex compatibility}
I \circ \sharp + \sharp \circ I^\ast = 0.
\end{equation}

The decomposition 
\begin{equation}
\label{Poisson decomposition}
\varpi = \varpi_{(2,0)} + \varpi_{(1,1)} + \varpi_{(0,2)}
\end{equation}
induced by the almost complex structure $I$ and the reality of $\varpi$, 
implies that the compatibility condition \eqref{almost complex
compatibility} is equivalent to
\begin{equation}
\label{almost complex compatibility tensor} 
\varpi_{(1,1)} = 0 \qquad \text{and} \qquad 
\overline{\varpi}_{(2,0)} 
= \varpi_{(0,2)}.
\end{equation}  
In view of \eqref{almost complex compatibility tensor}, 
$[\varpi, \varpi]_S = 0$ is equivalent to 
\begin{equation}
\label{complex Poisson}
[\varpi_{(2,0)}, \varpi_{(2,0)}]_S = 0 \qquad \text{and} \qquad [\varpi_
{(2,0)}, \overline{\varpi}_{(2,0)}]_S = 0.
\end{equation}

If \eqref{almost complex compatibility} holds, the triple 
$(P, \{\cdot, \cdot\}, I)$ is called an {\bfi almost complex Banach
Poisson manifold\/}. If $I$ is given by a complex analytic structure
$P_\mathbb{C}$ on $P$ it will be called a {\bfi complex Banach
Poisson manifold\/}. For finite dimensional complex manifolds these
structures were introduced and studied by Lichnerowicz [1988].

Denote by $\mathcal{O}\Omega^{(k,0)}(P_\mathbb{C})$ and
$\mathcal{O}\Omega_{(k,0)}(P_\mathbb{C})$ the space of holomorphic
$k$-forms and $k$-vector fields respectively. If
\begin{equation}
\label{complex compatibility}
\sharp\left(\mathcal{O}\Omega^{(1,0)}(P_\mathbb{C})\right) \subset 
\mathcal{O}\Omega_{(1,0)}(P_\mathbb{C}),
\end{equation}
that is, the Hamiltonian vector field $X_f$ is holomorphic if $f$ is a 
holomorphic function, then, in addition to \eqref{almost complex
compatibility tensor} and \eqref{complex Poisson}, one has $\varpi_{(2,0)} \in
\mathcal{O}\Omega_{(2,0)}(P_\mathbb{C})$. As expected, the compatibility
condition \eqref{complex compatibility} is stronger than \eqref{almost
complex compatibility}. Note that \eqref{complex compatibility} implies
the second condition in \eqref{complex Poisson}. Thus the compatibility
condition \eqref{complex compatibility} induces on the underlying complex
Banach manifold $P_\mathbb{C}$ a holomorphic Poisson tensor
$\varpi_\mathbb{C} := \varpi_{(2,0)}$. A pair $(P_\mathbb{C},
\varpi_\mathbb{C})$ consisting of an analytic complex manifold
$P_{\mathbb{C}}$ and a holomorphic skew symmetric contravariant two-tensor
field $\varpi_\mathbb{C}$ such that $[\varpi_\mathbb{C},
\varpi_\mathbb{C}]_S=0$ and \eqref{complex compatibility} holds will be
called a {\bfi holomorphic Banach Poisson manifold\/}.
\medskip

Consider now a holomorphic Poisson manifold $(P, \varpi)$. Denote by
$P_{\mathbb{R}}$ the underlying real Banach manifold and define the real two-vector
field $\varpi_\mathbb{R} :=\operatorname{Re} \varpi$. It is easy to see that
$(P_\mathbb{R},\varpi_\mathbb{R})$ is a real Poisson manifold compatible with the
complex Banach manifold structure of $P$ and $(\varpi_\mathbb{R})_\mathbb{C} =
\varpi$. Summarizing, we have shown that there are two procedures that are
inverses of each other: a holomorphic Poisson manifold corresponds in a
bijective manner to a real Poisson manifold whose Poisson tensor is
compatible with the underlying complex manifold structure.  One can call
these constructions the {\bfi complexification} and {\bfi realification}
of Poisson structures on complex manifolds.

\section{Classical reduction}
\label{section: classical reduction}

We shall review in this section the theory of classical Poisson
reduction for Banach Poisson manifolds. Let $(P, \{\cdot,\cdot\}_P)$ be a
real Banach Poisson manifold (in the sense of Definition \ref{Poisson manifold
definition}), $i: N \hookrightarrow P$ be a (locally
closed) submanifold, and $E \subset (TP)|_N$ be a subbundle of the
tangent bundle of $P$ restricted to $N$. For simplicity we make
the following topological regularity assumption throughout this
section: $E \cap TN$ is the tangent bundle to a foliation
$\mathcal{F}$ whose leaves are the fibers of a submersion $\pi : N
\rightarrow M: = N/ \mathcal{F}$, that is, one assumes that the
quotient topological space $N/ \mathcal{F}$ admits the quotient
manifold structure. The subbundle $E$ is said to be {\bfi
compatible with the Poisson structure\/} provided the following
condition holds: if $U \subset P$ is any open subset and $f, g \in
C^\infty (U)$ are two arbitrary functions whose differentials $df$
and $dg$ vanish on $E$, then $d\{f, g\}_P$ also vanishes on $E$. The
triple $(P, N, E)$ is said to be {\bfi reducible\/}, if $E$ is
compatible with the Poisson structure on $P$ and the manifold $M:
= N/\mathcal{F}$ carries a Poisson bracket $\{\cdot, \cdot\}_M$ (in the
sense of Definition \ref{Poisson manifold definition}) such that
for any smooth local functions $\bar{f}, \bar{g}$ on $M$ and any
smooth local extensions $f, g$ of $\bar{f} \circ \pi$,
$\bar{g}\circ \pi$ respectively, satisfying $df|_E = 0$, $dg|_E =
0$, the following relation on the common domain of definition of
$f$ and $g$ holds:
\begin{equation}
\label{reduction condition}
\{f, g\}_P \circ i = \{\bar{f},
\bar{g}\}_M \circ \pi.
\end{equation}
If $(P,N,E)$ is a reducible triple then $(M = N/\mathcal{F},
\{\cdot,\cdot \}_M)$ is called the {\bfi reduced manifold\/} of $P$ via
$(N, E)$. Note that \eqref{reduction condition} guarantees that if the reduced
Poisson bracket $\{\cdot, \cdot \}_M$ on $M$ exists, it is necessarily unique.

Given a subbundle $E \subset TP$, its {\bfi annihilator\/} is
defined as the subbundle of $T^\ast P$ given by $E^\circ : =
\{\alpha \in T^\ast P \mid \langle \alpha, v \rangle = 0 \text{~ for ~all~ } v \in
E\}$.

The following statement generalizing the finite dimensional
Poisson reduction theorem of Marsden and Ratiu [1986] is central
for our purposes. The proof in infinite dimensions is a
modification of the original one (see the above mentioned paper or
Vaisman [1994], \S 7.2, for the finite dimensional proof).

\begin{theorem}
\label{reduction theorem} 
Let $P$, $N$, $E$ be as above and assume that
$E$ is compatible with the Poisson structure on $P$. The triple
$(P,N,E)$ is reducible if and only if $\sharp(E^{\circ}_n) \subset
\overline{T_nN + E_n}$ for every $n \in N$.
\end{theorem}

\noindent\textbf{Proof}. Assume that $(P,N,E)$ is reducible. Thus
$M: = N/\mathcal{F}$ is a Banach Poisson manifold and (\ref{reduction condition})
holds. In addition, recall that $N$ is
a (locally closed) submanifold of $P$ and that $E \cap TN$ is the
tangent bundle of a foliation on $N$. For $n \in N$, choose a
chart domain $U$ of $n$ in $P$ with the submanifold property
relative to $N$ and such that $U \cap N$ is foliated.

Given $\alpha_n \in E^\circ _n$, find a smooth function $f$ on $U$
(shrunk if necessary), such that $df(n) = \alpha_n$ and $\langle
df, E \rangle = 0$. This is possible since $E$ is a subbundle of
$TP|_N$ and $E \cap TN$ is the tangent bundle to a foliation on
$N$. Let $\bar{f}$ be the smooth function on $\pi(U \cap N)
\subset M$ induced by $f$, that is, $f|_N = \bar{f} \circ \pi$.
Therefore, $f: U \rightarrow \mathbb{R}$ is a local extension of
$\bar{f} \circ \pi$.

Next, take an arbitrary $\beta_n \in (E_n + T_nN)^\circ =
E_n^\circ \cap (T_nN)^\circ$ and find a smooth function $g$ on $U$
such that $N\cap U = g^{-1}(0)$, $\langle dg, E \rangle = 0$, and
$dg(n) = \beta_n$. Again, the existence of $g$ is insured by the
hypothesis that $E$ is a subbundle of $TP|_N$ and that $E \cap TN$ is
the tangent bundle to a foliation on $N$. Thus $g: U \rightarrow
\mathbb{R}$ is a local extension of $0 \circ \pi$, where $0$ is
the identically zero function on $M$. Then we have by
(\ref{reduction condition})
\[
\langle \beta_n, \sharp(\alpha_n) \rangle = \langle dg(n), X_f(n)
\rangle = \{g, f \}_P(n) = \{0, f\}_M(\pi(n)) = 0.
\]
This shows that $\sharp(\alpha_n) \in (E_n + T_nN)^{\circ \circ} =
\overline{E_n + T_nN}$, that is, $\sharp(E^{\circ}_n) \subset
\overline{T_nN + E_n}$ for every $n \in N$.

Conversely, assume that $\sharp(E^{\circ}_n) \subset \overline{T_nN +
E_n}$ for every $n \in N$. For $\bar{f}, \bar{g}$ locally defined
smooth functions on $M$ we need to define their Poisson bracket,
show that (\ref{reduction condition}) holds, and that all
conditions in Definition \ref{Poisson manifold definition} are
satisfied. Let $f, g$ be local extensions of $\bar{f}\circ \pi$
and $\bar{g} \circ \pi$ respectively, such that $df$ and $dg$
vanish on $E$. Since $E$ is compatible with the Poisson bracket on
$P$, $d\{f, g\}_P$ also vanishes on $E$ and thus $\{f, g\}_P$ is
constant on the leaves of $\mathcal{F}$ thereby inducing a smooth
locally defined function on $M$. We take this function to be the
definition of $\{\bar{f}, \bar{g}\}_M$. If we show that this
function is well defined, that is, is independent on the
extensions chosen, then the axioms of a Poisson bracket (that is,
conditions (i) and (ii) in Definition \ref{Poisson manifold
definition}) are trivially verified and, by construction,
(\ref{reduction condition}) holds.

Since the Poisson bracket is skew symmetric it suffices to show
the independence of the extension only for the function $f$. So
let $f'$ be another local extension of $\bar{f} \circ \pi$ such
that $\langle df', E \rangle = 0$. On the common domain of
definition of $f$ and $f'$, we have hence $(f - f')|_N = 0$; in
particular, $d(f - f')$ vanishes on $TN$. However, since both $df$
and $df'$ vanish on $E$, it follows that $d(f - f')$ vanishes on
$E + TN$. Let $n \in N$ be an arbitrary point in the common domain
of definition of $f$ and $f'$. By continuity, $d(f - f')(n)$
vanishes on $\overline{E_n + T_nN}$. Since $X_g(n) \in \sharp(E^\circ
_n)$, using the working hypothesis $\sharp(E^{\circ}_n) \subset
\overline{T_nN + E_n}$, we conclude
\[
\{f - f', g\}_P(n) = \langle d(f- f')(n), X_g(n) \rangle = 0,
\]
that is, $\{f, g\}_P(n) = \{f', g\}_P(n)$.

It remains to verify condition (iii) of Definition \ref{Poisson
manifold definition}, that is, $\bar{\sharp}(T^\ast M) \subset TM$,
where $\bar{\sharp}: T^\ast M \rightarrow T^{\ast \ast}M$ is the
vector bundle map covering the identity defined by
$\bar{\sharp}_m(d\bar{f}(m)) := \{\cdot , \bar{f}\}_M(m)$ for any smooth
locally defined function $\bar{f}$ on $M$. The idea of the proof below is to
use \eqref{reduction condition} to show that $\bar{\sharp}_m(d\bar{f}(m)) =
T_n\pi\left(X_f(n)\right) \in T_mM$ for every $m \in M$ and every locally
defined function $\bar{f}$ around $m$.

To do this, let $\bar{f}, \bar{g}: W \rightarrow \mathbb{R}$ be two arbitrary
smooth functions, where $W$ is a chart domain on $M$ containing the point $m$.
We shall construct now local extensions of $\bar{f} \circ \pi$ and $\bar{g} \circ
\pi$ adapted to our needs. Since we have already shown that that the definition
of $\{\cdot , \cdot \}_M$ is independent on the extensions, we can work only with
these extensions and conclude the desired result.  Since
$E\cap TN$ is the tangent bundle to a foliation on $N$, if $n
\in N$ is such that $\pi (n) = m$, there is a foliated chart on
$N$ around $n$ whose domain is of the form $W \times W'$ (after an eventual
shrinking of $W$), that is, the leaves of the foliation are given by $\{w\} \times W'$
for all $w \in W$. Since $N$ is a submanifold of $P$ and since $E$ is
defined only along $N$, there is a chart on $P$ whose domain is of
the form $W \times W' \times V'$ (after shrinking, if necessary, both $W$ and
$W'$). Define the local extension $f: W \times W' \times V'
\rightarrow \mathbb{R}$ of $\bar{f}\circ \pi$ by $f(w, w', v') := \bar{f}(w)$.  By
condition (iii) of Definition
\ref{Poisson manifold definition}, $\sharp_n(df)(n)$ is a vector of the form
$X_f(n) \in T_nP$. Let us show that $X_f(n)$ is tangent to $N$. This is
equivalent to proving that for any linear continuous functional $\beta_n$ on the
ambient Banach space containing $V'$, we have $\langle \beta_n, X_f(n) \rangle = 0$.
However, $\beta_n = dk(n)$, for some smooth function $k: W \times W' \times V'
\rightarrow \mathbb{R}$ that does not depend on the variables from $W$ and $W'$.
But then
$k$ is a local extension of $0 \circ \pi$ and, using \eqref{reduction condition}, we
get
\[
\langle \beta_n, X_f(n) \rangle = \langle dk(n), X_f(n) \rangle = \{k, f\}_P(n) =
\{0, f\}_M(\pi(n)) = 0,
\]
which proves the claim.

Construct in the same fashion a local extension
of $\bar{g}\circ \pi$ to the same open neighborhood of $n$ in $P$.
Since $dg(n) \circ T_n i = d \bar{g}(m) \circ T_n \pi$, $m = \pi(n)$, we have by
\eqref{reduction condition},
\begin{align*}
\bar{\sharp}_m(d\bar{f}(m)) (d \bar{g}(m)) &= \{\bar{g},
\bar{f}\}_M(m) = \{g, f\}_P(n) = \sharp_n(df(n))(dg(n)) \\ &=
\left\langle dg(n), X_f(n) \right\rangle = \left\langle
d\bar{g}(m), T_n\pi\left(X_f(n)\right)\right\rangle.
\end{align*}
Since $\bar{g}$ is an arbitrary smooth function defined on a
neighborhood of $m$, the Hahn-Banach Theorem and the inclusion of
the Banach space into its bidual imply that
$\bar{\sharp}_m(d\bar{f}(m)) = T_n\pi\left(X_f(n)\right) \in T_mM$ for
every $m \in M$ , that is $\bar{\sharp}(T^\ast M) \subset TM$. \quad
$\blacksquare$
\medskip

The behavior of Poisson maps and Hamiltonian dynamics under
reduction is given by the following two theorems whose proofs are
identical to the ones in finite dimensions (Marsden and Ratiu
[1986] or Vaisman [1994], \S 7.4).

\begin{theorem}
\label{reduction of maps}
Let $(P_1, N_1, E_1)$ and $(P_2, N_2, E_2)$ be Poisson reducible triples and assume
that $\varphi: P_1 \rightarrow P_2$ is a Poisson map satisfying $\varphi(N_1) \subset
N_2$ and $T\varphi (E_1) \subset E_2$. Let $\mathcal{F}_i$ be the
regular foliation on $N_i$ defined by the subbundle $E_i$ and
denote by $\pi_i: N_i \rightarrow M_i:= N_i/\mathcal{F}_i$,
$i=1,2$, the reduced Poisson manifolds. Then there is a unique
induced Poisson map $\overline{\varphi}: M_1 \rightarrow M_2$,
called the {\bfi reduction\/} of $\varphi$, such that $\pi_2 \circ
\varphi = \overline{\varphi} \circ \pi_1$.
\end{theorem}

\noindent\textbf{Proof}. The hypotheses imply that $T\varphi(E_1
\cap TN_1) \subset E_2 \cap TN_2$ and hence $\varphi$ maps the
leaves of the foliation $\mathcal{F}_1$ to those of
$\mathcal{F}_2$. Therefore $\varphi$ is a projectable map, that
is, there exists a smooth map $\overline{\varphi}: M_1 \rightarrow M_2$
such that $\pi_2 \circ \varphi = \overline{\varphi} \circ \pi_1$.
It remains to be shown that $\overline{\varphi}$ is a Poisson map.

Let $\bar{f}$ and $\bar{g}$ be two smooth local functions on $M_2$
and let $f$ and $g$ be the local extensions of $\bar{f}\circ
\pi_2$ and $\bar{g}\circ \pi_2$ respectively, such that $df|_{E_2}
= dg|_{E_2} = 0$. Since $T\varphi(E_1) \subset E_2$ it follows
that $d(f \circ \varphi)|_{E_1} = d(g \circ \varphi)|_{E_1} = 0$.
Hence $f \circ \varphi$ is a smooth local extension of $\bar{f}
\circ \pi_2 \circ \varphi = \bar{f} \circ \bar{\varphi} \circ
\pi_1$. Similarly, $g \circ \varphi$ is a smooth local extension
of $\bar{g} \circ \pi_2 \circ \varphi = \bar{g} \circ
\bar{\varphi} \circ \pi_1$. Definition \eqref{reduction
condition} gives then
\begin{align*}
\{\bar{f} \circ \bar{\varphi}, \bar{g} \circ \bar{\varphi}\}_{M_1} \circ \pi_1 
&= \{ f \circ \varphi, g \circ
\varphi\}_{P_1} \circ i_1 = \{f , g \}_{P_2} \circ \varphi \circ i_1 \\
& = \{ f , g \}_{P_2} \circ i_2 \circ \varphi = \{ \bar{f} ,
\bar{g} \}_{M_2} \circ \pi_2 \circ \varphi = \{ \bar{f} , \bar{g}
\}_{M_2} \circ \bar{\varphi} \circ \pi_1,
\end{align*}
which implies that $\bar{\varphi}$ is a Poisson map by
surjectivity of $\pi_1$. \quad $\blacksquare$
\medskip

\begin{theorem}
\label{reduction of dynamics} Let $(P, N, E)$ be a Poisson
reducible triple and $\pi: N \rightarrow M$ be the corresponding
Poisson reduced manifold. Assume that $h \in C^\infty(P)$ and the
associated flow $\varphi_t$ of the Hamiltonian vector field $X_h$
satisfies the conditions
\begin{itemize}
\item[{\rm (i)}] $dh|_E = 0$,
\item[{\rm (ii)}] $\varphi_t(N) \subset N$,
\item[{\rm (iii)}] $T\varphi_t(E) \subset E$
\end{itemize}
for all $t$ for which the flow $\varphi_t$ is defined. Then
the reduction $\overline{\varphi}_t$ is the flow of the
Hamiltonian vector field on $M$ given by the function
$\bar{h}$ uniquely determined by the condition $\bar{h}
\circ \pi = h|_N$. The Hamiltonian vector fields $X_h$ on $N
\subset P$ and $X_{\bar{h}}$ on $M$ are $\pi$-related.
\end{theorem}

\noindent\textbf{Proof}. The hypotheses guarantee by Theorem
\ref{reduction of maps} that the flow $\varphi_t$ of $X_h$ reduces
to a smooth flow $\overline{\varphi}_t$ on the reduced manifold $M$. The
hypothesis on $h$ insures the existence of the smooth function
$\bar{h}$ on $M$. Let us prove that $X_h$ and $X_{\bar{h}}$ are
$\pi$-related. If this is done, their flows are necessarily $\pi$-related 
and hence, by surjectivity of $\pi$, it follows that the flow of 
$X_{\bar{h}}$ is $\overline{\varphi}_t$.

To prove that $X_h$ and $X_{\bar{h}}$ are
$\pi$-related, let $\bar{f}$ be a smooth locally defined function
in a neighborhood of $m \in M$ and let $f$ be a smooth local
extension of $\bar{f} \circ \pi$ satisfying $df|_E = 0$. Then,
since $X_h$ is tangent to $N$, using the defining identity of the
reduced bracket \eqref{reduction condition}, for any $n \in N$ satisfying $\pi(n) =
m$, we get
\begin{align*}
\langle d\bar{f} (m), X_{\bar{h}}(m) \rangle &= \{\bar{f},
\bar{h}\}_M(m) = \{f, h \}_P(n) = \langle df (n), X_{h}(n)
\rangle\\ & = \langle d(\bar{f} \circ \pi) (n), X_{h}(n) \rangle =
\langle d\bar{f} (m), T_n \pi\left( X_{h}(n) \right) \rangle,
\end{align*}
which proves that $X_{\bar{h}} \circ \pi = T\pi \circ X_h$.  \quad $\blacksquare$
\medskip

In this paper we shall not investigate the consistency of Poisson
reduction with other structures such as almost complex, complex, and
holomorphic structures. For finite dimensional Poisson manifolds the
consistency of Poisson reduction with the complex structure was presented
in Nunes da Costa [1997].

\section{Banach Lie-Poisson spaces}
\label{section: banach lie-poisson spaces}

It is well known that the dual of any Lie algebra admits a linear Poisson 
structure, called the Lie-Poisson structure.  In this section 
we shall extend the definition of this structure to the infinite dimensional 
case in agreement with Definition \ref{Poisson manifold definition}. We shall 
call such spaces Banach Lie-Poisson spaces and shall investigate their 
properties. 
\medskip

Recall that a {\bfi Banach Lie algebra} $(\mathfrak{g}, [\cdot,
\cdot])$ is a Banach space that is also a Lie algebra such that
the Lie bracket is a bilinear continuous map $\mathfrak{g} \times
\mathfrak{g} \rightarrow \mathfrak{g}$. Thus the adjoint and
coadjoint maps $\operatorname{ad}_x :\mathfrak{g} \rightarrow
\mathfrak{g}$, $\operatorname{ad}_x y:= [x,y]$, and
$\operatorname{ad}_x^\ast : \mathfrak{g}^\ast \rightarrow
\mathfrak{g}^\ast$ are also continuous for each $x \in
\mathfrak{g}$.

\begin{definition}
\label{definition of Lie-Poisson space}
 A {\bfi Banach Lie-Poisson
space} $(\mathfrak{b},\{\cdot, \cdot\})$ is a real or holomorphic Poisson
manifold such that $\mathfrak{b}$ is a Banach space and the dual
$\mathfrak{b}^\ast \subset C^\infty(\mathfrak{b})$ is a Banach Lie
algebra under the Poisson bracket operation.
\end{definition}

Throughout this section we shall treat the real and the holomorphic
cases simultaneously.
Denote by $[\cdot, \cdot]$ the restriction of the Poisson bracket
$\{\cdot, \cdot \}$ from $C^\infty(\mathfrak{b})$ to the Lie
subalgebra $\mathfrak{b}^\ast$. For any $x, y \in \mathfrak{b}^\ast$ and $b\in 
\mathfrak{b}$ we have
\begin{align*}
\langle y, \operatorname{ad}_x ^\ast b \rangle &=
\langle [x,y], b \rangle = \{x, y\}(b) = -\{y, x\}(b) \\
&=-X_x[y](b) = -\langle D y (b), X_x(b) \rangle =
-\langle y, X_x(b) \rangle,
\end{align*}
where we have used the linearity of $y \in C^\infty(\mathfrak{b})$
to conclude that the Fr\'echet derivative $Dy(b) = y$. Thus we
obtain the following identity in the bidual $\mathfrak{b}^{\ast \ast}$:

\begin{equation}
\label{identity in bidual}
 X_x (b) = - \operatorname{ad}_x ^\ast b \qquad
\text{for} \qquad x \in
\mathfrak{b}^\ast, \quad b \in \mathfrak{b}.
\end{equation}

\begin{theorem}
\label{general theorem}
The Banach space $\mathfrak{b}$ is a
Banach Lie-Poisson space $(\mathfrak{b}, \{\cdot, \cdot \})$ if
and only if its dual $\mathfrak{b}^\ast$ is a Banach Lie algebra
$(\mathfrak{b}^\ast, [\cdot, \cdot ])$ satisfying
$\operatorname{ad}_x^\ast \mathfrak{b} \subset \mathfrak{b}
\subset \mathfrak{b}^{\ast \ast}$ for all $x \in
\mathfrak{b}^\ast$. Moreover, the Poisson bracket of $f, g \in
C^\infty(\mathfrak{b})$ is given by
\begin{equation}\label{general LP}
\{f, g\}(b) = \langle [Df(b), Dg(b)], b \rangle,
\end{equation}
where $b \in \mathfrak{b}$ and $D$ denotes the Fr\'echet
derivative. If $h$ is a smooth function on $\mathfrak{b}$, the
associated Hamiltonian vector field is given by
\begin{equation}
\label{general Hamiltonian vector field}
X_h(b) = - \operatorname{ad}^\ast _{Dh(b)} b.
\end{equation}
\end{theorem}

\noindent \textbf{Proof}. Assume that $\mathfrak{b}$ is a Banach
Lie-Poisson space relative to the bracket $\{\cdot, \cdot\}$. By
Definition \ref{definition of Lie-Poisson space}, its dual
$\mathfrak{b}^\ast$ is a Banach Lie algebra relative to the
bracket $[\cdot, \cdot ] : = \{\cdot, \cdot
\}|_{\mathfrak{b}^\ast}$. However, $\mathfrak{b}$ is also a
Poisson manifold and thus, by definition,  $X_x (b) \in
\mathfrak{b}$ for all $x\in \mathfrak{b}^\ast$ and all $b \in
\mathfrak{b}$. Formula (\ref{identity in bidual}) implies then
that $\operatorname{ad}^\ast _x (b) \in \mathfrak{b}$ for all
$x\in \mathfrak{b}^\ast$ and all $b \in \mathfrak{b}$ which is the
required condition.

Conversely, assume that $(\mathfrak{b}^\ast, [\cdot, \cdot ])$ is a Banach Lie
algebra  satisfying
$\operatorname{ad}_x^\ast \mathfrak{b} \subset \mathfrak{b}
\subset \mathfrak{b}^{\ast \ast}$ for all $x \in
\mathfrak{b}^\ast$. Define the bracket $\{f, g\}$ of $f, g \in
C^\infty(b)$ by (\ref{general LP}). All properties of the Poisson
bracket are trivially satisfied by \eqref{general LP} except for
the Jacobi identity. For this, we note that from
$\ad_x^*\b\subset\b$, $x\in\b^*$, one has
\begin{equation}
\label{eq:2.3} D\{f, g\}(b) = \langle  [Df(b), Dg(b)], \cdot
\rangle - D^2 f (b) \left(\operatorname{ad}^\ast_{Dg(b)} b, \cdot
\right) +D^2 g (b)\left( \operatorname{ad}^\ast_{Df(b)} b, \cdot
\right).
\end{equation}
for $f, g \in C^\infty(\b)$. Using \eqref{eq:2.3} we obtain
\begin{equation*}
\begin{split}
\{\{f, g\}, h\}(b)&= \left \langle  \left[D\{f, g\}(b), Dh(b)
\right], b \right \rangle \\ & = \left \langle \left[ \left[
Df(b), Dg(b)\right], Dh(b) \right], b \right\rangle +D^2
f(b)\left(\operatorname{ad}^\ast_{Dg(b)} b,
\operatorname{ad}^\ast_{Dh(b)} b \right)\\ &  \qquad \qquad \qquad
- D^2 g(b)\left(\operatorname{ad}^\ast_{Df(b)} b,
\operatorname{ad}^\ast_{Dh(b)} b \right).
\end{split}
\end{equation*}
Taking the two other terms obtained by circular permutation of
$f$, $g$, and $h$, using the Jacobi identity for the commutator
bracket in the sum of the first three terms and the symmetry of
the second derivative in the sum of the remaining terms, proves
that \eqref{general LP} satisfies the Jacobi identity.

Since
\begin{equation*}
\langle Df(b), X_h(b) \rangle = \{f, h\}(b) = \langle [Df(b),
Dh(b)], b \rangle =  -\left\langle Df(b), \operatorname{ad}^\ast
_{Dh(b)} b \right\rangle
\end{equation*}
for every $f \in C^\infty(\mathfrak{b})$ and
$\operatorname{ad}_x ^\ast \mathfrak{b} \subset
\mathfrak{b}$ for every $x \in \mathfrak{b}^\ast$, it follows
that the Hamiltonian vector field $X_h$ is given by
(\ref{general Hamiltonian vector field}). $\quad \blacksquare$
\bigskip

\begin{example}
\label{reflexive}
\normalfont
Let $\mathfrak{b}$ be a reflexive Banach Lie algebra, that is, 
$\mathfrak{b}^{\ast \ast} = \mathfrak{b}$. Then its dual $\mathfrak{b}^\ast$ 
is a Banach Lie-Poisson space. To see this, note that
$\mathfrak{b}^{\ast \ast} = \mathfrak{b}$ is a
Banach Lie algebra and that $\operatorname{ad}_x^\ast(\mathfrak{b}^\ast) 
\subset \mathfrak{b}^\ast$ for all $x\in \mathfrak{b}$, so 
Theorem \ref{general theorem} applies. \quad $\blacklozenge$
\end{example}

\begin{example}
\label{finite dimensional}
\normalfont
Since every finite dimensional Lie algebra is reflexive Example 
\ref{reflexive} yields the following classical result: the dual of 
any finite dimensional Lie algebra is a Lie-Poisson space. \quad
$\blacklozenge$
\end{example}

\begin{definition}
\label{Lie-Poisson morphism}
A {\bfi morphism\/} between two
Banach Lie-Poisson spaces $\mathfrak{b}_1$ and $\mathfrak{b}_2$ is
a continuous linear map $\phi:\mathfrak{b}_1 \rightarrow
\mathfrak{b}_2$ that preserves the Poisson bracket structure, that is,
\[
\{f \circ \phi, g\circ \phi \}_1 = \{f, g \}_2 \circ \phi
\]
for any $f, g \in C^\infty(\mathfrak{b}_2)$. Such a map $\phi$
is  also called a {\bfi linear Poisson map}.
\end{definition}

We consider now the category $\mathfrak{B}$ whose objects are the
Banach Lie-Poisson spaces and whose morphisms are the linear
Poisson maps.

Let $\mathfrak{L}$ denote the category of Banach Lie algebras and
continuous Lie algebra homomorphisms. Denote by $\mathfrak{L}_0$
the following subcategory of $\mathfrak{L}$. An object of
$\mathfrak{L}_0$ is a Banach Lie algebra $\mathfrak{g}$ admitting
a predual $\mathfrak{g}_\ast$, that is, $(\mathfrak{g}_\ast)^\ast
= \mathfrak{g}$, and satisfying $\operatorname{ad}^\ast_{\mathfrak{g}}
\mathfrak{g}_\ast \subset \mathfrak{g}_\ast$ where $\operatorname{ad}^\ast$ is the
coadjoint  representation of $\mathfrak{g}$ on $\mathfrak{g}^\ast$; note that
$\mathfrak{g}_\ast \subset \mathfrak{g}^\ast$. A morphism in the
category $\mathfrak{L}_0$ is a Banach Lie algebra homomorphism
$\psi: \mathfrak{g}_1 \rightarrow \mathfrak{g}_2$ such that the
dual map $\psi^\ast: \mathfrak{g}_2^\ast \rightarrow
\mathfrak{g}_1^\ast$ preserves the corresponding preduals, that
is, $\psi^\ast: (\mathfrak{g}_2)_\ast \rightarrow
(\mathfrak{g}_1)_\ast$.

\begin{theorem}
\label{functor}
The category $\mathfrak{B}$ is isomorphic to the category
$\mathfrak{L}_0$. The category isomorphism is given by the contravariant
functor $\mathfrak{F}: \mathfrak{B} \rightarrow \mathfrak{L}_0$
defined by $\mathfrak{F}(\mathfrak{b}) = \mathfrak{b}^\ast$ and
$\mathfrak{F}(\phi) = \phi^\ast$. The inverse of $\mathfrak{F}$ is
given by $\mathfrak{F}^{-1}(\mathfrak{g}) = \mathfrak{g}_\ast$ and
$\mathfrak{F}^{-1}(\psi) = \psi^\ast|_{(\mathfrak{g}_2)_\ast}$,
where $\psi: \mathfrak{g}_1 \rightarrow \mathfrak{g}_2$.
\end{theorem}

\noindent \textbf{Proof}. If $\mathfrak{b}$ is a Banach
Lie-Poisson space, then $\mathfrak{F}(\mathfrak{b}) =
\mathfrak{b}^\ast$ is a Banach Lie algebra that admits
$\mathfrak{b}$ as a predual and, according to Theorem \ref{general
theorem}, $\operatorname{ad}^\ast_{\mathfrak{b}^\ast} \mathfrak{b}
\subset \mathfrak{b}$. Thus $\mathfrak{F}(\mathfrak{b})$ is indeed
an object in the category $\mathfrak{L}_0$. If $\phi:
\mathfrak{b}_1 \rightarrow \mathfrak{b}_2$ is a linear Poisson map
let us show that $\mathfrak{F}(\phi) = \phi^\ast:
\mathfrak{b}_2^\ast \rightarrow \mathfrak{b}_1^ \ast$ is a Banach
Lie algebra homomorphism. First, $\phi^\ast$ is a linear
continuous map between Banach spaces. Second, since the Lie
bracket on $\mathfrak{b}_2^\ast$ is defined by $[x, y]_2 = \{x,
y\}_2$ and similarly for $\mathfrak{b}_1$, we get
\[
\phi^\ast[x, y]_2 = \phi^\ast\{x, y\}_2 = \{ \phi^\ast x,
\phi^\ast y \}_1 = [\phi^\ast x, \phi^\ast y]_1,
\]
which shows that $\phi^\ast$ a homomorphism of Banach Lie
algebras. Finally, the dual of $\mathfrak{F}(\phi)$, that is,
$\phi^{\ast \ast}: \mathfrak{b}_1^{\ast \ast} \rightarrow
\mathfrak{b}_2^{\ast \ast}$ satisfies $\phi^{\ast
\ast}|_{\mathfrak{b}_1} = \phi$. Thus $\mathfrak{F}(\phi)$ is
indeed a morphism in the category $\mathfrak{L}_0$. Since duality
reverses the direction of the arrows and the order of the
composition, $\mathfrak{F}$ is a contravariant functor.

Conversely, consider the functor $\mathfrak{F}^{-1}:
\mathfrak{L}_0 \rightarrow \mathfrak{B}$ and let $\mathfrak{g}$ be
an object of $\mathfrak{L}_0$. By Theorem \ref{general theorem},
$\mathfrak{F}^{-1}(\mathfrak{g}) = \mathfrak{g}_\ast$ is a Banach
Lie-Poisson space, that is, an object of $\mathfrak{B}$. If $\psi:
\mathfrak{g}_1 \rightarrow \mathfrak{g}_2$ is a morphism in the
category $\mathfrak{L}_0$, then let us show that
$\mathfrak{F}^{-1}(\psi) = \psi^\ast|_{(\mathfrak{g}_2)_\ast}$ is
a linear Poisson map. Let $f, g$ be smooth functions on
$(\mathfrak{g}_1)_\ast$. From \eqref{general LP} and using the fact that
$\psi$ is morphism of Banach Lie algebras and that
$\psi^\ast|_{(\mathfrak{g}_2)_\ast}$ is a linear map, we get for every $b \in
(\mathfrak{g}_2)_\ast$
\begin{align*}
\{ f \circ \psi^\ast|_{(\mathfrak{g}_2)_\ast}, 
g \circ \psi^\ast|_{(\mathfrak{g}_2)_\ast} \}_2(b) 
&= \left\langle [D(f\circ \psi^\ast|_{(\mathfrak{g}_2)_\ast})(b), 
D(g \circ \psi^\ast|_{(\mathfrak{g}_2)_\ast})(b)]_2, b \right\rangle \\ &= 
\left\langle [Df(\psi^\ast(b)) \circ 
\psi^\ast|_{(\mathfrak{g}_2)_\ast}, Dg(\psi^\ast(b)) \circ
\psi^\ast|_{(\mathfrak{g}_2)_\ast} ]_2,  b \right\rangle \\
&= \left\langle [\psi\left(Df(\psi^\ast(b))\right),
\psi\left(Dg(\psi^\ast(b))\right)]_2, b \right\rangle \\
&= \left\langle \psi\left([Df(\psi^\ast(b)), Dg(\psi^\ast(b))]_1 \right), b 
\right\rangle \\
&= \left\langle [Df(\psi^\ast(b)), Dg(\psi^\ast(b))]_1, \psi^\ast(b)
\right\rangle \\&= \{ f, g \}_1(\psi^\ast|_{(\mathfrak{g}_2)_\ast}(b))
\end{align*}
which shows that $\psi^\ast|_{(\mathfrak{g}_2)_\ast}: 
{(\mathfrak{g}_2)_\ast} \rightarrow {(\mathfrak{g}_1)_\ast}$ is a morphism of 
Banach Lie-Poisson spaces. The functor $\mathfrak{F}^{-1}$ is contravariant since its
action on morphisms is given by duality.

Finally, it is clear the functors $\mathfrak{F}$ and
$\mathfrak{F}^{-1}$ are inverses of each other. $\quad
\blacksquare$
\medskip

We turn now to the study of the internal structure of morphisms of
Banach Lie-Poisson spaces.

\begin{proposition}
\label{decomposition}
Let $\phi:\mathfrak{b}_1 \rightarrow
\mathfrak{b}_2$ be a linear Poisson map between Banach Lie-Poisson
spaces and assume that $\operatorname{im} \phi$ is closed in
$\mathfrak{b}_2$. Then the Banach space $\mathfrak{b}_1/ \ker
\phi$ is predual to $\mathfrak{b}_2^\ast/ \ker \phi^\ast$, that
is, $(\mathfrak{b}_1/ \ker \phi)^\ast \cong \mathfrak{b}_2^\ast/
\ker \phi^\ast$. In addition, $\mathfrak{b}_2^\ast/ \ker
\phi^\ast$ is a Banach Lie algebra satisfying the condition
$\operatorname{ad}_{[x]}^\ast \left(\mathfrak{b}_1/ \ker \phi
\right) \subset \mathfrak{b}_1/ \ker \phi$ for all $[x] \in
\mathfrak{b}_2^\ast/ \ker \phi^\ast$ and $\mathfrak{b}_1/\ker
\phi$ is a Banach Lie-Poisson space. Moreover, the following
properties hold:
\begin{enumerate}
\item[{\rm (i)}] the quotient map $\pi: \mathfrak{b}_1
\rightarrow \mathfrak{b}_1/\ker \phi$ is a surjective
linear Poisson map;
\item[{\rm (ii)}] the map $\iota: \mathfrak{b}_1/\ker \phi
\rightarrow \mathfrak{b}_2$ defined by
$\iota([b]):= \phi(b)$, where $b \in \mathfrak{b}_1$ and
$[b] \in \mathfrak{b}_1/\ker \phi$ is an injective linear Poisson map;
\item[{\rm (iii)}] the decomposition $\phi= \iota \circ \pi$
into a surjective and an injective linear Poisson
map is valid.
\end{enumerate}
\end{proposition}

\noindent \textbf{Proof}. We define the pairing
$\langle \cdot , \cdot \rangle: \mathfrak{b}_2^\ast/ \ker
\phi^\ast \times \mathfrak{b}_1/ \ker \phi \rightarrow
\mathbb{C}$ (or $\mathbb{R}$) by
\begin{equation}
\label{quotient pairing}
\langle [x], [b] \rangle := \langle x, \phi(b) \rangle_2
= \langle \phi^\ast(x), b \rangle_1
\end{equation}
where $[x] \in \mathfrak{b}_2^\ast/ \ker \phi^\ast$,
$[b] \in \mathfrak{b}_1/ \ker \phi$, and $\langle \cdot ,
\cdot \rangle_i : \mathfrak{b}_i ^\ast \times
\mathfrak{b}_i \rightarrow \mathbb{C}$ (or $\mathbb{R}$), $i =
1,2$ are the pairings between the given Banach Lie-Poisson
spaces and their duals. This pairing is correctly
defined since it does not depend on the choice of the
representatives $x \in \mathfrak{b}_2$ and $b \in
\mathfrak{b}_1$. One has
\[
|\langle [x], [b] \rangle| \leq \|\phi\|\, \|[x]\|\, \|[b]\|
\]
and if $\langle [x], [b] \rangle = 0$ for each $x \in
\mathfrak{b}_2$ ($b \in \mathfrak{b}_1$) then $[b] = [0]$ ($[x] =
[0]$). Thus (\ref{quotient pairing}) defines a continuous weakly
non degenerate pairing and therefore the map
\[
[x] \in \mathfrak{b}_2^\ast/ \ker \phi^\ast \mapsto \langle [x],
[\cdot] \rangle = \langle \phi^\ast (x), \cdot \rangle _1 =
\langle x, \phi(\cdot ) \rangle _2 \in \left(\mathfrak{b}_1/ \ker
\phi \right)^\ast
\]
is a continuous linear injective map of Banach spaces. To show
that this map is surjective, we need to find for a given $\alpha
\in \left(\mathfrak{b}_1/ \ker \phi \right)^\ast$ an $[x] \in
\mathfrak{b}_2 ^\ast/ \ker \phi^\ast$ such that $\langle [x], [b]
\rangle = \langle \phi^\ast (x), b \rangle _1 = \langle x, \phi(b)
\rangle _1 = \alpha([b])$ for all $b \in \mathfrak{b}_1$. Since
the range $\operatorname{im} \phi$ is closed in $\mathfrak{b}_2$,
it is a Banach subspace and hence the map $\Phi: [b] \in
\mathfrak{b}_1/\ker \phi \mapsto \phi(b) \in \operatorname{im}
\phi$ is a Banach space isomorphism. Thus $\alpha \circ \Phi^{-1}
\in (\operatorname{im} \phi)^\ast$. Let $x \in \mathfrak{b}_2^\ast$ be an extension
of $\alpha \circ \Phi^{-1}$ to $\mathfrak{b}_2$. Then we have for any $b \in
\mathfrak{b}_1$
\[
\langle [x], [b] \rangle = \langle x, \phi(b) \rangle =
\langle \alpha \circ \Phi^{-1}, \phi(b) \rangle =
\alpha([b]).
\]
Thus the Banach space $\mathfrak{b}_2^\ast/ \ker \phi^\ast$
is isomorphic to the dual of  $\mathfrak{b}_1/ \ker \phi$.

The space $\mathfrak{b}_2^\ast/ \ker \phi^\ast$ is a Banach
Lie algebra because $\ker \phi^\ast$ is an ideal
in the Banach Lie algebra $\mathfrak{b}_2^\ast$
(since $\phi^\ast : \mathfrak{b}_2^\ast \rightarrow
\mathfrak{b}_1 ^\ast$ is a morphism of Banach Lie algebras).
Finally, since $\phi : \mathfrak{b}_1 \rightarrow
\mathfrak{b}_2$ is a linear Poisson map, we have
\begin{equation}
\label{linear Poisson morphism condition}
{\operatorname{ad}^2}_x^\ast \phi (b) =
\phi\left({\operatorname{ad}^1}_{x \circ \phi}^\ast b \right),
\end{equation}
for any $x \in \mathfrak{b}_2^\ast$, $b \in \mathfrak{b}_1$, and
where $\operatorname{ad}^i$ denotes the adjoint operator in the
Banach Lie algebra $\mathfrak{b}^\ast_i$, $i = 1,2$. Here we have
used the fact that ${\operatorname{ad}^1}_{x \circ \phi}^\ast
\mathfrak{b}_1 \subset \mathfrak{b}_1$ for any $x \in
\mathfrak{b}_2^\ast$. From (\ref{linear Poisson morphism
condition}) and ${\operatorname{ad}^2}^\ast _x \mathfrak{b}_2
\subset \mathfrak{b}_2$ for all $x \in \mathfrak{b}_2^\ast$ we
conclude that for all $b \in \mathfrak{b}_1$ we have
\begin{align*}
\left\langle [y], {\operatorname{ad}}^\ast _{[x]} [b]
\right\rangle &= \left\langle \left[ [x], [y] \right], [b]
\right\rangle = \left\langle \left[ [x, y] \right], [b]
\right\rangle = \left\langle \left[x, y \right],
\phi(b) \right\rangle_2 \\
&= \left\langle y, {\operatorname{ad}^2}^\ast _x \phi(b)
\right\rangle_2 =
\left\langle y, \phi\left({\operatorname{ad}^1}_{x \circ
\phi}^\ast b \right) \right\rangle_2 = \left\langle [y],
\left[{\operatorname{ad}^1}_{x \circ \phi}^\ast b \right]
\right\rangle
\end{align*}
for each $y \in \mathfrak{b}_2^\ast$. This implies that
\[
{\operatorname{ad}}^\ast _{[x]} [b] =
\left[{\operatorname{ad}^1}_{x \circ \phi}^\ast b \right] \in
\mathfrak{b}_1/ \ker \phi
\]
for all $[x] \in \mathfrak{b}_2 ^\ast/ \ker \phi^\ast$,
$[b] \in \mathfrak{b}_1/ \ker \phi$, and thus
\[
{\operatorname{ad}}^\ast _{[x]}\left(\mathfrak{b}_1/ \ker
\phi\right) \subset \mathfrak{b}_1/ \ker \phi
\]
for all $[x] \in \mathfrak{b}_2 ^\ast/ \ker \phi^\ast$.
Thus $\mathfrak{b}_2 ^\ast/ \ker \phi^\ast$ is an object 
in the category $\mathfrak{L}_0$. Theorem \ref{general theorem} (or Theorem
\ref{functor}) guarantees then that the quotient Banach space $\mathfrak{b}_1/ \ker
\phi$ is a Banach Lie-Poisson space.

Endow the Banach subspace $\operatorname{im} \phi$ with the Banach
Lie-Poisson structure making the Banach space isomorphism $\Phi:
\mathfrak{b}_1/\ker \phi \rightarrow \operatorname{im} \phi$ into
a linear Poisson isomorphism. Thus $\Phi^\ast: (\operatorname{im}
\phi)^\ast \rightarrow (\mathfrak{b}_1/\ker \phi)^\ast$ is an
isomorphism in the category $\mathfrak{L}_0$. Since $\phi = \Phi \circ \pi:
\mathfrak{b}_1 \rightarrow \operatorname{im} \phi$ is a linear
Poisson map by hypothesis, it follows that $\phi^\ast = \pi^\ast \circ
\Phi^\ast$ is a morphism in the category $\mathfrak{L}_0$ which then implies
that $\pi^\ast$ is also a morphism in the category $\mathfrak{L}_0$.  By
Theorem \ref{functor} this is equivalent to
the fact that $\pi$ is a linear Poisson map thereby proving
property (i) in the statement of the proposition.

Define $\iota: \mathfrak{b}_1/\ker \phi \rightarrow
\mathfrak{b}_2$ to be the composition of the inclusion
$\operatorname{im} \phi \hookrightarrow \mathfrak{b}_2$ with the
isomorphism $\Phi: \mathfrak{b}_1/\ker \phi \rightarrow
\operatorname{im} \phi$. The definition of $\Phi$ is equivalent 
to the equality $\phi = \Phi\circ \pi$ thought of as a map from
$\mathfrak{b}_1$ to $\operatorname{im}\phi$. Composing this identity
on the left with the inclusion $\operatorname{im} \phi \hookrightarrow
\mathfrak{b}_2$ yields $\phi = \iota \circ \pi$ which proves
property (iii). 

To prove part (ii), let $f,g \in C^\infty
(\mathfrak{b}_2)$. Then $f\circ \iota, g \circ \iota \in
C^\infty(\mathfrak{b}_1/\ker \phi)$. Since $\pi : \mathfrak{b}_1
\rightarrow \mathfrak{b}_1/ \ker \phi$ is a surjective linear
Poisson map and $\phi = \iota \circ \pi$, the relation
\[
\{f\circ \iota, g \circ \iota \} \circ \pi = \{f\circ \iota \circ
\pi, g \circ \iota \circ \pi \}_1 = \{f\circ \phi, g \circ \phi
\}_1 = \{f, g \}_2 \circ \phi = \{f, g \}_2 \circ \iota \circ \pi
\]
implies that $\{f\circ \iota, g \circ \iota \} =
\{f, g \}_2 \circ \iota$, that is, $\iota:
\mathfrak{b}_1/\ker \phi \rightarrow \mathfrak{b}_2$
is an injective linear Poisson map.
$\quad \blacksquare$
\medskip

Proposition \ref{decomposition} reduces the investigation of
linear Poisson maps with closed range between Banach Lie-Poisson
spaces to the study of surjective  and injective linear Poisson maps.

Consider therefore the surjective linear continuous map $\pi:
\mathfrak{b}_1 \rightarrow \mathfrak{b}_2$, where $\mathfrak{b}_1$
is a Banach Lie-Poisson space and $\mathfrak{b}_2$ is just a
Banach space with no additional structure. The dual map 
$\pi^\ast: \mathfrak{b}_2^\ast \rightarrow \mathfrak{b}_1^\ast$ is
therefore an injective continuous linear map of Banach spaces. The
space $\operatorname{im} \pi ^\ast$ coincides with the Banach
subspace of linear continuous functionals on $\mathfrak{b}_1$ that
vanish on $\ker \pi$, which is closed in $\mathfrak{b}_1$.
Thus $\im \pi^\ast$ is a closed subspace of $\mathfrak{b}_1^\ast$.

Assume next that $\im \pi^\ast$ is also closed under the Lie
bracket operation $[\cdot, \cdot ]_1$ of $\mathfrak{b}_1^\ast$.
Then $\pi^\ast : \mathfrak{b}_2^\ast \rightarrow \im \pi^\ast$ is
a Banach space isomorphism and, declaring it to be also a Lie
algebra morphism, it follows that there is a Banach Lie algebra
structure $[\cdot, \cdot] _2$ on $\mathfrak{b}_2^\ast$ and that
$\pi^\ast :(\mathfrak{b}_2^\ast, [\cdot, \cdot ]_2 ) \rightarrow
(\mathfrak{b}_1^\ast, [\cdot, \cdot ]_1)$ is a Banach Lie algebra
morphism.

Let $\tilde{x}, \tilde{y} \in \mathfrak{b}_2^\ast$ and let
$\pi(b) = \tilde {b} \in \mathfrak{b}_2$. Then
\begin{align*}
\left\langle \tilde{y}, {\operatorname{ad}^2}^\ast_{\tilde{x}}
\tilde{b} \right\rangle_2 &= \left\langle [\tilde{x},
\tilde{y}]_2, \tilde{b} \right\rangle_2 = \left\langle [\tilde{x},
\tilde{y}]_2, \pi(b) \right\rangle_2 = \left\langle
\pi^\ast([\tilde{x}, \tilde{y}]_2), b \right\rangle_1\\ & =
\left\langle \left[\pi^\ast(\tilde{x}),
\pi^\ast(\tilde{y})\right]_1, b \right\rangle_1 = \left\langle
\pi^\ast(\tilde{y}),
{\operatorname{ad}^1}^\ast_{\pi^\ast(\tilde{x})} b \right\rangle_1
= \left\langle \tilde{y}, \pi
\left({\operatorname{ad}^1}^\ast_{\pi^\ast(\tilde{x})} b \right)
\right\rangle_2;
\end{align*}
the last equality is a consequence of the inclusion
${\operatorname{ad}^1}^\ast_{\pi^\ast(\tilde{x})} \mathfrak{b}_1
\subset \mathfrak{b}_1$ for all $\tilde{x} \in
\mathfrak{b}_2^\ast$ which is insured by the fact that
$\mathfrak{b}_1$ is a Banach Lie-Poisson space. Since this
relation holds for any $\tilde{y} \in \mathfrak{b}_2$, we conclude
that
\[
{\operatorname{ad}^2}^\ast_{\tilde{x}} \tilde{b}
= \pi \left({\operatorname{ad}^1}^\ast_{\pi^\ast(\tilde{x})}
b \right)
\]
for any $\tilde{x} \in \mathfrak{b}_2^\ast$ and any $\tilde{b} \in
\mathfrak{b}_2$. This shows that
${\operatorname{ad}^2}^\ast_{\tilde{x}} \mathfrak{b}_2 \subset
\mathfrak{b}_2$ for any $\tilde{x} \in \mathfrak{b}_2^\ast$, and
hence, by Theorem \ref{general theorem}, $\mathfrak{b}_2$ is a
Banach Lie-Poisson space or, equivalently, $\mathfrak{b}_2^\ast$ is 
an object in the category $\mathfrak{L}_0$.

The map $\pi^\ast$ is a morphism of Banach Lie algebras. In addition, its dual
$\pi^{\ast \ast}: \mathfrak{b}_1^{\ast \ast} \rightarrow \mathfrak{b}_2^{\ast \ast}$
has the property $\pi^{\ast \ast}(\mathfrak{b}_1) \subset \mathfrak{b}_2$. Indeed,
for any $b_1 \in \mathfrak{b}_1$ and $\beta_2 \in \mathfrak{b}_2^\ast$, the
definition of the dual of a linear map gives
\[
\langle \pi ^{\ast \ast}(b_1), \beta_2 \rangle_2' = 
\langle b, \pi^\ast(\beta_2) \rangle_1 =
\langle \pi(b_1), \beta_2 \rangle_2
\]
where $\langle \cdot, \cdot \rangle_2' :\mathfrak{b}_2^{\ast \ast} 
\times \mathfrak{b}_2^\ast \rightarrow \mathbb{R}$ 
(or $\mathbb{C}$) is the canonical pairing between a 
Banach space and its dual and similarly for  
$\langle \cdot, \cdot\rangle_1 :\mathfrak{b}_1 
\times \mathfrak{b}_1^\ast \rightarrow \mathbb{R}$ 
(or $\mathbb{C}$) and $\langle \cdot, \cdot\rangle_2 :\mathfrak{b}_2
\times \mathfrak{b}_2^\ast \rightarrow \mathbb{R}$ 
(or $\mathbb{C}$). This shows that $\pi ^{\ast \ast}(b_1) = \pi(b_1) \in 
\mathfrak{b}_2$. Therefore $\pi^\ast$ is a morphism in the
category $\mathfrak{L}_0$ and, by Theorem \ref{functor},
$\pi: \mathfrak{b}_1 \rightarrow \mathfrak{b}_2$ is
a linear Poisson map. In addition, since $\pi$ is also surjective,
the Banach Lie-Poisson structure on $\mathfrak{b}_2$ is unique.
Therefore, following e.g. Vaisman [1994], we shall call this Banach Lie-Poisson
structure on $\mathfrak{b}_2$ {\bfi coinduced\/} by the surjective mapping $\pi$.

The above proves the ``only if" part of the following
proposition; the converse is an easy verification.

\begin{proposition}
\label{surjective LP map}
Let $(\mathfrak{b}_1, \{\cdot, \cdot \})$ be a Banach Lie-Poisson space 
and let $\pi : \mathfrak{b}_1 \rightarrow \mathfrak{b}_2$ be a continuous linear
surjective map onto the Banach space $\mathfrak{b}_2$. Then $\mathfrak{b}_2$
carries the Banach Lie-Poisson structure coinduced by 
$\pi$  if and only if $\im \pi^\ast \subset \mathfrak{b}_1^\ast$ is closed under the
Lie bracket $[\cdot, \cdot]_1$ of $\mathfrak{b}_1^\ast$.  The map
$\pi^\ast :\mathfrak{b}_2^\ast \rightarrow \mathfrak{b}_1^\ast$ is
a Banach Lie algebra morphism whose dual $\pi^{\ast \ast}: \mathfrak{b}_1^{\ast
\ast} \rightarrow \mathfrak{b}_2^{\ast \ast}$ maps $\mathfrak{b}_1$ into 
$\mathfrak{b}_2$.
\end{proposition}

\begin{example}
\label{example following Proposition}
\normalfont
Let $(\mathfrak{g}, [\cdot, \cdot])$ be a
complex Banach Lie algebra admitting a predual $\mathfrak{g}_\ast$
satisfying $\operatorname{ad}^\ast _x \mathfrak{g}_\ast \subset
\mathfrak{g}_\ast$ for every $x \in \mathfrak{g}$. Then, by Theorem \ref
{general theorem}, the predual $\mathfrak{g}_\ast$ admits a
holomorphic Banach Lie-Poisson structure, whose holomorphic Poisson tensor
$\varpi$ is given by \eqref{general LP}. We shall work with the
realification $({\mathfrak{g}_\ast}_\mathbb{R}, \varpi_\mathbb{R})$ of
$(\mathfrak{g}_\ast, \varpi)$ in the sense of \S \ref{section: banach
poisson manifolds}. We want to construct a real Banach space
$\mathfrak{g}_\ast ^\sigma$ with a real Banach Lie-Poisson structure
$\varpi_\sigma$ such that $\mathfrak{g}_\ast ^\sigma \otimes \mathbb{C} =
\mathfrak{g}_\ast$ and $\varpi_\sigma$ is coinduced from
$\varpi_\mathbb{R}$ in the sense of Proposition \ref{surjective LP map}.
To achieve this, introduce a continuous $\mathbb{R}$-linear map $\sigma:
{\mathfrak{g}_\ast}_\mathbb{R} \rightarrow {\mathfrak{g}_\ast}_\mathbb{R}$
satisfying the properties: 
\begin{enumerate}
\item[{\rm (i)}] $\sigma^2 = id$;
\item[{\rm (ii)}] the dual map $\sigma^\ast: \mathfrak{g}_\mathbb{R}
\rightarrow \mathfrak{g}_\mathbb{R}$ defined by 
\begin{equation}
\label{complex pairing}
\langle \sigma^\ast z, b \rangle = \overline{\langle  z, \sigma b
\rangle}
\end{equation}
for $z \in \mathfrak{g}_\mathbb{R}$, $b \in
{\mathfrak{g}_\ast}_\mathbb{R}$ and where $\langle \cdot, \cdot \rangle$
is the pairing between the complex Banach spaces $\mathfrak{g}$
and  ${\mathfrak{g}_\ast}$, is a homomorphism of the Lie algebra
$(\mathfrak{g}_\mathbb{R}, [\cdot, \cdot ])$;
\item[{\rm (iii)}] $\sigma \circ I + I \circ \sigma = 0$,
where $I:\mathfrak{g}_\mathbb{R} \rightarrow \mathfrak{g}_\mathbb{R}$ is
defined by
\begin{equation}
\label{definition of I}
\langle z, Ib \rangle:= \langle I^\ast z, b \rangle := i \langle z, b
\rangle
\end{equation}
for $z \in \mathfrak{g}_\mathbb{R}$, $b \in {\mathfrak{g}_\ast}_\mathbb{R}$.
\end{enumerate}

Consider the projectors
\begin{equation}
\label{example projectors}
R:= \frac{1}{2}( id + \sigma) \qquad R^\ast:= \frac{1}{2}( id +
\sigma^\ast)
\end{equation}
and define $\mathfrak{g}_\ast ^\sigma := \operatorname{im} R$, 
$\mathfrak{g}^\sigma := \operatorname{im} R^\ast$. Then one has the
splittings
\begin{equation}
\label{example decompositions}
{\mathfrak{g}_\ast}_\mathbb{R} = \mathfrak{g}_\ast ^\sigma \oplus
I\mathfrak{g}_\ast ^\sigma \quad \text{and} \quad 
\mathfrak{g}_\mathbb{R} = \mathfrak{g} ^\sigma \oplus
I\mathfrak{g}^\sigma
\end{equation}
into real Banach subspaces. One can identify canonically the splittings
\eqref{example decompositions} with the splittings
\begin{equation}
\label{example second splitting}
\mathfrak{g}_\ast ^\sigma \otimes_\mathbb{R} \mathbb{C}
= \left(\mathfrak{g}_\ast^\sigma
\otimes_\mathbb{R} \mathbb{R}\right) \oplus \left(\mathfrak{g}_\ast ^\sigma
\otimes_\mathbb{R} \mathbb{R} i\right).
\end{equation}
Thus one obtains isomorphisms $\mathfrak{g}_\ast ^\sigma \otimes_\mathbb{R}
\mathbb{C} \cong \mathfrak{g}_\ast$ and $\mathfrak{g} ^\sigma
\otimes_\mathbb{R} \mathbb{C} \cong \mathfrak{g}$ of complex Banach spaces. 

For any $x, y \in \mathfrak{g}_\mathbb{R}$ one has 
\begin{equation}
\label{example R star condition}
[R^\ast x, R^\ast y ] = R^\ast [x, R^\ast y]
\end{equation}
and thus $\mathfrak{g}^\sigma$ is a real Banach Lie subalgebra
of $\mathfrak{g}_\mathbb{R}$. From 
\begin{align}
\label{example real identity}
\operatorname{Re}\langle z, b \rangle &= \langle R^\ast z R b \rangle +
\langle I^\ast R^\ast I^\ast z, IRI b \rangle \nonumber \\
&= \langle R^\ast z R b \rangle + \langle (1 - R^\ast)z, (1-R) b \rangle
\end{align}
for all $z \in \mathfrak{g}_\mathbb{R}$ and all $b \in
{\mathfrak{g}_\ast}_\mathbb{R}$, where for the last equality we used $R =
1 + IRI$ and $R^\ast = 1 + I^\ast R^ \ast I^\ast$, one concludes that
the annihilator $(\mathfrak{g}_\ast^\sigma)^\circ$ of
$\mathfrak{g}_\ast^\sigma$ in $\mathfrak{g}_\mathbb{R}$ equals
$I^\ast\mathfrak{g}^\ast$. Therefore $\mathfrak{g}_\ast ^\sigma$ is the
predual of $\mathfrak{g}^\sigma$.

Taking into account all of the above facts we conclude from Proposition
\ref{surjective LP map} that $\mathfrak{g}_\ast^\sigma$ carries a real
Banach Lie-Poisson structure $\{\cdot, \cdot \}_{\mathfrak{g}_\ast^\sigma}$
coinduced by $R: {\mathfrak{g}_\ast}_\mathbb{R} \rightarrow
\mathfrak{g}_\ast^\sigma$. According to \eqref{example real identity}, the
bracket $\{\cdot, \cdot \}_{\mathfrak{g}_\ast^\sigma}$ is given by
\begin{equation}
\label{example real Poisson bracket}
\{f, g \}_{\mathfrak{g}_\ast^\sigma} (\rho) = \langle [df(\rho),
dg(\rho)], \rho \rangle, 
\end{equation}
where $\rho \in \mathfrak{g}_\ast^\sigma$ and the pairing on the right is
between $\mathfrak{g}_\ast^\sigma$ and
$\mathfrak{g}^\sigma$. In addition, for any real valued functions $f, g \in
C^\infty (\mathfrak{g}_\ast^\sigma)$ and any $b \in
{\mathfrak{g}_\ast}_\mathbb{R}$ we have 
\begin{align*}
&\{f \circ R , g \circ R\}_{\mathfrak{g}_\mathbb{R}} (b) 
= \operatorname{Re} \langle [d(f \circ R)(b), d(g \circ R)(b)], b
\rangle \\
& \; =\langle R^\ast[d(f \circ R)(b), d(g \circ R)(b)], R(b) \rangle  +
\langle (1- R^\ast) [d(f \circ R)(b), d(g \circ R)(b)], (1-R)b
\rangle \\
&\; = \langle R^\ast[R^\ast df(R(b)), R^\ast dg(R(b))], R(b)
\rangle  + \langle (1- R^\ast) [R^\ast df(R(b)), R^\ast dg(R(b))], (1-R)b
\rangle \\
& \; = \langle [df (R(b)), dg (R(b))], R(b)
\rangle = \{f, g\}_{\mathfrak{g}_\ast^\sigma}(R(b)),
\end{align*}
where we have used \eqref{example R star condition}. The above computation
proves, independently of Proposition \ref{surjective LP map}, that $R:
{\mathfrak{g}_\ast}_\mathbb{R} \rightarrow \mathfrak{g}_\ast^\sigma$ is a
linear Poisson map. \quad $\blacklozenge$
\end{example}

Next we investigate the case of injective linear Poisson maps.

\begin{proposition}
\label{injective LP map}
Let $\mathfrak{b}_1$ be a Banach space,
$(\mathfrak{b}_2, \{\cdot, \cdot \}_2)$ be a Banach Lie-Poisson
space, and $\iota : \mathfrak{b}_1 \rightarrow \mathfrak{b}_2$ be
an injective continuous linear map with closed range. Then
$\mathfrak{b}_1$ carries a unique Banach Lie-Poisson structure such that 
$\iota$ is a linear Poisson map if and only if
$\ker \iota^\ast$ is an ideal in the Banach Lie algebra
$\mathfrak{b}_2^\ast$.
\end{proposition}

\noindent \textbf{Proof}. Assume that $\ker \iota^\ast$ is an ideal in the 
Banach Lie algebra $\mathfrak{b}_2^\ast$. Denote by $[\cdot, \cdot]_2$ the Lie
bracket of the Banach Lie algebra $\mathfrak{b}_2^\ast$. Since $\iota :\mathfrak{b}_1
\rightarrow \mathfrak{b}_2$ is an injective linear continuous map,
its adjoint $\iota^\ast: \mathfrak{b}_2^\ast \rightarrow
\mathfrak{b}_1^\ast$ is a surjective linear continuous map inducing the Banach space
isomorphism $[\iota^\ast]: \mathfrak{b}_2^\ast/ \ker \iota^\ast
\tilde{\longrightarrow} \mathfrak{b}_1^\ast$. Since $\ker \iota^\ast$ is an ideal in
the Banach Lie algebra $\mathfrak{b}_2^\ast$, it follows that
$\mathfrak{b}_2^\ast/ \ker \iota^\ast$ is a Banach Lie algebra.
The isomorphism $[\iota^\ast]$ induces a Banach Lie algebra
structure $[\cdot, \cdot ]_1$ on $\mathfrak{b}_1^\ast$. The linear
map $\iota^\ast :\mathfrak{b}_2^\ast \rightarrow
\mathfrak{b}_1^\ast$ becomes a Banach Lie algebra homomorphism.

For each $x, y \in \mathfrak{b}_2^\ast$ and each $b \in
\mathfrak{b}_1$ we have
\begin{align*}
\left\langle y,\, \iota^{\ast \ast}
\left({\operatorname{ad}^1}^\ast_{\iota^\ast(x)} b \right)
\right\rangle_2 &= \left\langle \iota^\ast(y),\,
{\operatorname{ad}^1}^\ast_{\iota^\ast(x)} b \right\rangle_1 =
\left\langle [\iota^\ast(x), \iota^\ast(y)]_1, b \right\rangle_1 \\
&= \left\langle \iota^\ast \left([x, y]_2\right), b
\right\rangle_1 = \left\langle [x, y]_2, \iota(b) \right\rangle =
\left\langle y,\, {\operatorname{ad}^2}^\ast_x \iota(b)
\right\rangle
\end{align*}
and hence we get the following identity in $\mathfrak{b}_2^{\ast\ast}$
\begin{equation}
\label{double star}
\iota^{\ast \ast} \left({\operatorname{ad}^1}^\ast_{\iota^\ast(x)}
b \right) = {\operatorname{ad}^2}^\ast_x \iota(b)
\end{equation}
for any $x \in \mathfrak{b}_2^\ast$ and any $b \in
\mathfrak{b}_1$. 

Let us now prove that 
\begin{equation}
\label{second condition}
{\operatorname{ad}^2}^\ast_{\mathfrak{b}_2^\ast}
\iota(\mathfrak{b}_1) \subset \iota(\mathfrak{b}_1),
\end{equation}
where ${\operatorname{ad}^2}^\ast$ denotes the coadjoint action of
$\mathfrak{b}_2^\ast$ on $\mathfrak{b}_2^{\ast \ast}$. We begin by noticing  
that $\ker \iota^\ast = [\iota(\mathfrak{b}_1)]^\circ$, where
$[\iota(\mathfrak{b}_1)]^\circ$ is the annihilator of 
$\iota(\mathfrak{b}_1)$ in $\mathfrak{b}_2^\ast$. Taking the annihilator of 
this identity in $\mathfrak{b}_2$ one obtains
\begin{equation}
\label{double annihilator}
[\ker \iota^\ast]^\circ = [\iota(\mathfrak{b}_1)]^{\circ \circ}
= \iota(\mathfrak{b}_1),
\end{equation}
where the last equality follows by the closedness of $\iota(\mathfrak{b}_1)$ in
$\mathfrak{b}_2$. By the definition of Banach Lie-Poisson spaces, 
${\operatorname{ad}^2}^\ast _{\mathfrak{b}_2^\ast} \mathfrak{b}_2 \subset 
\mathfrak{b}_2$. Since $\ker \iota^\ast$ is an ideal in $\mathfrak{b}_2^\ast$ 
it follows that ${\operatorname{ad}^2}^\ast _{\mathfrak{b}_2^\ast} [\ker 
\iota^\ast]^\circ \subset [\ker \iota^\ast]^\circ$ and thus \eqref{double 
annihilator} implies \eqref{second condition}.

By \eqref{double star} and \eqref{second condition},  ${\operatorname{ad}^2}
^\ast_x \iota(b)
\in \iota(\mathfrak{b}_1)$ and thus
\[
\iota^{\ast \ast} \left({\operatorname{ad}^1}^\ast_{\iota^\ast(x)}
b \right) \in \iota(\mathfrak{b}_1)
\]
for any $x \in \mathfrak{b}_2^\ast$ and any $b \in
\mathfrak{b}_1$. The double adjoint $\iota^{\ast \ast}:
\mathfrak{b}_1^{\ast \ast} \rightarrow \mathfrak{b}_2^{\ast \ast}$
is a injective continuous linear map (since $\iota^\ast:
\mathfrak{b}_2^\ast \rightarrow \mathfrak{b}_1^\ast$ is a
surjective continuous linear map) and $\iota$ maps $\mathfrak{b}_1
\subset \mathfrak{b}_1^{\ast \ast}$ into $\mathfrak{b}_2 \subset
\mathfrak{b}_2^{\ast \ast}$. This shows that
${\operatorname{ad}^1}^\ast_{\iota^\ast(x)} \mathfrak{b}_1 \subset
\mathfrak{b}_1$. Applying Theorem \ref{general theorem} we
conclude that $\mathfrak{b}_1$ is a Banach Lie-Poisson space and
that $\iota: \mathfrak{b}_1 \rightarrow \mathfrak{b}_2$ is an
injective linear Poisson map. Uniqueness of the Poisson structure
on $\mathfrak{b}_1$ follows from the injectivity of $\iota$.

Conversely, let us assume that $\mathfrak{b}_1$ is a Banach
Lie-Poisson space and that $\iota: \mathfrak{b}_1 \rightarrow
\mathfrak{b}_2$ is a linear Poisson map. Then $\iota^\ast :
\mathfrak{b}_2^\ast \rightarrow \mathfrak{b}_1^\ast$ is a
homomorphism of Banach Lie algebras and therefore its kernel is an
ideal in $\mathfrak{b}_2^\ast$. $\quad \blacksquare$
\medskip

Proposition \ref{injective LP map} allows one 
to introduce a unique Banach Lie-Poisson structure on $\mathfrak{b}_1$ 
relative to which $\iota$ is a linear Poisson map. In analogy to the previous 
case, this Poisson structure on $\mathfrak{b}_1$ will be said to be the Banach
Lie-Poisson structure {\bfi induced\/} by the mapping $\iota$. Proposition
\ref{example of induced structure} in \S \ref{section: preduals} gives an 
example of an induced Poisson structure.

\medskip

\begin{proposition}
\label{stronger condition}
Let $\mathfrak{b}_1$ be a Banach space,
$(\mathfrak{b}_2, \{\cdot, \cdot \}_2)$ be a Banach Lie-Poisson
space, and $\iota : \mathfrak{b}_1 \rightarrow \mathfrak{b}_2$ be
an injective continuous linear map with closed range. Then the
equality
\begin{equation}
\label{cohomological condition}
{\operatorname{ad}^2}^\ast_{\mathfrak{b}_2^\ast}
\iota(\mathfrak{b}_1) =  \iota(\mathfrak{b}_1)
\end{equation}
implies that $\ker \iota^\ast$ is an ideal in the Banach Lie algebra
$\mathfrak{b}_2^\ast$ and thus the map $\iota: \mathfrak{b}_1 \rightarrow 
\mathfrak{b}_2$ induces a Banach Lie Poisson structure on $\mathfrak{b}_1$.
\end{proposition}

\noindent \textbf{Proof}. To show that \eqref{cohomological condition}
implies that $\ker \iota^\ast$ is an ideal, we prove first the following 
equality:
\[
\ker \iota^\ast = \left\{ x \in \mathfrak{b}_2^\ast \mid
{\operatorname{ad}^2}^\ast_x \iota(\mathfrak{b}_1) =0 \right\}.
\]
To see this, let $x, y \in \mathfrak{b}_2^\ast$ and note that
\[
\left\langle y, {\operatorname{ad}^2}^\ast_x \iota(\mathfrak{b}_1)
\right\rangle = - \left\langle x, {\operatorname{ad}^2}^\ast_y
\iota(\mathfrak{b}_1) \right\rangle.
\]
Thus, ${\operatorname{ad}^2}^\ast_x \iota(\mathfrak{b}_1)=0$ if
and only if
\[
\left\langle x, {\operatorname{ad}^2}^\ast_y \iota(\mathfrak{b}_1)
\right\rangle = 0 \quad \text{for~all}\quad  y \in \mathfrak{b}_2^\ast,
\]
which, by condition \eqref{cohomological condition}, is equivalent
to $0 = \langle x, \iota(\mathfrak{b}_1) \rangle = \langle
\iota^\ast x, \mathfrak{b}_1 \rangle$, that is, $\iota^\ast x = 0$.

Next we prove that
\[
\left\{ x \in \mathfrak{b}_2^\ast \mid
{\operatorname{ad}^2}^\ast_x \iota(\mathfrak{b}_1) =0 \right\}
\]
is an ideal. Indeed if $x$ is in this subspace and $y, z \in
\mathfrak{b}_2^\ast$ are arbitrary, then
\[
\left\langle z, {\operatorname{ad}^2}^\ast_{[x,y]}
\iota(\mathfrak{b}_1) \right\rangle = \left\langle z,
{\operatorname{ad}^2}^\ast_y {\operatorname{ad}^2}^\ast_x
\iota(\mathfrak{b}_1) \right\rangle -\left\langle z,
{\operatorname{ad}^2}^\ast_x {\operatorname{ad}^2}^\ast_y
\iota(\mathfrak{b}_1) \right\rangle = 0
\]
because in the second term ${\operatorname{ad}^2}^\ast_y
\iota(\mathfrak{b}_1) \subset \iota(\mathfrak{b}_1)$ by condition
\eqref{cohomological condition} and the element $x$ satisfies
${\operatorname{ad}^2}^\ast_x \iota(\mathfrak{b}_1)=0$.

These two steps show that $\ker \iota^\ast$ is an ideal in the
Banach Lie algebra $\mathfrak{b}_2^\ast$.

Therefore, if \eqref{cohomological condition} holds, by
Proposition \ref{injective LP map}, the space $\mathfrak{b}_1$
carries a unique Banach Lie-Poisson structure such that $\iota$ is
a linear Poisson map. $\quad \blacksquare$

\medskip

The previous two propositions give an algebraic characterization
of linear Poisson maps $\phi: \mathfrak{b}_1 \rightarrow
\mathfrak{b}_2$ between Banach Lie-Poisson spaces analogous to that
from linear algebra. We summarize this  in the following theorem.

\begin{theorem}
The linear continuous map $\phi:\b_1\to\b_2$ between the Banach
Lie-Poisson spaces $\b_1$ and $\b_2$, such that $\phi(\b_1)$ is
Banach subspace in $\b_2$, is a linear Poisson map if and only if
it has a decomposition $\phi=\iota \circ\pi$, where
\begin{enumerate}[{\rm (i)}]
\item $\pi:\b_1\to\b$ is a linear continuous surjective map
of Banach spaces such that $\im\pi^*\subset \b_1^*$ is closed with
respect to Lie bracket of $\b_1^*$;
\item $\iota:\b\to\b_2$ is a continuous injective linear map of
Banach spaces with closed range such that $\ker \iota^*$ is an
ideal in the Banach Lie algebra $\b_2^*$.
\end{enumerate}
\end{theorem}

Due to the isomorphism between the category $\mathfrak{B}$ of Banach
Lie-Poisson spaces and the category $\mathfrak{L}_0$ of Banach Lie
algebras admitting a predual, properties of objects in 
one category can be characterized by properties in the other. We
will be interested now in finding the properties of an object in
$\mathfrak{B}$ that correspond to the condition of being an ideal or a
subalgebra of the related object in $\mathfrak{L}_0$. To do this
we shall use Propositions
\ref{surjective LP map} and \ref{injective LP map}.

\begin{proposition}
\label{correspondence}
Let $\mathfrak{b} \in \operatorname{Ob}(\mathfrak{B})$ be a Banach Lie-Poisson 
space and let $\mathfrak{g} \in \operatorname{Ob}(\mathfrak{L}_0)$ be a Banach 
Lie algebra such that $\mathfrak{b}^\ast = \mathfrak{g}$. Then:
\begin{enumerate}
\item[{\rm (i)}] There exists a bijective correspondence between the coinduced 
Banach Lie-Poisson structures from $\mathfrak{b}$ and the Banach Lie 
subalgebras of $\mathfrak{g}$. If the surjective continuous linear map $\pi: 
\mathfrak{b} \rightarrow \mathfrak{c}$ coinduces a Banach Lie-Poisson structure
on $\mathfrak{c}$, the Banach Lie subalgebra of $\mathfrak{g}$ given by this 
correspondence is $\pi^\ast(\mathfrak{c}^\ast)$. 

Conversely, if $\mathfrak{k} 
\subset \mathfrak{g}$ is a Banach Lie subalgebra then the Banach Lie-Poisson 
space given by this correspondence is $\mathfrak{b}/\mathfrak{k}^\circ$, where 
$\mathfrak{k}^\circ$ is the annihilator of $\mathfrak{k}$ in $\mathfrak{b}$, 
and $\pi: \mathfrak{b} \rightarrow  \mathfrak{b}/\mathfrak{k}^\circ$ is the 
quotient projection.
\item[{\rm (ii)}] There exists a bijective correspondence between the induced 
Banach Lie-Poisson structures in $\mathfrak{b}$ (i.e., the Banach Lie-Poisson 
subspaces of $\mathfrak{b}$) and the Banach ideals of $\mathfrak{g}$. If the 
injection $\iota :\mathfrak{c} \rightarrow \mathfrak{b}$ with closed range 
induces a Banach Lie-Poisson structure on $\mathfrak{c}$, then the ideal in 
$\mathfrak{g}$ given by this correspondence is $\ker \iota^\ast$. 

Conversely, 
if $\mathfrak{i} \subset \mathfrak{g}$ is a Banach ideal, then the Banach Lie-
Poisson subspace of $\mathfrak{b}$ given by this correspondence is
$\mathfrak{i}
^\circ$, where $\mathfrak{i}^\circ$ is the annihilator of $\mathfrak{i}$ in 
$\mathfrak{b}$ and $\iota: \mathfrak{i}^\circ \rightarrow \mathfrak{b}$ is the 
inclusion.
\end{enumerate}
\end{proposition}

\noindent \textbf{Proof}. (i) If the surjective continuous linear map $\pi: 
\mathfrak{b} \rightarrow \mathfrak{c}$ coinduces a Banach Lie-Poisson structure
on $\mathfrak{c}$, Proposition \ref{surjective LP map} states that $\pi^\ast
(\mathfrak{c}^\ast)$ is a Banach Lie subalgebra of $\mathfrak{g}$. Conversely, 
if $\mathfrak{k} \subset \mathfrak{g}$ is a Banach Lie subalgebra then $\pi: 
\mathfrak{b} \rightarrow \mathfrak{b}/\mathfrak{k}^\circ$ is a surjective 
continuous linear map of Banach spaces. Consider the dual map $\pi^\ast: 
[\mathfrak{b}/\mathfrak{k}^\circ]^\ast \rightarrow \mathfrak{b}^\ast$. Since 
$\mathfrak{b}^\ast = \mathfrak{g}$ and since $[\mathfrak{b}/\mathfrak{k}^\circ]
^\ast \cong \mathfrak{k}^{\circ \circ} = \mathfrak{k}$, it follows that 
$\operatorname{im}\pi^\ast$ is a Banach Lie subalgebra of $\mathfrak{g}$. 
Therefore, by Proposition \ref{surjective LP map}, there is a unique coinduced 
Banach Lie-Poisson structure on $\mathfrak{b}/\mathfrak{k}^\circ$.

(ii) This is a direct consequence of Proposition \ref{injective LP map}. $\quad
\blacksquare$

\medskip

Consider two Banach Lie-Poisson spaces $(\mathfrak{b}_1, \{\cdot,\cdot \}_1)$ and 
$(\mathfrak{b}_1, \{\cdot, \cdot\}_2)$. According to Theorem \ref{product theorem}, the 
product $(\mathfrak{b}_1 \times \mathfrak{b}_2, \{\cdot, \cdot\}_{12})$ is a Banach 
Poisson manifold. The Banach space isomorphism $(\mathfrak{b}_1 \times
\mathfrak{b}_2)^\ast \cong \mathfrak{b}_1^\ast \times \mathfrak{b}_2^\ast$ and formula 
\eqref{product bracket} show that $(\mathfrak{b}_1 \times \mathfrak{b}_2)
^\ast$ is closed under the product Poisson bracket $\{\cdot, \cdot\}_{12}$ which proves 
that $(\mathfrak{b}_1 \times \mathfrak{b}_2, \{\cdot,\cdot\}_{12})$ is also a Banach 
Lie-Poisson space. 

As opposed to the general case of Poisson manifolds, the inclusions $i_k:
\mathfrak
{b}_k \rightarrow 
\mathfrak{b}_1 \times \mathfrak{b}_2$, $k = 1,2$, defined by
$i_1(b_1) := (b_1, 0)$ and $i_2(b_2) := (0, b_2)$ are Poisson maps. Indeed, by 
\eqref{product bracket}, we get
\[
(i_1^\ast \{f, g\}_{12})(b_1) = \{f, g\}_{12}(b_1, 0) = 
\{f_0, g_0\}_1(b_1) + \{f_{b_1}, g_{b_1}\}_2(0) = 
\{i_1^\ast f, i_1^\ast g\}_{1}(b_1), 
\] 
where the term $\{f_{b_1}, g_{b_1}\}_2(0)$ vanishes because
in this case the Poisson bracket is linear. The proof for 
$i_2$ is similar.

Regarding the product construction, the following question arises naturally. 
When does a Banach Lie-Poisson space 
$(\mathfrak{b}, \{\cdot, \cdot\})$ allow a decomposition as a product of two Banach 
Lie-Poisson spaces $(\mathfrak{b}_1, \{\cdot,\cdot\}_1)$ and 
$(\mathfrak{b}_2, \{\cdot, \cdot\}_2)$?

In the category of Banach spaces this means that there is a splitting, i.e., 
$\mathfrak{b} = \mathfrak{b}_1 \oplus \mathfrak{b}_2$, for two Banach subspaces
$\mathfrak{b}_1$ and $\mathfrak{b}_2$, which is equivalent to $\mathfrak{b} 
\cong \mathfrak{b}_1 \times \mathfrak{b}_2$. In view of the previous properties
of the product of two Banach Lie-Poisson spaces, this suggests the following 
definition.

\begin{definition}
\label{splitting definition}
Let $(\mathfrak{b}, \{\,,\})$ be a Banach Lie-Poisson space. The splitting 
$\mathfrak{b} = \mathfrak{b}_1 \oplus \mathfrak{b}_2$ into two Banach subspaces
$\mathfrak{b}_1$ and $\mathfrak{b}_2$ is called a {\bfi Poisson splitting\/} if
\begin{enumerate}
\item[{\rm (i)}] $\mathfrak{b}_1$ and $\mathfrak{b}_2$ are Banach Lie Poisson spaces 
whose brackets shall be denoted by $\{\,,\}_1$ and $\{\,,\}_2$ respectively;
\item[{\rm (ii)}] the projections $\pi_k: \mathfrak{b} \rightarrow \mathfrak{b}_k$ 
and the inclusions $i_k: \mathfrak{b}_k \rightarrow \mathfrak{b}$, $k = 1,2$, 
consistent with the above splitting, are Poisson maps;
\item[{\rm (iii)}] if $f \in \pi_1^\ast (C^\infty(P_1))$ and 
$g \in \pi_2^\ast (C^\infty(P_2))$, then $\{f, g\} = 0$. 
\end{enumerate}
\end{definition}

The following proposition gives equivalent conditions for the existence of 
Poisson splittings.

\begin{proposition}
\label{splitting proposition}
    The following conditions are equivalent:
\begin{enumerate}
\item[{\rm (i)}] the Banach Lie-Poisson space 
$(\mathfrak{b}, \{\cdot, \cdot\})$ admits a Poisson splitting into the 
two Banach Lie-Poisson subspaces $(\mathfrak{b}_1, \{\cdot,\cdot\}_1)$ 
and $(\mathfrak{b}_2, \{\cdot, \cdot\}_2)$;
\item[{\rm (ii)}] the Banach Lie-Poisson space $(\mathfrak{b}, \{\cdot,\cdot\})$ is 
isomorphic to the product Banach Lie-Poisson space
$(\mathfrak{b}_1 \times \mathfrak{b}_2, \{\cdot,\cdot\}_{12})$;
\item[{\rm (iii)}] the components $\mathfrak{b}_1^\ast$ and $\mathfrak{b}_2^\ast$ of 
the dual splitting $\mathfrak{b}^\ast = \mathfrak{b}_1^\ast \oplus \mathfrak{b}
_2^\ast$ are
ideals of the Banach Lie algebra $\mathfrak{b}^\ast$, where one identifies 
$\mathfrak{b}_1^\ast$ and  $\mathfrak{b}_2^\ast$ with the annihilators of 
$\mathfrak{b}_2$ and $\mathfrak{b}_1$ in $\mathfrak{b}^\ast$ respectively.
\end{enumerate}
\end{proposition} 

\noindent \textbf{Proof}. The equivalence of (i) and (ii) is a direct 
consequence of Theorem \ref{product theorem} and the subsequent comments. 
Conditions (i) and (iii) are equivalent by applying Propositions 
\ref{surjective LP map} and \ref{injective LP map}. 
$\quad \blacksquare$
\medskip

In this section we established an equivalence between the category on Banach 
Lie algebras admitting  a predual and the category of Banach Lie-Poisson 
spaces. The statements proved above give examples how this equivalence can be 
used in the study of these two categories. For example, the simplicity of a 
Banach Lie algebra from the category
$\mathfrak{L}_0$ is equivalent to the 
non existence of Banach Lie-Poisson subspaces of its predual from
the category 
    $\mathfrak{B}$.

\section{Preduals of $W^\ast$-algebras as Banach Lie-Poisson
spaces}
\label{section: preduals}

In this section we shall consider the important class of Banach
Lie-Poisson spaces related to the category of $W^*$-algebras.

Recall that a $W^*$-algebra is a $C^*$-algebra $\mathfrak{m}$
which posses a predual Banach space $\mathfrak{m}_\ast$, i.e.
$\m=(\mathfrak{m}_\ast)^*$; this predual is unique (Sakai [1971]). Since
$\mathfrak{m}^\ast=(\mathfrak{m}_\ast)^{**}$, the predual Banach
space $\mathfrak{m}_\ast$ canonically embeds into the Banach space
$\mathfrak{m}^\ast$ dual to $\mathfrak{m}$. Thus we shall always
think of $\mathfrak{m}_\ast$ as a Banach subspace of
$\mathfrak{m}^\ast$. The existence of $\mathfrak{m}_\ast$ allows
the introduction of the 
$\sigma(\mathfrak{m},\mathfrak{m}_\ast)$-topology on the 
$W^*$-algebra
    $\mathfrak{m}$; for simplicity we shall call it the $\sigma$-topology in
the sequel. Recall that a net $\{x_\alpha\}_{\alpha \in A} \subset
\mathfrak{m}$ converges to $x \in \mathfrak{m}$ in the $\sigma$-topology
if, by definition, $\lim_{\alpha \in A} \langle x_\alpha, b \rangle =
\langle x, b \rangle$ for all $b \in \mathfrak{m}_\ast$. The
$\sigma$-topology is Hausdorff. Alaoglu's theorem states that the
unit ball of $\mathfrak{m}$ is compact in the $\sigma$-topology.
One can characterize the predual space $\mathfrak{m}_\ast$ as the
subspace of $\mathfrak{m}^\ast$ consisting of all
$\sigma$-continuous linear functionals, see Sakai [1971]. A theorem of 
Diximier (see Sakai [1971], \S 1.13)
states that a positive linear functional $\nu\in\mathfrak{m}^\ast$
is $\sigma$-continuous if and only if it is normal, i.e. it
satisfies
\[
\langle \nu,\, l.u.b. x_\alpha \rangle =l.u.b. \langle \nu,\,
x_\alpha \rangle
\]
for every uniformly bounded increasing direct set $\{x_\alpha\}$
of positive elements in $\mathfrak{m}$. The normality is
determined by the ordering on $\mathfrak{m}$ only. So, the predual
space $\mathfrak{m}_\ast$ and thus the pairing
\[
\mathfrak{m}_\ast\times\mathfrak{m}\ni(\nu,x)\mapsto\langle
x,\nu\rangle:=x(\nu)\in\C
\]
are defined by the algebraic structure of $\mathfrak{m}$ in a unique
way.

\begin{theorem}
\label{w star algebra theorem}
Let $\mathfrak{m}$ be a
$W^*$-algebra and $\mathfrak{m}_\ast$ be the predual of
$\mathfrak{m}$. Then $\mathfrak{m}_\ast$ is a Banach Lie-Poisson
space with the Poisson bracket $\{f,g\}$ of $f,g\in
C^\infty(\mathfrak{m}_\ast)$ given by \eqref{general LP}. The
Hamiltonian vector field $X_f$ defined by the smooth function $f\in
C^\infty(\mathfrak{m}_\ast)$ is given by \eqref{general Hamiltonian vector field}. 
\end{theorem}
\noindent\textbf{Proof}. We shall prove the theorem by checking
the conditions of Theorem \ref{general theorem}. Since the
$W^*$-algebra $\mathfrak{m}$ is an associative Banach algebra we
can define the Lie bracket in $\mathfrak{m}$ as the commutator
$$[x,y]=xy-yx$$ of $x,y\in\mathfrak{m}$. Left and right multiplication
by $a\in\mathfrak{m}$ define uniformly and $\sigma$-continuous
maps
\[
 L_a: \mathfrak{m}\ni x\mapsto ax\in\mathfrak{m}
\]
\[R_a: \mathfrak{m} \ni x\mapsto xa\in\mathfrak{m},
\]
see Sakai [1971]. Let
$L_a^*:\mathfrak{m}^\ast\to\mathfrak{m}^\ast$ and
$R_a^*:\mathfrak{m}^\ast\to\mathfrak{m}^\ast$ denote the dual maps
of $L_a$ and $R_a$ respectively. If $v\in\mathfrak{m}_\ast$, then
$L_a^*(v)$ and $R_a^*(v)$ are $\sigma$-continuous functionals and
therefore, by the characterization of the predual
$\mathfrak{m}_\ast$  as the subspace of $\sigma$-continuous
functionals in $\mathfrak{m}^\ast$, it follows that $L_a^*(v),
R_a^*(v)\in\mathfrak{m}_\ast$.
One has $\ad_a=[a,\cdot]=L_a-R_a$
and thus, $\ad_a^*=L_a^*-R_a^*$. We conclude from the above that
$\mathfrak{m}$ is a Banach Lie algebra and
$\ad_a^*\mathfrak{m}_\ast\subset\mathfrak{m}_\ast$ for each
$a\in\mathfrak{m}$, which are the conditions of Theorem
\ref{general theorem}. $\quad \blacksquare$
\medskip

\begin{corollary}
\label{corollary to w star theorem}
Let $\mathfrak{a}$ be a $C^\ast$-algebra. Then its dual 
$\mathfrak{a}^\ast$ is a Banach Lie-Poisson space.
\end{corollary}

\noindent\textbf{Proof}. The bidual $\mathfrak{a}^{\ast \ast}$ is 
isomorphic to the universal enveloping von Neumann algebra of
$\mathfrak{a}$ and by the canonical inclusion $\mathfrak{a}
\hookrightarrow \mathfrak{a}^{\ast \ast}$ the $C^\ast$-algebra
$\mathfrak{a}$ can be considered as a $C^\ast$-subalgebra of
$\mathfrak{a}^{\ast \ast}$ (see Sakai [1971] \S17.1, or Takesaki [1979]).
Since $\mathfrak{a}^\ast$ is predual to $\mathfrak{a}^{\ast \ast}$,
Theorem \ref{w star algebra theorem} guarantees that it is a Banach
Lie-Poisson space. $\quad \blacksquare$

\medskip

Any $\sigma$-closed $C^\ast$-subalgebra $\mathfrak{n}\subset \mathfrak{m}$ 
has the predual given by $\mathfrak{m}_\ast/\mathfrak{n}^\circ$, where
$\mathfrak{n}^\circ$ is the annihilator of $\mathfrak{n}$ in
$\mathfrak{m}_\ast$. Thus $\mathfrak{n}$ is a  Banach Lie subalgebra of
$\mathfrak{m}$ admitting  a predual. By Proposition \ref{correspondence}
(i), the quotient map $\pi:  \mathfrak{m}_\ast \rightarrow
\mathfrak{m}_\ast/\mathfrak{n}^\circ$ coinduces a Lie-Poisson structure on
the quotient Banach space $\mathfrak{m}_\ast/\mathfrak{n}^\circ$.
Therefore, there is a bijective correspondence between
$W^\ast$-subalgebras of $\mathfrak{m}$ and a subclass of  Banach
Lie-Poisson spaces coinducedÊ from $\mathfrak{m}_\ast$. It would be
interesting to characterize this subclass in Poisson geometrical terms.

If $\mathfrak{n}$ is a hereditary subalgebra of $\mathfrak{m}$, then 
there exists a projector $p \in \mathfrak{m}$ such that $\mathfrak{n} =
\operatorname{im} P$ (see Murphy [1990]), where the map $P: \mathfrak{m}
\rightarrow \mathfrak{m}$ is defined by 
\begin{equation}
\label{projector}
P(x): = pxp.
\end{equation}
The map $P$ is a $\sigma$-continuous projector with $\|P \| = 1$. Thus 
its dual $P^\ast : \mathfrak{m}^\ast \rightarrow \mathfrak{m}^\ast$
preserves $\mathfrak{m}_\ast$. We therefore conclude that $P_\ast: =
P^\ast |_{\mathfrak{m}_\ast}: \mathfrak{m}_\ast \rightarrow
\mathfrak{m}_\ast$ is a projector with $\|P_\ast \| = 1$. Thus there is a
splitting $\mathfrak{m}_\ast = \im P_\ast \oplus \ker P_\ast$ which allows
one to canonically identify $\mathfrak{m}_\ast/ \mathfrak{n}^\circ$ with
$\im P_\ast$ since $\ker P_\ast = \mathfrak{n}^\circ$.
\begin{proposition}
Let $\mathfrak{n}$ be a hereditary $W^\ast$-subalgebra of $\mathfrak{m}$. 
Then the projector $P_\ast: \mathfrak{m}_\ast \rightarrow
\mathfrak{m}_\ast$ coinduces a Banach Lie-Poisson structure on $\im
P_\ast$.
\end{proposition} 

So, a necessary condition for $\mathfrak{n}$ to be a hereditary subalgebra 
is that the Lie-Poisson structure on $\mathfrak{m}_\ast$ coinduces one on
a Banach subspace of $\mathfrak{m}_\ast$.

\medskip

Let us recall that $\mathfrak{m}_\ast$ and $\mathfrak{m}^\ast$ have
natural Banach Lie-Poisson structures according to Theorem \ref{w star
algebra theorem} and Corollary \ref{corollary to w star theorem}
respectively.

\begin{proposition}
\label{example of induced structure}
Let $\mathfrak{m}_\ast$ be the predual of the $W^\ast$-algebra
$\mathfrak{m}$ and $\iota: \mathfrak{m}_\ast \hookrightarrow
\mathfrak{m}^\ast$ be the canonical inclusion. Then $\iota$ is an injective
linear Poisson map and the Poisson structure induced by it from
$\mathfrak{m}^\ast$ coincides with the original Lie-Poisson structure on
$\mathfrak{m}_\ast$.
\end{proposition}

\noindent\textbf{Proof}. Since $\|\iota(b)\| = \|b\|$ for $b \in
\mathfrak{m}_\ast$, the range of $\iota: \mathfrak{m}_\ast
\hookrightarrow \mathfrak{m}^\ast$ is closed in $\mathfrak{m}^\ast$. The dual map
$\iota^\ast: \mathfrak{m}^{\ast \ast } \rightarrow \mathfrak{m} =
(\mathfrak{m}_\ast)^\ast$ is a projection of the universal enveloping
$W^\ast$-algebra $\mathfrak{m}^{\ast \ast}$ onto $\mathfrak{m}$ of norm
one. One has the equality $\ker \iota^\ast = (\mathfrak{m}_\ast)^\circ$,
where $(\mathfrak{m}_\ast)^\circ$ is the annihilator of
$\mathfrak{m}_\ast$ in $\mathfrak{m}^{\ast \ast}$. In addition, $L_x^\ast
\mathfrak{m}_\ast \subset \mathfrak{m}_\ast$ and $R_x^\ast
\mathfrak{m}_\ast \subset \mathfrak{m}_\ast$ for any $x \in
\mathfrak{m}^{\ast \ast}$ (see Sakai [1971]). Thus, $\ker \iota^\ast$ is a
$\sigma(\mathfrak{m}^{\ast \ast}, \mathfrak{m}^\ast)$-closed ideal of
$\mathfrak{m}^{\ast \ast}$. Therefore, $\ker \iota^\ast$ is also an ideal in the
Banach Lie algebra structure of $\mathfrak{m}^{\ast \ast}$ defined by
$[x,y] = xy - yx$. Proposition \ref{injective LP map} implies that
$\iota$ induces a Banach Lie-Poisson structure on $\mathfrak{m}_\ast$.
Since $\mathfrak{m}^{\ast \ast}/\ker \iota^\ast$ is isomorphic to
$\mathfrak{m}$, this structure coincides with the original Banach
Lie-Poisson structure of $\mathfrak{m}_\ast$ defined by Theorem \ref{w
star algebra theorem}. \quad $\blacksquare$
\medskip

Consider now a $\sigma$-closed two sided ideal $\mathfrak{i} \subset \mathfrak
{m}$. Then it equals $\im P$, where $P$ is given by \eqref{projector} where $p$
is a central projector in $\mathfrak{m}$ (see Sakai [1971] \S 1.10).
The projector $P_\perp: \mathfrak{m} \rightarrow \mathfrak{m}$ defined by 
$P_\perp (x) := (1-p)x (1-p)$ is also a $\sigma$-continuous linear map with 
$\|P_\perp\| = 1$ and projects $\mathfrak{m}$ onto a $\sigma$-closed two sided 
ideal $\mathfrak{i}_\perp$. Since $P + P_\perp = I$, we have the splitting
\begin{equation}
\label{ideal splitting}
\mathfrak{m} = \mathfrak{i} \oplus \mathfrak{i}_\perp
\end{equation}
of $\mathfrak{m}$ into two sided ideals. The decomposition \eqref{ideal 
splitting} is also a splitting into ideals in the category of Banach Lie 
algebras. By Proposition \ref{splitting proposition}, the direct sum \eqref
{ideal splitting} induces a Poisson splitting
\begin{equation}
\label{dual ideal splitting}
\mathfrak{m}_\ast = \mathfrak{i}^\circ \oplus 
\mathfrak{i}_\perp ^\circ, 
\end{equation} 
where $\mathfrak{i}^\circ$ and $\mathfrak{i}_\perp ^\circ$ are the annihilators
in $\mathfrak{m}_\ast$ of $\mathfrak{i}$ and $\mathfrak{i}_\perp$ 
respectively.

As a special case, one can consider the universal enveloping
$W^\ast$-algebra $\mathfrak{m}^{\ast \ast}$ of the $W^\ast$-algebra
$\mathfrak{m}$ with predual $\mathfrak{m}_\ast$. Then $\mathfrak{m}^\ast$ is the
predual to $\mathfrak{m}^{\ast \ast}$ and $\mathfrak{m}_\ast
\subset \mathfrak{m}^\ast$ is a $L^\ast _{\mathfrak{m}^{\ast \ast}}$ and 
$R^\ast _{\mathfrak{m}^{\ast \ast}}$ invariant Banach subspace. In this
case, the splitting \eqref{dual ideal splitting} gives the Poisson
splitting 
\[
\mathfrak{m}^\ast = \mathfrak{m}_\ast \oplus \mathfrak{m}_\ast ^\perp
\]
of $\mathfrak{m}^\ast$ into the normal and singular functionals (see
Takesaki [1979] for this terminology).

\medskip

In order to illustrate Theorem \ref{w star algebra theorem} let us
take a complex Hilbert space \M. By $\one$, $\two$, and $\bounded$
we shall denote the involutive Banach algebras of the trace class
operators, the Hilbert-Schmidt operators, and the bounded
operators on \M respectively. Recall that $\one$ and $\two$ are
self adjoint ideals in $\bounded$. Let
$\mathcal{K(M)}\subset\bounded$ denote the ideal of all compact
operators on \M. Then
\begin{equation}
\label{inclusions} \one \subset \two \subset \compact\subset
\bounded
\end{equation}
and the following remarkable dualities hold (see e.g. Murphy [1990]):
\begin{equation}
\label{dualities}
\mathcal{K(M)}^*\cong\one, \quad \two^*\cong\two, \quad \textrm{
and } \quad \one^*\cong\bounded.
\end{equation}
These are implemented by the strongly non-degenerate pairing
\begin{equation}
\label{eq:tr}\langle x,\rho \rangle= \tr(x\rho)
\end{equation}
where $x\in\one$, $\rho \in\mathcal{K(M)}$ for the first
isomorphism, $\rho,x\in\two$ for the second isomorphism and
$x\in\bounded$, $\rho \in\one$ for the third isomorphism. The isomorphism
$\one^*\cong\bounded$ gives the crucial example of the
$W^*$-algebra of bounded operators on the complex Hilbert space \M. So, we
recover the result of Bona [2000] as a corollary of Theorem
\ref{w star algebra theorem}.
\begin{corollary}
\label{LP structure on bounded operators}
The Banach space $\one$ of trace class operators on the Hilbert
space \M is a Banach Lie-Poisson space relative to the Poisson
bracket given by
\begin{equation}
\label{eq:pb}
\{f,g\}(\rho)=\tr([Df(\rho),Dg(\rho)]\rho)
\end{equation}
where $\rho\in\one$ and the bracket $[Df(\rho),Dg(\rho)]$ denotes
the commutator of the bounded operators $Df(\rho),
Dg(\rho)\in\bounded\cong\one^*$. The Hamiltonian vector fieldassociated to $f\in
C^\infty(\one)$ is given by
\begin{equation}
\label{eq:ham} 
X_f(\rho)=[Df(\rho),\rho].
\end{equation}
\end{corollary}
\noindent\textbf{Proof}. Formula \eqref{eq:pb} follows from
\eqref{general LP} by using \eqref{eq:tr} for the pairing between
$\one$ and $\bounded$ . In order to obtain \eqref{eq:ham} from
\eqref{general Hamiltonian vector field}, let us notice that
\begin{equation}
\langle y,-\ad_x^*\rho\rangle = -\langle[x,y],\rho\rangle=
-\tr([x,y]\rho)=\tr(y[x,\rho])=\langle y,[x,\rho]\rangle
\end{equation}
for $\rho\in\one$ and $x,y\in\bounded$. Thus
$-\ad^*_x\rho=[x,\rho]\in\one$, since $\one$ is an ideal in
$\bounded$. (We have identified here $\{\rho\}\times\one$ with the
tangent space $T_\rho\one$.) $\quad \blacksquare$
\medskip

The other two isomorphisms in \eqref{dualities} also give rise to
Banach Lie-Poisson spaces, but as a corollary of Theorem \ref{general
theorem}; Theorem \ref{w star algebra theorem} cannot be applied because
$\two$ and $\one$ are not $W^\ast$-algebras.

\begin{example}
\label{LP structure on two}
\normalfont
The Banach space $\two$ of Hilbert-Schmidt operators on the Hilbert
space $\mathcal{M}$ is a Banach Lie-Poisson space. Indeed,
we use the isomorphism $\two^*\cong\two$ given by the pairing
\eqref{eq:tr} and notice that $\two$ is a reflexive (that is,
$\two^{**}=\two$) Banach algebra. The formulas for the Poisson bracket
and for the Hamiltonian vector field are \eqref{eq:pb} and
\eqref{eq:ham} respectively with $\rho \in \two$. \quad $\blacklozenge$
\end{example}

\begin{example}
\label{LP structure on one}
\normalfont
The Banach space $\compact$ of compact operators on the Hilbert space
$\mathcal{M}$, as a predual of $\one$, is a Banach Lie-Poisson space. The
proof is identical to that of Corollary \ref{LP structure on bounded
operators}. The formulas  for the Poisson bracket and for the
Hamiltonian vector field are
\eqref{eq:pb} and \eqref{eq:ham} respectively with $\rho \in \compact$.
\quad $\blacklozenge$
\end{example}

\begin{example}
\label{l one l infinity}
\normalfont

Let $L^\infty(M, \mu)$ be the $W^\ast$-algebra of all essentially
bounded $\mu$-locally measurable functions on a localizable measure
space $M$. Then its predual is the Banach space $L^1(M, \mu)$ of all
$\mu$-integrable functions on $M$. Since $L^\infty(M, \mu)$
is commutative, the Banach Lie-Poisson structure on $L^1(M, \mu)$ is
trivial, that is, $\{f, g\} = 0$ for all $f, g \in C^\infty(L^1(M,
\mu))$. 

However, one can take the $W^\ast$-algebra tensor
product $L^\infty(M, \mu) \bar{\otimes} \mathfrak{m}$, where
$\mathfrak{m}$ is a $W^\ast$-algebra with predual $\mathfrak{m}_\ast$.
Then (see, e.g. Sakai [1971], \S 1.22)  $L^\infty(M, \mu) \bar{\otimes}
\mathfrak{m}$ is isomorphic with the Banach space $L^\infty(M, \mu,
\mathfrak{m})$ of all $\mathfrak{m}$-valued essentially bounded weakly
${}^\ast$
$\mu$-locally measurable functions on $M$. Moreover, $L^\infty(M, \mu,
\mathfrak{m})$ is a $W^\ast$-algebra under pointwise multiplication and
its predual is the Banach space $L^1(M, \mu, \mathfrak{m}_\ast)$ of all
$\mathfrak{m}_\ast$-valued Bochner $\mu$-integrable functions on $M$.
For details see Sakai [1971], Takesaki [1979], or Bourbaki [1959] \S2.6.
The duality  pairing between $b \in L^1(M, \mu,
\mathfrak{m}_\ast)$ and $x \in L^\infty (M, \mu, \mathfrak{m})$ is given
by 
\begin{equation}
\label{pairing with measure spaces}
\langle b, x \rangle = \int_M\langle b(m), x(m) \rangle
_{\mathfrak{m}} d\mu(m),
\end{equation}
where $\langle \cdot, \cdot \rangle _{\mathfrak{m}}$ is the duality
 pairing between $\mathfrak{m}_\ast$ and $\mathfrak{m}$. Thus,
by Theorem \ref{w star algebra theorem} and formula \eqref{general LP}
the Lie-Poisson  bracket of $f, g \in C^\infty(L^1(M, \mu,
\mathfrak{m}_\ast))$ is given by
\begin{equation}
\label{LP with measure spaces}
\{f, g\}(b) = \int_M \left\langle b(m), \left[\frac{\delta f}{\delta
b}(m), \frac{\delta g}{\delta b}(m)   \right]
\right\rangle _{\mathfrak{m}} d\mu(m),
\end{equation}
where $\delta f/ \delta b, \delta g/ \delta b \in L^\infty (M, \mu,
\mathfrak{m})$ are the representatives via the paring \eqref{pairing with
measure spaces} of the Fr\'echet derivatives $Df(b)$ and $Dg(b) \in
L^1(M, \mu, \mathfrak{m}_\ast)^\ast$ respectively.

Applying to $L^1(M, \mu, \mathfrak{m}_\ast)$ the quantum reduction
procedure (see Section \ref{section: quantum reduction}), one obtains the
Banach Lie-Poisson space
$L^1(M,
\mu,
\mathfrak{g}_\ast)$, where $\mathfrak{g}_\ast$ is the predual space of
the reduced Banach Lie algebra $\mathfrak{g} = R^\ast (\mathfrak{m})$.
In the finite dimensional case, for example when $\mathfrak{m} =
\mathfrak{gl}(N, \mathbb{C})$ and $M$ is a smooth manifold, we will
consider the Banach Lie algebra  $L^\infty(M, \mu, \mathfrak{g})$ as
the Lie algebra of the current group $C^\infty (M, G)$, where $G$ is a
Lie group with Lie algebra $\mathfrak{g}$ and the group structure on
$C^\infty(M, G)$ is defined by pointwise multiplication of maps. Usually
the Lie algebra of $C^\infty (M, G)$ is taken to be $C^\infty (M,
\mathfrak{g})$  (see, e.g. Kirillov [1993]); in our approach we shall
work with the
$L^\infty$ completion of this Lie algebra. For
$M = S^1$ one has the loop group case. So, we could consider the Banach
Lie-Poisson space $L^1(M, \mu, \mathfrak{g}^\ast)$ with the bracket
\begin{equation}
\label{LP finite dimensional on measure spaces}
\{f, g\}(\alpha) = \int_M C^i_{jk} \alpha_i(m) \frac{\delta f}{\delta
\alpha_j}(m) \frac{\delta f}{\delta \alpha_k}(m) d\mu(m)
\end{equation}
as one related to the current group. In order to clarify \eqref{LP
finite dimensional on measure spaces}, let us mention that, since
$\mathfrak{g}^{\ast \ast} = \mathfrak{g}$, we identified
$\mathfrak{g}_\ast$ with $\mathfrak{g}^\ast$. The scalar functions
$\alpha_1, \dots, \alpha_s$, where $s=\operatorname{dim}\mathfrak{g}$,
denote the components of $\alpha \in L^1(M, \mu, \mathfrak{g})$ in a
basis of $\mathfrak{g}^\ast$ dual to a given basis of $\mathfrak{g}$
relative to which the structure constants $C^i_{jk}$, $i, j, k = 1,
\dots , s$ are determined. \quad $\blacklozenge$
\end{example}

Let us now discuss the realifications $\mathfrak{m}_\mathbb{R}$ and
${\mathfrak{m}_\ast}_\mathbb{R}$ of the $W^\ast$-algebras $\mathfrak{m}$
and its predual $\mathfrak{m}_\ast$. As was mentioned in \S \ref{section:
banach lie-poisson spaces}, ${\mathfrak{m}_\ast}_\mathbb{R}$ has a real
Banach Lie-Poisson structure. For a fixed Hermitian element $\eta \in
\mathfrak{m}$ satisfying $\eta^2 = 1$, one defines the involutions
\begin{equation}
\label{involution one}
\sigma(b) = -\eta b^\ast \eta =: -\operatorname{Ad}^\ast_\eta (b^\ast),
\qquad b \in {\mathfrak{m}_\ast}_\mathbb{R}
\end{equation}
\begin{equation}
\label{involution two}
\sigma^\ast(x) = -\eta x^\ast \eta  =: -\operatorname{Ad}_\eta (x^\ast),
\qquad x \in \mathfrak{m}_\mathbb{R}
\end{equation}
in the sense of Example \ref{example following Proposition}, i.e., they
satisfy conditions (i), (ii), and (iii) given there. To check them, one
uses the defining identity for the conjugation in the predual
$\mathfrak{m}_\ast$: 
\begin{equation}
\label{definition of star in the predual}
\overline{\langle x, b\rangle} =
\langle x^\ast, b^\ast \rangle,
\end{equation}
where $x^\ast$ and $b^\ast$ are the conjugates of
$x \in \mathfrak{m}$ and $b \in \mathfrak{m}_\ast$ respectively. The
real Banach Lie algebra $\mathfrak{m}^\sigma: = \{x \in \mathfrak{m} \mid \sigma^\ast
x = x \} = \{x \in \mathfrak{m} \mid \eta x^\ast + x \eta = 0 \}$ has underlying
Banach Lie group
\[
U(\mathfrak{m}, \eta) := \{ g \in \mathfrak{m} \mid g^\ast \eta g = \eta \}
\]
consisting of the set of pseudounitary elements (see Bourbaki
[1972], Chapter 3, \S3.10, Proposition 37). For
$\eta = 1$ one  obtains the group of unitary elements of $\mathfrak{m}$. From the
considerations presented in Example \ref{example following
Proposition}, one can conclude that $\mathfrak{m}_\ast ^\sigma : = \{b \in
\mathfrak{m}_\ast \mid \sigma(b) = b\}$ has the real Banach Lie-Poisson
structure coinduced from $\mathfrak{m}_\mathbb{R}$ by the projector $R=(id
+ \sigma)/2$. This structure is given by \eqref{example real Poisson bracket} and is
$\operatorname{Ad}^\ast_{U(\mathfrak{m}, \eta)}$-invariant. In \S \ref{section:
symplectic leaves}, we will discuss the orbits of this action.

The above more general constructions are of course valid if one considers
the special case $\mathfrak{m} = \bounded$ and $\mathfrak{m}_\ast = \one$.
\medskip

As we have seen, Poisson geometry naturally arises in the  
the theory of operator algebras. The links between these theories 
established above show the importance of the fact that the category of
$W^\ast$-algebras can be considered as a subcategory of the category of
Banach Lie-Poisson spaces. Finally, let us mention that Poisson structures
that are fundamental for classical phase spaces appear in a natural way on
quantum phase spaces, i.e. duals to $C^\ast$-algebras.

\section{Quantum reduction} 
\label{section: quantum reduction}

Recall that $\one$ contains the subset of mixed states $\rho$ of the quantum
mechanical physical system, i.e., $\rho^\ast = \rho \geq 0$ and
$\operatorname{tr} \rho = \|\rho \|_1 = 1$. If the system under
consideration is an isolated quantum system, its dynamics is reversible
and is  described by the Liouville-von Neumann equation
\begin{equation}
\label{von Neumann equation}
\dot \rho = [H, \rho]
\end{equation}
which is a Hamiltonian equation on $(\one, \{\cdot , \cdot \})$ with
Hamiltonian $\operatorname{tr}(H\rho)$. For simplicity, let us assume that
$H^\ast = H \in \bounded$ is a given ($\rho$-independent) operator.
Therefore the Schr\"odinger flow $U(t) = e^{itH}$ is a Poisson flow on the
Banach Lie-Poisson space $(\one, \{\cdot, \cdot \})$.

Let us now apply a measurement operation to the system corresponding to
the discrete orthonormal decomposition of the unit
\begin{equation}
\label{decomposition of unity}
P_n P_m = \delta_{nm} P_n, \qquad \sum_{n=1}^\infty P_n = 1.
\end{equation}For example, this is the case when one measures the physical quantity
given by the operator $X= \sum_{n=1}^\infty x_n P_n$, with $x_n \in
\mathbb{R}$. Then, according to the well known von Neumann projection
postulate, the density operator $\rho$ of the state before measurement is
transformed by the measurement to the density operator $R(\rho)$ given by
\begin{equation}
\label{measurement}
R(\rho) := \sum_{n=1}^\infty P_n \rho P_n
\end{equation}

\begin{proposition}
\label{measurement properties}
The measurement operator $R: \one \rightarrow \one$ has the following
properties:
\begin{enumerate}
\item[{\rm (i)}] $R$ is a continuous norm one projector, i.e., $R^2 = R$ and
$\|R\|= 1$;
\item[{\rm (ii)}] it preserves the space of states, i.e., if $\rho^\ast =
\rho >0$, then $R(\rho)^\ast =  R(\rho^\ast)= R(\rho) >0$;\item[{\rm (iii)}] the
range $\operatorname{im} R^\ast$ of its dual $R^\ast
:\bounded \rightarrow \bounded$ is a Banach Lie subalgebra of $\bounded$.
\end{enumerate}
\end{proposition}

\noindent\textbf{Proof}. Using the natural pairing between $\one$ and
$\bounded$ given by the trace of the product, it follows that $R^\ast$ is
given by
\begin{equation}
\label{dual measurement}
R^\ast(X) := \sum_{n=1}^\infty P_n X P_n
\end{equation}
for $X \in \bounded$. Then, for $v \in \mathcal{M}$ one concludes
\[
\|R^\ast(X)v\|^2 = \sum_{n=1}^\infty \|P_nXP_nv\|^2 \leq 
\|X\|^2 \sum_{n=1}^\infty \|P_n v\|^2 = \|X\|^2 \|v\|^2
\]
which proves that $\|R^\ast\| \leq 1$. Since $\|R^\ast\| = \|{R^\ast}^2\|
\leq \|R^\ast\|^2$ it follows that $\|R^\ast \| = 1$. Now, using the
defining identity $\operatorname{Tr} R(\rho) X = \operatorname{Tr}
\rho R^\ast (X)$ of $R^\ast$ it follows that $\|R\|=1$. This proves (i).
Property (ii) follows directly from \eqref{dual measurement}. Finally, in
order to prove (iii) it is enough to remark that
\begin{equation}
\label{product of r stars}
R^\ast(X) R^\ast(Y) = R^\ast(R^\ast(X) R^\ast(Y)).
\qquad \blacksquare
\end{equation}

We conclude from Propositions \ref{measurement properties} and
\ref{surjective LP map} that the quantum measurement procedure gives a
Poisson projection $R: \one \rightarrow \operatorname{im}R$ of $\one$ on
the Banach subspace $\operatorname{im}R = \ker(1- R)$ endowed with the
Poisson bracket $\{\cdot, \cdot\}_{\operatorname{im}R}$ coinduced from
$\one$. Clearly, opposite to the $U$-procedure, i.e., the unitary time
evolution $U(t)$, $t \in \mathbb{R}$, the $R$-procedure is not reversible.
However, both the $U$-procedure and the $R$-procedure share an essential
common feature: they are linear Poisson maps.
\medskip

After this physical introduction, let us now come back to the case when the
Banach Lie-Poisson space is the predual space $\mathfrak{m}_\ast$ of a
$W^\ast$-algebra $\mathfrak{m}$. In the theory of quantum physical systems
(including statistical physics) the $W^\ast$-algebra is the algebra of
observables and the norm one  positive elements of
$\mathfrak{m}_\ast\subset\mathfrak{m}^\ast$ are the normal states of
the considered system, see e.g. Bratteli and Robinson [1979], [1981], or
Emch [1972].

The norm one map $E:\m\to\m$ which is idempotent ($E^2=E$)
and maps $\mathfrak{m}$ onto a $C^*$-subalgebra $\mathfrak{n}$ is called
{\bfi conditional expectation\/}. If $E$ is $\sigma$-continuous then
$\mathfrak n$ is a $W^*$-subalgebra of $\mathfrak{m}$. In that case, the
adjoint map $E^*:\mathfrak{m}^\ast \rightarrow \mathfrak{m}^\ast$ preserves
$\mathfrak{m}_* \subset \mathfrak{m}^*$ and maps $\mathfrak{m}_*$ onto
$\mathfrak{n}_*$. The conditional expectation is said to be 
{\bfi compatible with the state\/}
$\mu\in\mathfrak{m}_\ast$ if $E^*(\mu)=\mu$.

The concept of conditional expectation comes from probability
theory where it is very important in martingale theory. The definition of
the conditional expectation as the linear map
$E:\mathfrak{m} \rightarrow \mathfrak{m}$ on the $W^*$-algebra
$\mathfrak{m}$ with the properties mentioned above is the generalization of
the conditional expectation concept in non-commutative probability theory.
The role of conditional expectation in the theory of quantum measurement
theory and in quantum statistical physics and their remarkable
mathematical properties were elucidated in many remarkable publications
such as Takesaki [1972], Accardi, Frigerio, and Gorini [1984],
and Accardi and von Waldenfels [1988]. See Holevo [2001] for an extended
list of references to publications concerning conditional expectations.

Resuming, we see that the restriction $R:= 
E^*|_{\mathfrak{m}_\ast}:\mathfrak{m}_\ast\to\mathfrak{m}_\ast$
of the map dual to a conditional expectation $E:\m\to\m$
is a continuous projector. Since $\operatorname{im}R^\ast =
\operatorname{im}E = \mathfrak{n}$, the range of the projector $R^\ast :
\mathfrak{m} \rightarrow \mathfrak{m}$ is a Banach Lie subalgebra
$(\mathfrak{n}, [\cdot, \cdot])$ of $(\mathfrak{m}, [\cdot, \cdot])$. So,
like in the case of the measurement map \eqref{measurement}, one can apply
Proposition \ref{surjective LP map} in order to coinduce a Banach
Lie-Poisson structure on $\operatorname{im}R$.

Motivated by the above two examples, we introduce the following definition.

\begin{definition}
\label{quantum reduction definition}
A {\bfi quantum reduction map\/} is a continuous projector 
$R: \mathfrak{b} \rightarrow \mathfrak{b}$ on a Banach Lie-Poisson 
space $(\mathfrak{b}, \{\cdot, \cdot\})$ such that the range 
$\operatorname{im}R^\ast$ of the dual map $R^\ast: \mathfrak{b}^\ast 
\rightarrow \mathfrak{b}^\ast$ is a Banach Lie
subalgebra of $\mathfrak{b}^\ast$.
\end{definition}

 This immediately implies that  $R$
coinduces a Poisson structure on $\operatorname{im}R$ (see Proposition
\ref{surjective LP map}) and, in particular, $R: (\mathfrak{b}, \{\cdot,
\cdot \}) \rightarrow (\operatorname{im}R, \{\cdot,
\cdot\}_{\operatorname{im}R})$ is a Poisson map.
\medskip

Let us now give some important examples of quantum reduction.

\begin{example} 
\label{section 7, first example}
\normalfont
Every self-adjoint projector $p$ in the
$W^*$-algebra \m defines a uniformly and $\sigma$-continuous  projector
\begin{equation}
\m\ni x\mapsto P(x):=pxp\in\m 
\end{equation} 
of \m, see Sakai [1971] or Takesaki [1979]. Let 
$P^*:\mathfrak{m}^\ast\to\mathfrak{m}^\ast$ be the
projector dual to $P$, i.e.
$$
\langle P^*\mu,x\rangle=\langle \mu,P x\rangle
$$ 
for any $\mu\in \mathfrak{m}^\ast$ and $x\in\m$, where
$\langle\mu,x\rangle:=\mu(x)$. Since $P$ is $\sigma$-continuous, the
predual space $\mathfrak{m}_\ast\subset \mathfrak{m}^\ast$ is
preserved by $P^*$. Let $P_*$ be the restriction of $P^*$ to
$\mathfrak{m}_\ast$. The dual $(P_*)^*$ of the projector $P_*$ is
equal to $P$. The range $\im P$ of the projector $P$ is a
$W^*$-subalgebra of \m (see Sakai [1971]). Recalling that
$\ad^*_x\mathfrak{m}_\ast\subset \mathfrak{m}_\ast$ for $x\in\m$,
we see that $\ad_x\im P_*\subset \im P_*$, for $x\in\im P$.

Summarizing, we have proved the following.

\begin{proposition}
\label{example one proposition}
The projector $P_\ast :\mathfrak{m}_\ast\to\mathfrak{m}_\ast$ has the
following properties:
\begin{enumerate}[{\rm (i)}]
\item $\norm{P_*}=1$;
\item $\im(P_*)^*$ is a Banach-Lie algebra;
\item $\ad_x\im P_*\subset \im P_*$, for $x\in\im P$.
\end{enumerate}
Therefore $P_*:\mathfrak{m}_\ast\to\mathfrak{m}_\ast$ is a quantum
reduction map.
\end{proposition}

If $\m=\bounded$ and $\mathfrak{m}_\ast=\one$ the
projector $P_*:\one\to\one$ reduces the mixed state $\rho$ of the
quantum system to the state $p\rho p=P_*\rho$ localized on the
subspace $L^1(p\M)\subset\one$. In the quantum mechanical
formalism the projector $p:\M\to\M$ represents the so called elementary
observable ``proposition" (or ``question") which can have
only two alternative outcomes: ``yes" or ``no". The measurement of
the ``proposition" $p$ reduces the state $\rho$ to the state
$P_*\rho$ and the non-negative number $\tr(P_*\rho)$ is the
probability of the yes-answer. Since $P_*$ is a projector, the
repetition of the measurement does not change the state $P_*\rho$.
This is the mathematical expression of the von Neumann reproducing
postulate (von Neumann [1955]). \quad $\blacklozenge$
\end{example} 

\begin{example}
\label{section 7, second example}
\normalfont 

Let $\mathfrak{m}$ be a $W^\ast$-algebra and $\{p_\alpha\}_{\alpha\in I}\subset
\mathfrak{m}$, be a family of self-adjoint
mutually orthogonal projectors (i.e., $p_\alpha
p_\beta=\delta_{\alpha\beta}p_\alpha$, and $p^*_\alpha=p_\alpha$) such
that $\sum_{\alpha \in I} p_\alpha = 1$; the index set $I$ is not necessary
countable. Define the map $R^\ast : \mathfrak{m} \rightarrow \mathfrak{m}$ by
\begin{equation}
\label{r star in example two}
R^\ast (x) := \sum_{\alpha \in I} p_\alpha x p_\alpha
\end{equation}
for $x \in \mathfrak{m}$, where the summation is taken in the sense of the
$\sigma$-topology.

\begin{proposition}
\label{example two proposition}
The map $R^\ast : \mathfrak{m} \rightarrow \mathfrak{m}$ is a
$\sigma$-continuous linear projector with $\|R^\ast\| = 1$. Moreover,
$\operatorname{im}R^\ast$ is a $W^\ast$-subalgebra of $\mathfrak{m}$ and
hence a Banach Lie subalgebra of $(\mathfrak{m}, [\cdot, \cdot])$. 
Additionally, one has
\begin{equation}
\label{r star product identity}
R^\ast(R^\ast(x) R^\ast(y)) = R^\ast(x) R^\ast(y)
\quad \text~{and~} \quad
R^\ast(x^\ast) = (R^\ast(x))^\ast
\end{equation}
for all $x, y, \in \mathfrak{m}$.
\end{proposition}

\noindent\textbf{Proof}. We can always consider $\mathfrak{m}$ as a von
Neumann algebra of operators in the Hilbert space $\mathcal{M}$. If $v \in
\mathcal{M}$, then
\[
\|R^\ast(x)v\|^2 = \Big\|\sum_{\alpha \in I} p_\alpha x p_\alpha\Big
\|^2 = 
\sum_{\alpha \in I} \|p_\alpha x p_\alpha v\|^2 \leq \|x\|^2 \sum_{\alpha
\in I}\|p_\alpha v\|^2 = \|x\|^2 \|v\|^2
\]
which shows that $\|R^\ast (x)\| \leq 1$. From \eqref{r star in example
two} we have
\[
{R^\ast}^2(x) = \sum_{\beta \in I} p_\beta \left(\sum_{\alpha \in I}
p_\alpha x p_\alpha \right) p_\beta = 
\sum_{\beta \in I} \sum_{\alpha \in I} \delta_{\alpha \beta}
p_\alpha x p_\alpha  \delta_{\alpha \beta} =
\sum_{\alpha \in I} p_\alpha x p_\alpha = R^\ast (x);
\]
in this computation the $\sigma$-continuity of left and right
multiplication with an element $p_\beta$ was used. Thus ${R^\ast}^2 =
R^\ast$ and $\|R^\ast\| = 1$.

For any $b \in \mathfrak{m}_\ast$, there is an element $\rho \in \one$
such that $\langle x, b \rangle = \operatorname{tr}(x \rho)$. Thus
\[
\langle R^\ast(x), b \rangle = \operatorname{tr}(R^\ast(x)\rho) =
\sum_{\alpha \in I}\operatorname{tr}(p_\alpha x p_\alpha \rho) =
\sum_{\alpha \in I}\operatorname{tr}(x p_\alpha \rho p_\alpha) = 
\operatorname{tr}\left(x \sum_{\alpha \in I} p_\alpha x p_\alpha \right).
\]
We want to check that $x_i \stackrel{\sigma} \rightarrow x$ implies that
$R^\ast(x_i) \stackrel{\sigma} \rightarrow R^\ast(x)$. To do this,
substitute $x_i$ in the previous identity to get
\[
\langle R^\ast(x_i), b \rangle =
\operatorname{tr}\left(x_i \sum_{\alpha \in I} p_\alpha x p_\alpha \right)
\stackrel{\sigma} \longrightarrow
\operatorname{tr}\left(x \sum_{\alpha \in I} p_\alpha x p_\alpha \right)
= \operatorname{tr}(R^\ast(x) b) = \langle R^\ast(x), b \rangle
\]
for any $b \in \mathfrak{m}_\ast$. So $R^\ast$ is a $\sigma$-continuous
linear map. 

The defining formula for $R^\ast$ shows that for $x, y \in 
\mathfrak{m}$ one has $R^\ast(R^\ast(x)
R^\ast(y)) = R^\ast(x) R^\ast(y)$  and $R^\ast(x^\ast) = (R^\ast(x))^\ast$.
Thus $\operatorname{im}R^\ast$ is a
$W^\ast$-subalgebra of $\mathfrak{m}$ which implies that it is also a
Banach Lie subalgebra of $(\mathfrak{m}, [\cdot, \cdot])$. \quad
$\blacksquare$
\end{example} 

We conclude from Proposition \ref{example two proposition} that
$(R^\ast)^\ast: \mathfrak{m}^\ast
\rightarrow \mathfrak{m}^\ast$ preserves the predual subspace
$\mathfrak{m}_\ast \subset (\mathfrak{m}_\ast)^{\ast \ast} =
\mathfrak{m}^\ast$ and hence $R:= (R^\ast)^\ast |_{\mathfrak{m}_\ast}$ is
a quantum reduction.  Note that one has the splitting $\mathfrak{m} 
= \im R^\ast \oplus \ker R^\ast$. \quad $\blacklozenge$

\begin{example}
\label{section 7, third example}
\normalfont
Take the decomposition of the unit \eqref{decomposition of unity} and
define the map $R_-: \one \rightarrow \one$  by
\begin{equation}
\label{r minus}
R_-(\rho) := \sum_{n=1}^\infty \sum_{m=1}^n p_n \rho p_m = 
\sum_{n=1}^\infty p_n \rho q_n
\end{equation} 
where $q_n := \sum_{m=1}^n p_m$. It is clear that $R_-$ is a
linear projector on $\one$ whose range is the linear subspace of all
``lower triangular" trace class operators $\one_-$. From
\[
R_-(\rho)^\ast R_-(\rho) = \sum_{\ell = 1}^\infty \sum_{n = 1}^\infty
q_\ell \rho^\ast p_\ell p_n \rho q_n = \sum_{n=1}^\infty q_n \rho^\ast p_n
\rho q_n \leq \sum_{n=1} ^\infty \rho^\ast p_n \rho = \rho^\ast \rho
\]
we have
\[
\operatorname{tr} \sqrt{R(\rho)^\ast R(\rho)} \leq 
\operatorname{tr} \sqrt{\rho^\ast \rho}
\]
which shows that $\|R_-(\rho)\|_1 \leq \|\rho\|_1$. Thus, $R_-:\one
\rightarrow \one$ is a continuous projector with $\|R_-\| = 1$. So, the
dual map $R_-^\ast: \bounded \rightarrow \bounded$ defined by
\[
R^\ast_-(x): = \sum_{n=1}^\infty {q_n}x p_n
\]
also satisfies $\|R_-^\ast\|=1$ and projects $\bounded$ onto the ``upper
triangular" Banach Lie subalgebra $\bounded_+ \subset \bounded$. In this way,
$R_-: \one \rightarrow \one$ is a quantum reduction. Note that $R^\ast_-$ satisfies
\eqref{r star product identity}. In 
\S \ref{section: momentum maps and reduction} we will use the
quantum reduction $R_-$
in the description of the Toda lattice. \quad
$\blacklozenge$
\end{example} 

The discussion below uses the standard notion of Banach Lie group and its
associated Banach Lie algebra. Recall that a (real or complex) 
{\bfi Banach Lie group\/} is a (real or complex) smooth
Banach manifold $G$ with a group structure such that the multiplication
and inversion are smooth maps.  As in the
finite dimensional case, the tangent space at the identity, $T_e G$,
carries a Lie algebra structure which makes it isomorphic to the Lie
algebra of left invariant vector fields on $G$  endowed with the usual
bracket operation on vector fields. Due to the smoothness of the group
operations, this Lie bracket is a continuous bilinear map on the Banach
space $T_e G$. Thus $T_e G$ is a {\bfi Banach Lie algebra\/} which 
will be denoted, as customary, by $\mathfrak{g}$; it is called the 
Lie algebra of $G$.

Let $G(\mathfrak{m})$ be the group of invertible elements of a 
$W^\ast$-algebra $\mathfrak{m}$; 
it is an open subset (in the norm topology) of $\mathfrak{m}$ 
and is therefore a (real or complex) Banach Lie group whose Lie algebra 
is $\mathfrak{m}$ relative to the commutator bracket $[\cdot, \cdot]$
(Bourbaki [1972], Chapter III, \S3.9). Its exponential map is the usual
exponential function on $\mathfrak{m}$.

\begin{proposition}
\label{example two  proposition two}
Let $R:\mathfrak{m} \rightarrow \mathfrak{m}$ be a quantum reduction as
given in Definition 
\ref{quantum reduction definition} that also satisfies 
properties \eqref{r star product identity} and $\|R^\ast \| = 1$. Then
the set
$G(\mathfrak{m})
\cap
\im R^\ast =G(\im R^\ast)$ is a closed  Banach Lie subgroup of 
$G(\mathfrak{m})$ whose Lie algebra is the  Banach Lie subalgebra $\im R^\ast$. 
\end{proposition}

\noindent\textbf{Proof}. The equality in the statement is obvious. We next prove 
that $G(\im R^\ast)$ is a subgroup of $G(\mathfrak{m})$.

From \eqref{r star product identity} and the fact that $R^\ast$ maps the identity
to the identity, it follows that $G(\im R^\ast)$ is closed under multiplication 
in $\mathfrak{m}$ and that it contains the identity element. We shall prove now 
that if $R^\ast(x)$ is invertible, then its inverse is also an element of
$G(\im R^\ast)$. To see this, we assume without loss of generality that 
$\| 1- x\|<1$. Since $\|R^\ast\| = 1$ one has $\|R^\ast(1-x)\| <1$ and therefore
\begin{align*}
(R^\ast(x))^{-1} &= (1 - R^\ast(1-x))^{-1} 
= \sum_{k=0} ^ \infty \left[R^\ast (1-x)\right]^k \\
&= \sum_{k=0} ^ \infty R^\ast\left(\left[R^\ast (1-x)\right]^k\right) 
= R^\ast \left(\sum_{k=0} ^ \infty \left[R^\ast (1-x)\right]^k \right)
= R^\ast \left( (R^\ast(x))^{-1}\right).
\end{align*} 
Thus $G(\im R^\ast)$ is also closed under inversion and is therefore a subgroup 
of $G(\mathfrak{m})$. 

Since $\im R^\ast$ is closed in $\mathfrak{m}$, it follows that $G(\im R^\ast)$
is a closed subgroup of $G(\mathfrak{m})$. As the group of invertible elements
of the $W^\ast$-algebra $\im R^\ast$, $G(\im R^\ast)$ is a Lie group in its 
own right whose Lie algebra equals $\im R^\ast$. Because $\im R^\ast$ splits in
$\mathfrak{m}$ it follows that the inclusion of $G(\im R^\ast)$ into 
$G(\mathfrak{m})$ is an immersion. However, the topologies on $G(\im R^\ast)$
and $G(\mathfrak{m})$ are both induced by the norm topology of $\mathfrak{m}$ and
thus the inclusion $G(\im R^\ast) \hookrightarrow G(\mathfrak{m})$ is a 
homeomorphism onto its image which shows that this inclusion is an embedding and
hence $G(\im R^\ast)$ is a Lie subgroup of $G(\mathfrak{m})$.
\quad $\blacksquare$
\medskip

We shall return to this proposition in \S \ref{section: symplectic leaves}.
Note that both Examples \ref{section 7, second example} and
\ref{section 7, third example} satisfy the hypotheses of Proposition
\ref{example two  proposition two}.

\section{Symplectic leaves and coadjoint orbits}
\label{section: symplectic leaves}

A smooth map $f: M \rightarrow N$ between \textit{finite\/}
dimensional manifolds is called an immersion, if for every $m \in
M$ the derivative $T_mf : T_m M \rightarrow T_{f(m)}N$ is
injective. In infinite dimensions there are various notions
generalizing this concept.

\begin{definition}
\label{weak immersion}
A smooth map $f: M \rightarrow N$ between
Banach manifolds is called a
\begin{enumerate}
\item[{\rm (i)}] {\bfi immersion} if for every $m \in M$ the
tangent map $T_mf : T_mM \rightarrow T_{f(m)}N$ is injective with closed split range;
\item[{\rm (ii)}] {\bfi quasi immersion} if for every $m \in M$ the
tangent map $T_mf : T_mM \rightarrow T_{f(m)}N$ is injective with
closed range;
\item[{\rm (iii)}] {\bfi weak immersion} if for every $m \in M$ the
tangent map $T_mf : T_mM \rightarrow T_{f(m)}N$ is injective.
\end{enumerate}
\end{definition}

An immersion between Banach manifolds has the same properties as
an immersion between finite dimensional manifolds. For example, it
is characterized by the property that locally it is given by a map of the form $u
\mapsto (u, 0)$, where the model space of the chart
on $N$ necessarily splits. This is the concept widely used in the
literature; see e.g. Abraham, Marsden, and Ratiu [1988], or
Bourbaki [1967]. The notion of quasi immersion is modeled on the
concept of quasi regular submanifold introduced in Bourbaki
[1967]. Unfortunately, in the study of Banach Poisson manifolds
not even this weaker concept of quasi immersion is satisfactory
and one is forced to work with genuine weak immersions, as we
shall see in this section. 
\medskip

If $(P, \{\cdot, \cdot\}_P)$ is a Banach Poisson manifold (in the
sense of Definition \ref{Poisson manifold definition}), the vector
subspace $S_p := \{X_f(p) \mid f \in C^\infty(P) \}$ of $T_pP$ is
called the {\bfi characteristic subspace} at $p$. Note that $S_p$
is, in general, not a closed subspace of the Banach space $T_pP$.
The union $S:= \cup_{p\in P} S_p \subset TP$ is called the {\bfi
characteristic distribution} of the Poisson structure on $P$. Note
that even if  $S_p$ were closed and split in $T_pP$ for every
$p\in P$, $S$ would not necessarily be a subbundle of $TP$.
However, the characteristic distribution $S$ is always {\bfi
smooth}, in the sense that for every $v_p \in S_p \subset T_pP$
there is a locally defined smooth vector field (namely some $X_f$)
whose value at $p$ is $v_p$. Assume that the characteristic 
distribution is integrable. For finite dimensional manifolds this is
automatic by the Stefan-Sussmann theorem (see, e.g. Libermann and Marle
[1987], Appendix 3, Theorem 3.9) which, to our knowledge, is
not available in infinite dimensions.

Let $\mathcal{L}$ be a leaf of the characteristic distribution, that is,
\begin{itemize}
\item $\mathcal{L}$ is a connected smooth Banach manifold,
\item the inclusion $\iota: \mathcal{L} \hookrightarrow P$ is a
weak injective immersion,
\item $T_q \iota(T_q\mathcal{L}) = S_q$ for each $q \in
\mathcal{L}$,
\item if the inclusion $\iota':\mathcal{L}' \hookrightarrow P$ is
another weak injective immersion satisfying the three conditions
above and $\mathcal{L} \subset \mathcal{L}'$, then necessarily
$\mathcal{L}' = \mathcal{L}$, that is, $\mathcal{L}$ is maximal.
\end{itemize}
If we assume in addition that on the leaf $\mathcal{L}$
\begin{itemize}
\item there is a weak symplectic form $\omega_{\mathcal{L}}$  
consistent with the Poisson structure on $P$,
\end{itemize}
then $\mathcal{L}$ will be called a {\bfi symplectic leaf}.

In order to explain what this consistency means, consider from Definition 
\ref{Poisson manifold definition} the bundle map 
$\sharp : T^\ast P \rightarrow TP$ associated to the Poisson tensor $\varpi$ on $P$
and note that for each $p \in P$, the linear continuous map 
$\sharp_p :T^\ast _p P \rightarrow T_pP$ induces a bijective continuous map
$[\sharp_p] :T^\ast _p P / \ker \sharp_p \rightarrow S_p$. By definition,  
$\omega_{\mathcal{L}}$ is {\bfi consistent\/} with the Poisson structure on $P$ if  
\begin{equation}
\label{leaf form}
\omega_{\mathcal{L}}(q) (u_q, v_q) = 
\varpi(\iota(q))\left(([\sharp_{\iota(q)}]^{-1} \circ T_q\iota)(u_q), 
([\sharp_{\iota(q)}]^{-1} \circ T_q\iota)(v_q) \right).
\end{equation}
This shows that the weak symplectic form $\omega_{\mathcal{L}}$ consistent 
with the Poisson structure on $P$ is unique.  

For finite dimensional Poisson manifolds, it is known that all leaves are symplectic
(see Weinstein [1983]) and so the last assumption above is not necessary. In the
infinite dimensional case this question is open, even in the case of a Banach 
Lie group $G$ whose Lie algebra $\mathfrak{g}$ has a predual 
$\mathfrak{g}_\ast$ invariant under the coadjoint action. In this case, 
$\mathfrak{g}_\ast$ is a Banach Lie-Poisson space and we will
characterize a large class of points in $\mathfrak{g}_\ast$ whose
coadjoint orbits are all  weak symplectic manifolds. Their connected
components are therefore symplectic  leaves. Among these, we will also
identify a class of coadjoint orbits who are strong symplectic
manifolds.  

\begin{proposition}
\label{leaf bracket}
Let $\iota: (\mathcal{L}, \omega_{\mathcal{L}}) \hookrightarrow P$ be a symplectic leaf of the characteristic 
distribution of $P$. 
Then 
\begin{enumerate}
\item[{\rm (i)}] for any $f \in C^\infty (U)$,  $q \in \iota^{-1}( U) \cap \mathcal{L}$, 
where $U$ is an open subset of $P$, one has
\begin{equation}
\label{compatibility}
d\left((f \circ \iota)|_{ \iota^{-1}(U)} \right)(q) 
= \omega_\mathcal{L}(q)\left(\left(T_q\iota\right)^{-1} ( X_f(\iota(q))), \cdot \right).
\end{equation}

\item[{\rm (ii)}] the subspace $\iota^\ast(C^\infty(P))$ of $C^\infty(\mathcal{L})$ consisting of
functions that are obtained as restrictions of smooth functions from $P$
is a Poisson algebra relative to the bracket $\{\cdot, \cdot \}_{\mathcal{L}}$
given by
\begin{equation}
\label{bracket on leaf}
\{f \circ \iota, g \circ \iota \}_{\mathcal{L}}(q) 
: = \omega_{\mathcal{L}}(q) \left(\left(T_q\iota\right)^{-1} ( X_f(\iota(q))), 
\left(T_q\iota\right)^{-1} ( X_g(\iota(q))) \right);
\end{equation}

\item[{\rm (iii)}]  $\iota^\ast: C^\infty(P) \rightarrow \iota^\ast(C^\infty(P))$ is a 
homomorphism of Poisson algebras, that is,
\begin{equation}
\label{leaf homomorphism}
\{f \circ \iota, g \circ \iota \}_{\mathcal{L}} = \{f, g\}_P \circ \iota
\end{equation}
for any $f, g \in C^\infty(P)$.
\end{enumerate}
\end{proposition}

\noindent\textbf{Proof}. We begin with the proof of formula (\ref{compatibility}). For any $q \in 
\mathcal{L} \cap \iota^{-1}(U)$, 
$v_q \in T_q \mathcal{L}$, $f \in C^\infty(U)$, $U$ open in $P$, we have by (\ref{leaf form})
and the definition of $\sharp$
\begin{align*}
\omega_\mathcal{L}(q)\left(\left(T_q\iota\right)^{-1} ( X_f(\iota(q))), v_q \right) &= 
\varpi(\iota(q))\left(df(\iota(q)), 
([\sharp_{\iota(q)}]^{-1} \circ T_q\iota (v_q) \right) \\
&= \langle df(\iota(q)), T_q\iota (v_q) \rangle 
= d(f \circ \iota)(q)(v_q),
\end{align*} 
which proves (\ref{compatibility}). Now replace in the above computation $v_q$ by
$\left(T_q\iota\right)^{-1} (X_g(\iota(q)))$ and get
\begin{align*}
\omega_\mathcal{L}(q)\left(\left(T_q\iota\right)^{-1} ( X_f(\iota(q))), 
\left(T_q\iota\right)^{-1} (X_g(\iota(q))) \right) = 
\langle df(\iota(q)), X_g(\iota(q)) \rangle 
= \{f, g \}_P (\iota(q)),
\end{align*}
which, in view of (\ref{bracket on leaf}), proves (\ref{leaf homomorphism}). 
Finally, (\ref{leaf homomorphism}) shows that (\ref{bracket on leaf}) defines a Poisson
bracket on $\mathcal{L}$.
\quad $\blacksquare$
\medskip

Formula (\ref{compatibility}) is remarkable since it guarantees the existence of Hamiltonian vector fields on the weak symplectic manifold $\mathcal{L}$
for a large class of functions, namely those that are pull backs to the symplectic leaf.   

\medskip

We shall give below a class of Banach Poisson manifolds for 
which some of the symplectic leaves can be explicitly determined. 
In what follows $G$ denotes a (real or complex) Banach Lie group, 
$L_g$ and $R_g$ denote the diffeomorphisms of $G$ 
given by left and right translations by $g \in G$, and $\mathfrak{g}$
denotes the (left) Lie algebra of $G$.

\begin{theorem}
\label{quotient symplectic structure}
Let $G$ be a (real or complex) Banach Lie group with Lie algebra 
$\mathfrak{g}$. Assume that:
\begin{enumerate}
\item[{\rm (i)}] $\mathfrak{g}$ admits a predual $\mathfrak{g}_\ast$, that is, 
$\mathfrak{g}_\ast$ is a Banach space whose dual is 
$\mathfrak{g}$;
\item[{\rm (ii)}] the coadjoint action of $G$ on the dual 
$\mathfrak{g}^\ast$ leaves the predual $\mathfrak{g}_\ast$ invariant, that is, 
$\operatorname{Ad}_{g}^\ast (\mathfrak{g}_\ast) \subset
\mathfrak{g}_\ast$, for any $g\in G$;
\item[{\rm (iii)}] for a fixed $\rho \in \mathfrak{g}_\ast$ the coadjoint isotropy
subgroup $G_\rho := \{ g \in G \mid \operatorname{Ad}_{g}^\ast \rho = \rho\}$, which 
is closed in $G$, is a Lie subgroup of $G$ (in the sense that it is a submanifold
of $G$ and not just injectively immersed).
\end{enumerate}
Then the Lie algebra of $G_\rho$ equals
$\mathfrak{g}_\rho 
:= \{\xi \in \mathfrak{g} \mid
\operatorname{ad}^\ast _\xi \rho = 0\}$ and the quotient topological space
$G/G_\rho : = \{gG_\rho \mid g \in G\}$ admits a unique smooth (real or complex)
Banach manifold structure making the canonical projection $\pi:G \rightarrow
G/G_\rho$ a surjective submersion.  The manifold
$G/G_\rho$ is weakly symplectic relative to
the two form $\omega_\rho$ given by
\begin{equation}\label{KKSL}
\omega_\rho([g])(T_g\pi(T_e L_g \xi),T_g\pi(T_e L_g \eta))
:=\left\langle\rho,[\xi,\eta]\right\rangle,
\end{equation}
where $\xi, \eta \in\mathfrak{g},\,g \in G$, 
$[g]:=\pi(g)=g G_\rho$, and 
$\langle\cdot,  \cdot \rangle: \mathfrak{g}_\ast \times 
\mathfrak{g} \rightarrow \mathbb{R}$ (or $\mathbb{C}$) is the canonical pairing
between $\mathfrak{g}_\ast$ and $\mathfrak{g}$. Alternatively, this
form can be expressed as
\begin{equation}\label{KKSR}
\omega_\rho([g])(T_g\pi(T_e R_g \xi),T_g\pi(T_e R_g \eta))
:=\left\langle \operatorname{Ad}_{g^{-1}}^\ast \rho, 
[\xi,\eta]\right\rangle.
\end{equation}
The two form $\omega_\rho$ is invariant under the left action of
$G$ on $G/G_\rho$ given by $g \cdot [h] := [gh]$, for $g,h \in G$.
\end{theorem}

\noindent\textbf{Proof}. The subgroup $G_\rho$ is clearly 
closed. For Banach Lie groups it is no longer true that closed subgroups
are Lie subgroups (for a counterexample see Bourbaki [1972],
Chapter III, Exercise 8.2). However, as in the finite dimensional case, if $G_\rho$ 
is assumed to be a Lie subgroup of $G$, then (Bourbaki [1972], Chapter III, \S 6.4, 
Corollary 1) $\xi \in \mathfrak{g}$ is an element of 
the Lie algebra of $G_\rho$ if and only if $\exp t\xi \in G_\rho$ for all
$t \in \mathbb{R}$ (or $\mathbb{C}$ depending on whether $G$ is a real or complex 
Banach Lie group). Thus, since (see, e.g. Marsden and Ratiu [1994], Chapter 9)
\[
\frac{d}{dt} \operatorname{Ad}_{\exp t\xi}^\ast \rho = 
\operatorname{Ad}_{\exp t\xi}^\ast \operatorname{ad}^\ast _\xi \rho,
\]
it follows that 
\begin{align*}
\exp t\xi \in G_\rho &\Longleftrightarrow \operatorname{Ad}_{\exp t\xi}^\ast \rho 
= \rho \Longleftrightarrow 0 = \frac{d}{dt} \operatorname{Ad}_{\exp t\xi}^\ast \rho = 
\operatorname{Ad}_{\exp t\xi}^\ast \operatorname{ad}^\ast _\xi \rho \\
&\Longleftrightarrow \operatorname{ad}^\ast _\xi \rho = 0 \Longleftrightarrow 
\xi \in \mathfrak{g}_\rho,
\end{align*}
which shows that the Lie algebra of $G_\rho$ is $\mathfrak{g}_\rho$.

Since $G_\rho$ is assumed to be a Lie subgroup of $G$, the set
$G/G_\rho$ has a unique smooth manifold structure such that the canonical 
projection $\pi : G \rightarrow G/G_\rho$ is a submersion. The underlying 
manifold topology of $G/G_\rho$ is the quotient topology. The left action
$(g, [h]) \in G \times G/G_\rho\mapsto g\cdot [h]: = [gh]$
is smooth (see Bourbaki [1972], Chapter III, \S 1.6, Proposition 11, for 
a proof of these statements).

In what follows we shall need the following observation. Condition
(ii) implies that 
$\operatorname{ad}^\ast_\xi (\mathfrak{g}_\ast) \subset \mathfrak{g}_\ast$
for any $\xi \in \mathfrak{g}$.

The two-forms defined in formulas (\ref{KKSL}) and (\ref{KKSR})
are equal. Indeed, taking (\ref{KKSL}) as the definition but
applying it to tangent vectors of the form $T_g\pi(T_e R_g \xi)$,
$T_g\pi(T_e R_g \eta)$, we get
\begin{align*}
&\omega_\rho([g])(T_g\pi(T_e R_g \xi),T_g\pi(T_e R_g \eta)) =
\omega_\rho([g])\left(T_g\pi(T_e L_g (\operatorname{Ad}_{g^{-1}}\xi)),
T_g\pi(T_e L_g(\operatorname{Ad}_{g^{-1}}\xi))\right)\\ &
\qquad = \left\langle \rho, [\operatorname{Ad}_{g^{-1}}\xi, 
\operatorname{Ad}_{g^{-1}}\eta] \right\rangle =
\left\langle \rho, \operatorname{Ad}_{g^{-1}}[\xi,\eta] \right\rangle = 
\left\langle \operatorname{Ad}_{g^{-1}}^\ast \rho,
[\xi, \eta] \right\rangle
\end{align*}
which gives formula (\ref{KKSR}).

We shall prove now that the two-form (\ref{KKSR}) is well defined. Indeed, if $[g] =
[g']$ and $T_g\pi(T_e R_g \xi) = T_{g'}\pi(T_e R_{g'} \xi')$, then there is some $h
\in G_\rho$ such that $g' = gh$ and hence
\[
T_g\pi (T_e R_g \xi) = T_{g'}\pi(T_e R_{g'} \xi') 
= T_{gh}\pi(T_e R_{gh} \xi')= T_g (\pi \circ R_h)(T_e R_g \xi') = 
T_g\pi(T_e R_g \xi'),
\]
which means that $T_g\pi(T_e R_g (\xi- \xi')) = 0$. Due to the fact that the
fibers of $\pi$ are of the form $gG_\rho $, this is
equivalent to $T_e R_g (\xi- \xi') \in T_e L_g(\mathfrak{g}_\rho)$, that is,
$\xi - \xi' \in \operatorname{Ad}_g(\mathfrak{g}_\rho)$. Thus 
there is some $\zeta \in \mathfrak{g}_\rho$ such that 
$\xi' = \xi + \operatorname{Ad}_g \zeta$. Similarly, if 
$T_g\pi(T_e R_g \eta) = T_{g'}\pi(T_e R_{g'} \eta')$ there
is  some $\zeta' \in \mathfrak{g}_\rho$ such that 
$\eta' = \eta + \operatorname{Ad}_g \zeta'$. Therefore, since
$\operatorname{ad}_{\zeta} ^\ast \rho 
= \operatorname{ad}_{\zeta'} ^\ast \rho = 0$, it follows that
\begin{align*}
\omega_\rho([g'])&(T_{g'}\pi(T_e R_{g'} \xi'),
T_{g'}\pi(T_e R_{g'} \eta'))
=\left\langle \operatorname{Ad}_{{g'}^{-1}}^\ast \rho,[\xi',\eta']\right\rangle
\\ &\;
= \left\langle\operatorname{Ad}_{(gh)^{-1}}^\ast \rho,\, 
\left[\xi + \operatorname{Ad}_g \zeta,\, 
\eta + \operatorname{Ad}_g \zeta' \right]\right\rangle
\\ &\; 
= \left\langle\operatorname{Ad}_{g^{-1}}^\ast \rho,\, 
\left[\xi + \operatorname{Ad}_g \zeta,\, 
\eta + \operatorname{Ad}_g \zeta' \right]\right\rangle
\\ &\;
= \left\langle\operatorname{Ad}_{g^{-1}}^\ast \rho,\, 
[\xi, \eta ]\right\rangle
+ \left\langle\operatorname{Ad}_{g^{-1}}^\ast \rho,\, 
\left[ \operatorname{Ad}_g \zeta,\, 
\eta \right]\right\rangle \\
&\qquad \qquad 
+ \left\langle\operatorname{Ad}_{g^{-1}}^\ast \rho,\, 
\left[\xi, \operatorname{Ad}_g \zeta'\right]\right\rangle
+ \left\langle\operatorname{Ad}_{g^{-1}}^\ast \rho,\, 
\operatorname{Ad}_g[\zeta, \zeta']\right\rangle \\ & \;
= \left\langle\operatorname{Ad}_{g^{-1}}^\ast \rho,\, 
[\xi, \eta ]\right\rangle 
+ \left\langle \rho,\, 
\left[ \zeta,\, 
\operatorname{Ad}_{g^{-1}}\eta \right]\right\rangle
+ \left\langle \rho,\, 
\left[\operatorname{Ad}_{g^{-1}}\xi,\zeta'\right]\right\rangle
+ \left\langle \rho,\, [\zeta, \zeta']\right\rangle \\ & \;
= \left\langle\operatorname{Ad}_{g^{-1}}^\ast \rho,\, 
[\xi, \eta ]\right\rangle
+ \left\langle \operatorname{ad}_\zeta ^\ast \rho,\,  
\operatorname{Ad}_{g^{-1}}\eta \right\rangle
- \left\langle \operatorname{ad}_{\zeta'} ^\ast \rho,\,  
\operatorname{Ad}_{g^{-1}}\xi \right\rangle
+ \left\langle \operatorname{ad}_\zeta ^\ast \rho,\,  
\zeta' \right\rangle  \\ & \;
= \left\langle\operatorname{Ad}_{g^{-1}}^\ast \rho,\, 
[\xi, \eta ]\right\rangle
= \omega_\rho([g])(T_g\pi(T_e R_g \xi),T_g\pi(T_e R_g \eta)).
\end{align*}

The two-form $\omega_\rho$ is weakly non degenerate. Indeed if
$[g] \in G/G_\rho$ is given and if
\[
\omega_\rho([g])(T_g\pi(T_e L_g \xi), T_g\pi(T_e L_g \eta)) = 0
\] 
for all $T_g\pi(T_e L_g \eta)$
then, by (\ref{KKSL}), $\langle \operatorname{ad}_\xi ^\ast \rho, 
\eta \rangle = 0$ for all $\eta
\in \mathfrak{g}$ and thus $\xi\in \mathfrak{g}_\rho$ (since $\operatorname{ad}
^\ast _\xi \rho \in \mathfrak{g}_\ast$) which is
equivalent to $T_g\pi(T_e L_g \xi) = 0$.

We shall prove that $\omega_\rho$ is a smooth closed two form on
$G/G_\rho$ by showing that the smooth one form
on $G$ given by $\nu_\rho(g)(T_e L_g \xi): = 
-\langle \rho, \xi \rangle$
satisfies $d\nu_\rho = \pi^\ast\omega_\rho$. Since $\pi$ is a
surjective submersion this immediately implies that $\omega_\rho$
is smooth and closed. To compute the exterior derivative of
$\nu_\rho$, we denote by $X$, $Y$ the vector
fields on $G$ given by $X(g) = T_e L_g \xi$ and
$Y(g) = T_e L_g \eta$ and note that $\nu_\rho(X)(g) =
-\langle\rho, \xi \rangle$ is constant and $[X, Y ] (g) =
T_e L_g  [\xi, \eta]$. Therefore, by Cartan's formula,
\begin{align*}
&d\nu_\rho(g)(T_e L_g \xi, T_e L_g \eta) = 
d \nu_\rho(X, Y)(g) =
X[\nu_\rho (Y)](g) - Y[\nu_\rho
(\mathcal{X})](g) - \nu_\rho([X, Y])(g) \\
&\qquad = - \nu_\rho(g) (T_eL_g[\xi,\eta]) 
= \langle \rho, [\xi, \eta] \rangle 
= (\pi^\ast \omega_\rho)(g)(T_e L_g \xi, T_e L_g \eta).
\end{align*}

To show that $\omega_\rho$ is $G$-invariant, we note
that $\pi$ is $G$-equivariant, $\omega_\rho = \pi^\ast
\nu_\rho$, and that $\nu_\rho$ is $G$-invariant. $\quad
\blacksquare$
\medskip

We shall see a concrete example of a symplectic form 
$\omega_\rho$ that is weak and not strong after 
Example \ref{trace class operators example}.
\medskip

Next we study the coadjoint orbits of $G$ through points of $\mathfrak{g}_\ast$.

\begin{theorem}
\label{stratification theorem}
Let the Banach Lie group $G$ and the element $\rho \in \mathfrak{g}_\ast$ 
satisfy the hypotheses of Theorem \ref{quotient symplectic structure}.
Then the map
\begin{equation}
\label{leaf injection}
\iota:[g] \in G/G_\rho \mapsto \operatorname{Ad}^\ast _{g^{-1}} \rho 
\in \mathfrak{g}_\ast
\end{equation}
is an injective weak immersion of the quotient manifold 
$G/G_\rho$ into the predual space $\mathfrak{g}_\ast$.
Endow the coadjoint orbit $\mathcal{O} := \{\operatorname{Ad}^\ast _g \rho \mid
g \in G \}$ with the smooth 
manifold structure making $\iota$ into a diffeomorphism. The push forward 
$\iota_\ast(\omega_\rho)$ of the weak symplectic form 
$\omega_\rho \in \Omega^2(G/G_\rho)$ to $\mathcal{O}$ has 
the expression
\begin{equation}\label{KKS}
\omega_\mathcal{O}(\operatorname{Ad}^\ast_{g^{-1}} \rho)
\left(\operatorname{ad}^\ast_{\operatorname{Ad}_g \xi}
\operatorname{Ad}^\ast_{g^{-1}} \rho, 
\operatorname{ad}^\ast_{\operatorname{Ad}_g \eta}
\operatorname{Ad}^\ast_{g^{-1}} \rho \right)
= \langle\rho, [\xi,\eta]\rangle,
\end{equation}
for $g \in G$, $\xi, \eta \in \mathfrak{g}$, and $\rho \in \mathfrak{g}_\ast$. 
Relative to this symplectic form
the connected components of the coadjoint orbit 
$\mathcal{O}$ are 
symplectic leaves of the Banach Lie-Poisson space 
$\mathfrak{g}_\ast$.
\end{theorem}

\noindent\textbf{Proof}.   By Theorem \ref{functor} the predual $\mathfrak{g}_\ast$ is a (real
or holomorphic) Banach Lie-Poisson space whose (real or complex) Poisson
bracket is given by
\eqref{general LP}. For each $\rho \in \mathfrak{g}_\ast$ its
characteristic subspace is therefore given by 
$S_\rho = \{\operatorname{ad}^\ast_\xi \rho \mid \xi 
\in \mathfrak{g}\} \subset \mathfrak{g}_\ast$ since 
$\operatorname{ad}^\ast_\xi (\mathfrak{g}_\ast) \subset 
\mathfrak{g}_\ast$ for any $\xi \in \mathfrak{g}$.

The map $\iota: [g] \in G /G_\rho \mapsto \operatorname{Ad}^\ast_{g^{-1}}\rho 
\in \mathcal{O}$ is a bijection, so
we put on $\mathcal{O}$ the smooth Banach manifold
structure making this bijection into a diffeomorphism. Since the map $g
\in G \mapsto \operatorname{Ad}^\ast_{g^{-1}}\rho 
\in \mathfrak{g}_\ast$ is continuous, it thus follows that the inclusion
$\mathcal{O} \subset \mathfrak{g}_\ast$ is also continuous.

We shall prove now that the map $g \in G \mapsto 
\operatorname{Ad}^\ast_{g^{-1}}\rho
\in \mathfrak{g}_\ast$ is smooth. Indeed its derivative 
\begin{equation}
\label{derivative of diffeo implementing the orbit structure}
T_e L_g \xi \in  T_g G  \mapsto
- \operatorname{Ad}^\ast_{g^{-1}} \left(\operatorname{ad}^\ast_\xi
\rho \right)
= - \operatorname{ad}^\ast_{\operatorname{Ad}_g \xi} 
\operatorname{Ad}^\ast_{g^{-1}} \rho \in \mathfrak{g}_\ast
\end{equation}
is a continuous linear map from $T_gG$ to $\mathfrak{g}_\ast$, that is, 
the map $g \in G \mapsto \operatorname{Ad}^\ast_{g^{-1}} \rho
\in \mathfrak{g}_\ast$ is $C^1$. Inductively, it follows that it is
$C^\infty$. In addition, the range of the derivative at $g$ is the
characteristic subspace 
$S_{\operatorname{Ad}^\ast_{g^{-1}} \rho}$ at 
$\operatorname{Ad}^\ast_{g^{-1}} \rho$.

Since the map $g \in G \mapsto 
\operatorname{Ad}^\ast_{g^{-1}} \rho \in \mathfrak{g}_\ast$ is
$G_\rho$-invariant, it follows that $\iota: [g ]
\in G/G_\rho \mapsto \operatorname{Ad}^\ast_{g^{-1}} \rho 
\in \mathfrak{g}_\ast$ is
smooth and that the range of its derivative at $[g]$, given by
$T_{[g]} \iota:
T_g \pi(T_e L_g \xi) \in T_{[g]}(G/G_\rho) \mapsto
- \operatorname{ad}^\ast_{\operatorname{Ad}_g \xi} 
\operatorname{Ad}^\ast_{g^{-1}} \rho \in 
\mathfrak{g}_\ast$, 
equals $S_{\operatorname{Ad}^\ast_{g^{-1}} \rho}$.
The map $T_{[g]} \iota$ is injective. Indeed, if 
\[
0 = T_{[g]} \iota(T_g \pi(T_e L_g \xi)) =  
- \operatorname{Ad}^\ast_{g^{-1}} \left(\operatorname{ad}^\ast_\xi \rho \right) 
\] 
then $\xi \in \mathfrak{g}_\rho$ and hence $T_g \pi(T_e L_g \xi) = 0$. This shows
that $\iota$ is an injective weak immersion.

Let us endow the manifold $\mathcal{O} \subset \mathfrak{g}_\ast$ with the
push forward weak symplectic form $\omega_\mathcal{O}$ given by the
diffeomorphism $[g] \in G/G_\rho \mapsto 
\operatorname{Ad}^\ast_{g^{-1}} \rho \in
\mathcal{O} \subset \mathfrak{g}_\ast$. From the formula for
its derivative given by \eqref{derivative of diffeo implementing the orbit
structure} and (\ref{KKSL}), it immediately
follows that this weak symplectic form on $\mathcal{O}
\subset \mathfrak{g}_\ast$ has the expression
\begin{equation}\label{KKS1}
\omega_\mathcal{O}(\operatorname{Ad}^\ast_{g^{-1}} \rho)
\left(\operatorname{ad}^\ast_{\operatorname{Ad}_g \xi}
\operatorname{Ad}^\ast_{g^{-1}} \rho, 
\operatorname{ad}^\ast_{\operatorname{Ad}_g \eta}
\operatorname{Ad}^\ast_{g^{-1}} \rho \right)
= \langle\rho, [\xi,\eta]\rangle.
\end{equation}

Let now $f \in C^\infty(\mathfrak{g}_\ast)$ and 
$\rho \in \mathfrak{g}_\ast$. Since
\[
X_f(\operatorname{Ad}^\ast_{g^{-1}} \rho) =
- \operatorname{ad}^\ast_{Df(\operatorname{Ad}^\ast_{g^{-1}} \rho)}
\left(\operatorname{Ad}^\ast_{g^{-1}} \rho\right)  \in 
S_{\operatorname{Ad}^\ast_{g^{-1}} \rho},
\]
and $S_{\operatorname{Ad}^\ast_{g^{-1}} \rho}$ is the tangent space at
$\operatorname{Ad}^\ast_{g^{-1}} \rho$ to the
orbit $\mathcal{O}$, it follows that this orbit is a
Poisson submanifold of $\mathfrak{g}_\ast$. Since
\[
Df(\operatorname{Ad}^\ast_{g^{-1}} \rho) =
\operatorname{Ad}_g\left(D(f \circ 
\operatorname{Ad}^\ast_{g^{-1}})(\rho) \right)
\] 
for $g \in G$ and $\rho \in \mathfrak{g}_\ast$, it follows that 
\[
X_f(\operatorname{Ad}^\ast_{g^{-1}} \rho)
= - \operatorname{ad}^\ast_{\operatorname{Ad}_g\left(D(f \circ 
\operatorname{Ad}^\ast_{g^{-1}})(\rho) 
\right)}\left(\operatorname{Ad}^\ast_{g^{-1}} \rho\right)
\]
and hence for any $\eta \in \mathfrak{g}$ we have by \eqref{KKS}
\begin{multline*}
\omega_\mathcal{O} (\operatorname{Ad}^\ast_{g^{-1}} \rho) 
\left(X_f(\operatorname{Ad}^\ast_{g^{-1}} \rho), 
- \operatorname{ad}^\ast_{\operatorname{Ad}_g \eta}
 \left(\operatorname{Ad}^\ast_{g^{-1}} \rho\right) \right) \\ 
= \omega_\mathcal{O} (\operatorname{Ad}^\ast_{g^{-1}} \rho)
\left(- \operatorname{ad}^\ast_{\operatorname{Ad}_g
\left(D(f \circ \operatorname{Ad}^\ast_{g^{-1}})(\rho) 
\right)}\left(\operatorname{Ad}^\ast_{g^{-1}} \rho\right), 
- \operatorname{ad}^\ast_{\operatorname{Ad}_g \eta}
\left(\operatorname{Ad}^\ast_{g^{-1}} \rho\right) \right)\\
= \langle\rho, [D(f \circ 
\operatorname{Ad}^\ast_{g^{-1}})(\rho), \eta]\rangle
= -\langle \operatorname{ad}^\ast_\eta \rho, D(f \circ 
\operatorname{Ad}^\ast_{g^{-1}})(\rho)\rangle \\
= - Df(\operatorname{Ad}^\ast_{g^{-1}} \rho) \left(
\operatorname{Ad}^\ast_{g^{-1}}\left( 
\operatorname{ad}^\ast_\eta \rho \right)\right)
= - Df(\operatorname{Ad}^\ast_{g^{-1}}
\rho) \left(\operatorname{ad}^\ast_{\operatorname{Ad}_g \eta} \left(
\operatorname{Ad}^\ast_{g^{-1}} \rho\right)\right),
\end{multline*}
which shows that the Hamiltonian vector field $X_f$ relative to the Lie-Poisson
structure (\ref{general LP}) computed at a point
of the orbit $\mathcal{O}$ is also Hamiltonian relative
to the weak symplectic form (\ref{KKS}). Thus the Lie-Poisson
structure on $\mathfrak{g}_\ast$ and the weak symplectic form 
on the orbit are compatible, i.e., the Lie-Poisson structure (\ref{general
LP}) induces the weak symplectic form (\ref{KKS}) on the orbit.

Summarizing, we have shown
that each connected component of a coadjoint orbit is a connected smooth Banach
manifold, that the inclusion of the orbit in $\mathfrak{g}_\ast$ is a weak injective
immersion such that its tangent map has at each point as range the
characteristic subspace, and that the Lie-Poisson structure
induces the weak symplectic form on the orbit given by the
canonical diffeomorphism of the orbit with the quotient
$G/G_\rho$. In addition, since the orbits are
a partition of $\mathfrak{g}_\ast$, the maximality condition 
holds automatically. $\quad \blacksquare$
\medskip 

Next we analyze a remarkable particular situation that will give 
much stronger conclusions.

\begin{theorem}
\label{strong stratification theorem}
Let the Banach Lie group $G$ and the element $\rho \in \mathfrak{g}_\ast$ 
satisfy the hypotheses of Theorem \ref{quotient symplectic structure}. 
The following conditions are equivalent:
\begin{enumerate}
\item[{\rm (i)}] $\iota:
G/G_\rho \rightarrow \mathfrak{g}_\ast$ is an injective immersion;
\item[{\rm (ii)}] the characteristic subspace 
$S_\rho: = \{\operatorname{ad}_\xi \rho \mid 
\xi \in \mathfrak{g} \}$ is closed in $\mathfrak{g}_\ast$;
\item[{\rm (iii)}] $S_\rho  = \mathfrak{g}_\rho ^\circ$, where 
$\mathfrak{g}_\rho ^\circ $ is the annihilator of $\mathfrak{g}_\rho$ 
in $\mathfrak{g}_\ast$.
\end{enumerate}
Endow the coadjoint orbit $\mathcal{O}: = \{ \operatorname{Ad}^\ast _g \rho 
\mid g \in G\}$ with the manifold structure making $\iota$ a 
diffeomorphism. Then, under any of the hypotheses {\rm (i)--(iii)}, the two-form 
given by \eqref{KKS} is a strong symplectic form.
\end{theorem}

\noindent\textbf{Proof}. It is a general fact that any set is included in its
double annihilator, so $S_\rho \subset S_\rho ^{\circ \circ}$. We shall prove now
that if $S_\rho$ is closed, then this inclusion is actually an equality (this is
also a general fact for closed subspaces in any Banach space). Indeed, if $S_\rho \neq 
S_\rho ^{\circ \circ} \subset \mathfrak{g}^\ast$, then  closedness of 
$S_\rho$ guarantees by the Hahn-Banach theorem that there exists 
$0 \neq \varphi \in \mathfrak{g}^{\ast \ast}$ such that 
$\varphi \in S_\rho ^\circ$ and $\varphi \notin S_\rho^{\circ \circ\circ}$. The 
inclusion $S_\rho \subset S_\rho ^{\circ \circ}$ implies $S_\rho^{\circ \circ \circ}
\subset S_\rho ^\circ$. Since it is in general true that $S_\rho ^\circ \subset
S_\rho ^{\circ \circ \circ}$ we get $S_\rho ^\circ =
S_\rho ^{\circ \circ \circ}$ , which contradicts the existence of $\varphi$. Thus if
$S_\rho$ is closed, then $S_\rho  = S_\rho ^{\circ \circ}$.

Using the identity $\langle \eta, \operatorname{ad}^\ast_\xi \rho \rangle  
= - \langle \xi, \operatorname{ad}^\ast _\eta \rho \rangle$ for any
$\xi, \eta
\in \mathfrak{g}$, it follows that $S_\rho^\circ = \mathfrak{g}_\rho$. Taking
the annihilator in this relation and using the equality $S_\rho ^\circ = 
S_\rho ^{\circ \circ}$ yields $S_\rho = \mathfrak{g}_\rho^\circ$. Thus 
$S_\rho$ is closed if and only if $S_\rho = \mathfrak{g}_\rho^\circ$. This
proves the equivalence of (ii) and (iii). 

Assume now that (iii) holds. Since $G_\rho$ is a Banach Lie subgroup one has
the splitting $\mathfrak{g} = \mathfrak{g}_\rho \oplus \mathfrak{g}_\rho ^c$, where
$\mathfrak{g}_\rho ^c$ is a closed subspace.  This induces the splitting of the
dual space $\mathfrak{g}^\ast = \mathfrak{g}_\rho ^\circ \oplus
\left(\mathfrak{g}_\rho ^c \right)^\circ = 
S_\rho \oplus \left(\mathfrak{g}_\rho ^c \right)^\circ$. Thus, using the inclusion
$S_\rho \subset \mathfrak{g}_\ast$ we obtain the splitting 
\[
\mathfrak{g}_\ast = 
S_\rho \oplus \left[ \left(\mathfrak{g}_\rho^c\right) ^\circ \cap \mathfrak{g}_\ast
\right].
\] 
The identity
\begin{equation}
\label{moving the characteristic distribution} 
T_{[g]}\iota \left(T_{[g]}(G/G_\rho) \right) =
S_{\operatorname{Ad}^\ast _{g^{-1}}\rho} = \operatorname{Ad}^\ast _{g^{-1}} S_\rho,
\end{equation}
and the fact that $\operatorname{Ad}^\ast _{g^{-1}}$ is a Banach space isomorphism 
show that $\iota$ is a immersion. So (i) holds. Conversely, if (i) is valid then
$S_\rho$ is closed by definition, so (iii) is satisfied.

In order to prove that \eqref{KKS} is a strong symplectic form, 
let us notice that since $S_\rho$ is a closed subspace of $\mathfrak{g}_\ast$, 
by the Hahn-Banach Theorem, for any $f \in S_\rho ^\ast$ there exists an element
$\eta \in (\mathfrak{g}_\ast)^\ast  = \mathfrak{g}$ such that 
\[
f(\operatorname{ad}^\ast_\xi \rho ) =
\langle\operatorname{ad}^\ast_\xi \rho, \eta \rangle = 
\langle\rho, [\xi, \eta] \rangle 
\]
for any $\xi \in \mathfrak{g}$. So the linear map $\operatorname{ad}^\ast_\eta \rho
\in S_\rho \mapsto \omega_{\mathcal{O}}(\rho) (\cdot, 
\operatorname{ad}^\ast_\eta \rho) \in S_\rho ^\ast$ is surjective. Due
to left invariance of $\omega_\mathcal{O}$ this surjectivity will hold at all
points of the orbit $\mathcal{O}$. Since
$\omega_\mathcal{O}$ is in general a weak symplectic form according to Theorem
\ref{stratification theorem}, it follows that $\omega_\mathcal{O}$ is a strong 
symplectic form.
\quad $\blacksquare$

\begin{corollary}
\label{quasi immersion corollary}
Let the Banach Lie group $G$ and the element $\rho \in \mathfrak{g}_\ast$ 
satisfy the hypotheses of Theorem \ref{quotient symplectic structure}. Then
$\iota :G/G_\rho \rightarrow \mathfrak{g}_\ast$ is a quasi immersion 
if and only if it is an immersion.
\end{corollary}

We now apply the above theorems to the important class of $W^\ast$-algebras.
From Theorem \ref{w star algebra theorem} one knows that the
predual $\mathfrak{m}_\ast$ of the (complex) $W^\ast$-algebra
$\mathfrak{m}$ is a holomorphic Banach Lie-Poisson space. Recall
that the set $G(\mathfrak{m})$ of invertible elements of $\mathfrak{m}$ is a
Banach Lie group who acts on  $\mathfrak{m}_\ast$ by the coadjoint action.

\begin{corollary}
\label{quantum reduction corollary}
Let $\rho \in \mathfrak{m}_\ast$ be such that $G(\mathfrak{m})_\rho 
= \im R^\ast \cap G(\mathfrak{m})$, where $R^\ast$ is given by 
\eqref{r star in example two} (so $R:= (R^\ast)^\ast |_{\mathfrak{m}_\ast}$ is
a quantum reduction). Then the connected components of the coadjoint orbit  
 through $\rho$ are weakly immersed weak symplectic manifolds that are
symplectic leaves of the Banach Lie-Poisson space 
$\mathfrak{m}_\ast$.
\end{corollary}

\noindent\textbf{Proof}. By Proposition \ref{example two  proposition two},
$G(\mathfrak{m})_\rho$ is a Lie subgroup of $G(\mathfrak{m})$, so the hypotheses
of Theorems \ref{quotient symplectic structure} and \ref{stratification theorem} 
hold and the conclusion follows. \quad $\blacksquare$
\medskip

\begin{corollary}
\label{strong quantum reduction corollary}
Let $\rho \in \mathfrak{m}_\ast$ be such that $G(\mathfrak{m})_\rho 
= \im R^\ast \cap G(\mathfrak{m})$, where $R^\ast$ is given by 
\eqref{r star in example two}. Then the following conditions are equivalent:
\begin{enumerate}
\item[{\rm (i)}] $S_\rho = \ker R$;
\item[{\rm (ii)}] the map $\iota:
G/G_\rho \rightarrow \mathfrak{g}_\ast$ is an injective immersion.
\end{enumerate}
Under any of these conditions, the coadjoint orbit $\mathcal{O}$ 
endowed with the smooth manifold structure making $\iota$ a diffeomorphism 
onto its image is strong symplectic. Additionally, the Poisson tensor 
on $S_\rho$ defined by the strong symplectic form 
$\omega_{\mathcal{O}}(\rho)$ coincides with the one coinduced 
by the quantum reduction $R : \mathfrak{m}_\ast \rightarrow \im R$.
\end{corollary}

\noindent\textbf{Proof}. If $S_\rho = \ker R$, then $S_\rho$ is closed and 
Theorem \ref{strong stratification theorem} applies thus guaranteeing that 
$\iota$ is an injective immersion. 

Conversely, if $\iota$ is an immersion,
Theorem \ref{strong stratification theorem} states that $S_\rho = 
\mathfrak{m}_\rho ^\circ$.  However, the hypothesis and Proposition 
\ref{example two proposition two} guarantee that $G_\rho$ is a Lie subgroup
of $G(\mathfrak{m})$ whose Lie algebra is $\im R^\ast$. On the other hand
it is clear that the Lie algebra of $G(\mathfrak{m})_\rho$ is $\mathfrak{m}_\rho$
since the $\exp (\lambda x) \in G(\mathfrak{m})_\rho$ for all 
$x \in \mathfrak{m}_\rho$ and all $\lambda \in \mathbb{C}$ (see Bourbaki [1972], 
Chapter III, \S 6.4, Corollary 1). Therefore $\mathfrak{m}_\rho = \im R^\ast$.

The decomposition $\mathfrak{m} = 
\im R^\ast \oplus \ker R^\ast$ and the one induced on the dual imply the 
general identities
\[
\left( \im R^\ast \right)^\circ = \ker R^{\ast \ast} \quad \text{~and~}
\quad \left( \ker R^\ast \right) ^ \circ= \im R^{\ast \ast}.
\]
Therefore, using  $\mathfrak{m}_\rho = \im R^\ast$, we get
$\ker R = \ker R^{\ast \ast} \cap \mathfrak{m}_\ast = 
\left(\im R^\ast \right)^\circ \cap \mathfrak{m}_\ast
= \mathfrak{m}_\rho^\circ = S_\rho$ and the equivalence of (i) and (ii) is 
proved.

The other statements follow by directly applying Theorem \ref{strong 
stratification theorem} and the definition of quantum reduction. 
 \quad $\blacksquare$
\medskip

\begin{example}
\label{trace class operators example}
\normalfont
Take in the previous considerations $\mathfrak{m} = \bounded$, 
$\mathfrak{m}_\ast = \one$, and the quantum reduction map $R: \one 
\rightarrow \one$ defined by \eqref{measurement}, where $\sum_{k=1}^N
P_k = 1$, with $N \in \mathbb{N}$ or $N = \infty$. If $N = \infty$, let
\begin{equation}
\label{spectral decomposition of rho infinite}
\rho = \sum _{k=1} ^\infty \lambda_k P_k, \quad \lambda_k \in \mathbb{C}, 
\quad \lambda_k \neq \lambda_\ell \neq 0 \;\text{~for~} 
\; k \neq \ell, \quad \operatorname{rank}P_k < \infty \; \text{~for~} k \geq 1
\end{equation}
and if $N < \infty$ let
\begin{equation}
\label{spectral decomposition of rho finite}
\rho = \sum _{k=1} ^N \lambda_k P_k, \; \lambda_k \in \mathbb{C}, 
\; \lambda_k \neq \lambda_\ell \; \text{~for~} 
\; k \neq \ell, \quad \lambda_1 = 0, \quad \operatorname{rank}P_k < \infty
\;\text{~for~} k \geq 2. 
\end{equation}
Thus $\rho \in \one$. It is easy to check that 
\[
\bounded_\rho = \im R^\ast = \left\{\sum_{k=1}^NP_k X P_k \mid X \in 
\bounded  \right\}.
\]
So all conclusions of Theorem \ref{stratification theorem} hold, that is,
the coadjoint orbit  $\mathcal{O}:= \{ g\rho g^{-1} \mid g \in 
\invertible \} \subset \one$ through $\rho$ is weakly immersed in $\one$ 
and is weakly symplectic relative to the two-form
\begin{equation}
\label{KKS for one}
\omega_\mathcal{O}(g\rho g^{-1})([gXg^{-1}, g\rho g^{-1}], [gYg^{-1}, g\rho
g^{-1}]) = \tr(\rho[X,Y])
\end{equation}
for $\rho \in \one$ given by \eqref{spectral decomposition of rho
infinite} or \eqref{spectral decomposition of rho finite}, $g \in
\invertible$, $X,Y \in
\bounded$. In addition,  for $\mathcal{M}$ a
complex separable Hilbert space, recall that the group
$\invertible$ is path connected (see e.g. Boos and Bleecker [1983]
\S I.6) and hence the coadjoint orbit is also connected; thus it is a 
symplectic leaf of the Banach Lie-Poisson space $\one$.

The characteristic subspace $S_\rho=\{[X, \rho] \mid X \in \bounded\}$ 
is contained in 
\[
\ker R = \left\{\sum_{k \neq \ell} P_kXP_\ell \mid X \in \bounded \right\}
\]
and if $N \in \mathbb{N}$ one has 
\[
\sum_{k \neq \ell} P_k X P_\ell = [\rho, Y]
\]
for some $Y \in \bounded$ which is related to $X$ through the system 
of equations
\[
P_k X P_\ell = (\lambda_k - \lambda _\ell) P_k Y P_\ell
\]
for all $k \neq \ell$. Note that if $N = \infty$, the above system is not
solvable for some $Y \in \bounded$. Therefore $S_\rho = \ker R$ and by
Theorem \ref{strong stratification theorem} one concludes that the connected
coadjoint orbit is immersed in $\one$ and that it is strongly symplectic.
\quad $\blacklozenge$
\end{example} 

\noindent \textbf{Remark}. The weak symplectic form 
$\omega_\rho$ given by \eqref{KKSL} or \eqref{KKSR} is, in general, not strong since  
$\omega_\rho([e])$ is, in general, not a strong bilinear form on 
$T_{[e]}(G/G_\rho)$. To begin with, one notices that $T_e \pi(\xi)
\in T_{[e]}(G/G_\rho) \mapsto [\xi] \in 
\mathfrak{g}/\mathfrak{g}_\rho$ is a linear continuous bijective map and hence 
a Banach space isomorphism. Thus 
$\omega_\rho([e])$ can be viewed as a bilinear continuous map
$\omega_\rho([e]) : \mathfrak{g}/\mathfrak{g}_\rho \times
\mathfrak{g}/\mathfrak{g}_\rho \rightarrow \mathbb{R}$ given by
$\omega_\rho([e])([\xi], [\eta]) = \langle \rho, [\xi, \eta]\rangle$. The map 
$[\xi] \in \mathfrak{g}/\mathfrak{g}_\rho \mapsto \omega_\rho([e])([\xi], 
[\cdot]) = \langle \operatorname{ad}_\xi ^\ast \rho, \cdot \rangle
\in (\mathfrak{g}/\mathfrak{g}_\rho)^\ast$ is clearly linear continuous and 
injective. Thus, if the symplectic form $\omega_\rho$ were strong, then the
Banach spaces $\mathfrak{g}/\mathfrak{g}_\rho$ and 
$(\mathfrak{g}/\mathfrak{g}_\rho)^\ast$ would necessarily be isomorphic. Here is
a concrete situation in which this cannot occur.

Consider the case described in Example \ref{trace class operators example} 
for the trace class operator $\rho$ given by \eqref{spectral decomposition of rho infinite}. Then 
$\bounded / \bounded_\rho \cong \ker R^\ast$ since $\bounded_\rho = \im R^\ast$.
In the proof of Corollary \ref{strong quantum reduction corollary} we have seen
that $(\bounded_\rho)^\circ = \ker R$ and hence $(\bounded/\bounded_\rho)^\ast
\cong \ker R$. Thus, strongness of the symplectic form would imply the 
isomorphism $\ker R \cong \ker R^\ast$. We shall show now that if $\rho$ is
an infinte rank operator in $\one$, then this is impossible. Indeed, there exists 
$Y \in \bounded$, $Y \notin \one$ such that the expression
$\tr(YX) = \tr \sum_{k \neq \ell} P_k Y P_\ell X P_k$ does not make sense for all $X
\in \ker R^\ast \subset \bounded$. The functional $\sigma \in \ker R \subset
\one \mapsto \tr(Y\sigma) \in \mathbb{C}$ is, however, an element of $(\ker R)^\ast$.
Thus, for such $\rho$ the 
orbit symplectic form \eqref{KKS for one} is not strong.

If $\rho$ is given by \eqref{spectral decomposition of rho finite}, the above
argument breaks down. Moreover, as was shown in Corollary 
\ref{strong quantum reduction corollary} the orbit symplectic form 
\eqref{KKS for one} is in this case strong.
\quad $\blacklozenge$
\medskip

It was shown in Bona [2000] that unitary group coadjoint
orbits through Hermitian finite rank operators are always strong symplectic 
manifolds. We present  this case below.

\begin{example}
\label{unitary operators example}
\normalfont
We apply the considerations of this section to the real closed
Banach Lie subgroup  $\unitary : = \{U \in \bounded \mid
UU^\ast = U^\ast U = I \}$ of unitary elements of $\invertible$ (Bourbaki [1972],
Chapter III, \S 3.10, Corollary 2). Its Lie
algebra consists of the skew Hermitian bounded operators $\sherm := \{ X
\in \bounded \mid X + X^\ast = 0 \}$; this is a closed split real Banach Lie
subalgebra of $\bounded$. To study this case, we place ourselves in the 
context of Example \ref{example following Proposition}, 
that is, we take $\mathfrak{g} = \bounded$, $\mathfrak{g}_\ast = \one$,
${\mathfrak{g}_\ast}_\mathbb{R} = \one_\mathbb{R}$,
$\mathfrak{g}_\mathbb{R} = \bounded_\mathbb{R}$ (in other words, 
the Banach spaces $\one$ and $\bounded$ thought of as a real Banach spaces), 
the continuous $\mathbb{R}$-linear involution $\sigma: \one_\mathbb{R} 
\rightarrow \one_\mathbb{R}$ given by $\sigma \rho = -\rho^\ast$, for 
$\rho \in \one$,
and the complex structure $I$ given by $I\rho = i\rho$. It is easily
verified that the involution $\sigma$ satisfies the conditions (i), 
(ii), and (iii) of Example \ref{example following Proposition}. Then, by
construction, $\mathfrak{g}_\ast^\sigma = \{ \rho \in \one \mid
\rho + \rho^\ast =0 \} = :\shermone$ and, as was shown in
Example \ref{example following Proposition}, $\shermone$ is a real 
Banach Lie-Poisson space and the map $R: \one_\mathbb{R} \rightarrow
\shermone$ given by $R = (id + \sigma)/2$ is a linear Poisson map. 

The same type of argument as in Example \ref{trace class operators example}
shows that one can directly apply Theorems \ref{stratification theorem} and 
\ref{strong stratification theorem} to $G =\unitary$ and $\mathfrak{g}_\ast =
\shermone$ . The symplectic leaves in this 
case correspond to the  infinite dimensional flag manifolds and the strong
symplectic form given by  \eqref{KKS for one} ($\rho$ of finite rank and the 
arguments in the correct spaces) 
coincides with the imaginary part of the natural K\"ahler metric on these
manifolds. A particular example of such an infinite dimensional flag manifold
is the projectivized Hilbert space $\mathbb{CP}(\mathcal{M})$ thought of as
immersed in $\one$ as the coadjoint orbit through $\rho :=
|\psi \rangle \langle \psi |/ \langle \psi | \psi \rangle $ for any 
$|\psi \rangle \in \mathcal{M}$. 
\quad $\blacklozenge$
\end{example}

We next discuss the cotangent bundle of a
Banach Lie group and introduce a remarkable submanifold, called 
in the sequel the precotangent bundle. Consider a Banach Lie group 
$G$ with Banach Lie algebra $\mathfrak{g}$ admitting a predual 
$\mathfrak{g}_\ast$ and assume that 
$\operatorname{Ad}^\ast_G(\mathfrak{g}_\ast) \subset \mathfrak{g}_\ast$. If
$L_g$ and $R_g$ denote the left and right translations by $g \in G$
respectively, it follows that
$T_g L_{g^{-1}}:T_gG \rightarrow T_e G = \mathfrak{g}$ and $T_g
R_{g^{-1}}:T_gG \rightarrow T_e G = \mathfrak{g}$ are a Banach space
isomorphisms.  Their duals $T_g^\ast L_{g^{-1}} : \mathfrak{g}^\ast
\rightarrow T_g^\ast G$ and $T_g^\ast R_{g^{-1}} : \mathfrak{g}^\ast
\rightarrow T_g^\ast G$ are therefore also Banach space isomorphisms.
Define ${T_g}_\ast G : = T_g^\ast L_{g^{-1}} \mathfrak{g}_\ast$, $T_\ast G
:=  \cup_{g \in G} {T_g}_\ast G$, and conclude, as usual, that $T_\ast G$ is
a vector bundle over $G$ which is also a subbundle of $T^\ast G$ (see, e.g.
Abraham, Marsden, Ratiu [1988] for such an argument); it will be called
the {\bfi precotangent bundle\/} of $G$. This construction could have been equally
well done using right translations since $T_g^\ast R_{g^{-1}}
\operatorname{Ad}_{g^{-1}}^\ast \mathfrak{g}_\ast =
T_g^\ast L_{g^{-1}} \mathfrak{g}_\ast$ and, by hypothesis,
$\operatorname{Ad}_g^\ast \mathfrak{g}_\ast = \mathfrak{g}_\ast$ for any $g \in G$.
The precotangent bundle $T_\ast G$ has been constructed
using the left trivialization $L:  T_\ast G \rightarrow G \times
\mathfrak{g}_\ast$, $L(\rho_g) := (g, T_e^\ast L_{g} \rho_g)$
with inverse $L^{-1} (g, \rho) = T_g^\ast L_{g^{-1}} \rho$, for $\rho_g
\in {T_g}_\ast G$ and $\rho \in \mathfrak{g}_\ast$. Completely analogous
formulas hold for  the right trivialization $R:T_\ast G \rightarrow G \times
\mathfrak{g}_\ast$, $R(\rho_g) := (g, T_g^\ast R_{g^{-1}} \rho_g)$,
$R^{-1} (g, \rho) = T_e^\ast R_{g^{-1}}\rho$; $L$ and $R$ are vector bundle
isomorphisms covering the identity of $G$. Denote by $\pi:
T^\ast G \rightarrow G$ the cotangent bundle projection and use the
same letter to denote its restriction to $T_\ast G$.

The usual construction of the canonical one-form on the cotangent bundle
$T^\ast G$ works also in the case of the precotangent bundle $T_\ast G$.
Indeed, define the one form $\Theta$  on $T^\ast G$ or on $T_\ast G$ by
\begin{equation}
\label{canonical one form}
\Theta(\rho_g)(v) : = \langle \rho_g, T_{\rho_g} \pi (v) \rangle
\end{equation}
for any $\rho_g \in T_g^\ast G$ (respectively ${T_g}_\ast G$), $v \in
T_{\rho_g}(T^\ast G)$ (respectively $ T_{\rho_g}(T_\ast G)$) and where the
pairing is between ${T^\ast_g} G$ and $T_g G$ (respectively
${T_g}_\ast G$ and $T_g G$). Left trivialized, this formula reads
\begin{equation}
\label{left trivialized canonical one form}
\Theta_L(g, \rho)(u_g, \mu):=(L_\ast\Theta)(g, \rho)(u_g, \upsilon) 
= \langle \rho, T_gL_{g^{-1}}u_g \rangle
\end{equation}
for $g \in G$, $u_g \in T_gG$, and $\rho, \mu \in
\mathfrak{g}^\ast$ (respectively $\mathfrak{g}_\ast$), where the pairing is
now between $\mathfrak{g}^\ast$ (respectively
$\mathfrak{g}_\ast$) and $\mathfrak{g}$. Define the canonical symplectic
form on $T^\ast G$ or $T_\ast G$ by $\Omega: = - d
\Theta$ and let $\Omega_L: = L_\ast \Omega$. A computation identical to the one in the
finite dimensional case using the identity $\Omega_L = - \mathbf{d} \Theta_L$ (see
Abraham and Marsden [1978],
\S 4.4), leads to the expression of the canonical two-form in the left trivialization 
\begin{align}
\label{left trivialized canonical two form}
&\Omega_L(g, \rho)\left((u_g, \mu), (v_g, \nu) \right) 
:= (L_\ast \Omega)(g, \rho)\left((u_g, \mu), (v_g, \nu) \right) 
\nonumber \\
& \qquad = \langle \nu, T_gL_{g^{-1}}u_g \rangle - \langle \mu,
T_gL_{g^{-1}}v_g \rangle - \langle \rho, [T_gL_{g^{-1}}u_g ,
T_gL_{g^{-1}}v_g ] \rangle,
\end{align}
where $g \in G$, $u_g, v_g \in T_gG$, and $ \rho, \mu,\nu \in
\mathfrak{g}^\ast$ (respectively $\mathfrak{g}_\ast$). This formula
immediately shows that $\Omega_L$ and hence
$\Omega$ is a weak symplectic form on both $T^\ast G$ and $T_\ast G$. We
shall see below that it is not strong in general, for different reasons,
on both $T^\ast G$ and  $T_\ast G$. 

To show that $\Omega_L$ is strong on $G \times \mathfrak{g}^\ast$, one
needs to prove that for fixed $(g, \rho) \in G \times \mathfrak{g}^\ast$
the linear continuous map 
$(u_g, \mu) \in T_g G \times \mathfrak{g}^\ast \mapsto \Omega_L(g,
\rho)\left((u_g, \mu), (\cdot, \cdot) \right) \in (T_g G \times
\mathfrak{g}^\ast)^\ast \cong T^\ast_g G \times
\mathfrak{g}^{\ast \ast}$ is surjective, that is, given 
$\alpha_g \in T_g^\ast G$ and $\Gamma \in \mathfrak{g}^{\ast \ast}$ one can find $u_g
\in T_g G$ and $\mu \in \mathfrak{g}^\ast$ such that
\[
\langle \nu, T_gL_{g^{-1}}u_g \rangle - \langle \mu +
\operatorname{ad}^\ast_{T_gL_{g^{-1}}u_g} \rho, T_gL_{g^{-1}}v_g \rangle  =
\langle \alpha_g, v_g \rangle + \langle \Gamma, \nu \rangle
\]
for all $v_g \in T_g G$, $\nu \in \mathfrak{g}^\ast$.
If this were possible, a necessary condition is that $\Gamma =
T_gL_{g^{-1}}u_g$ which is not the case if $\mathfrak{g}$ is strictly
included in $\mathfrak{g}^{\ast \ast}$. If the Banach space
$\mathfrak{g}$ is reflexive then $\Gamma \in \mathfrak{g}$ and one can choose $\mu =
-\operatorname{ad}^\ast _\Gamma \rho - T_e^\ast L_g \alpha_g$. Thus, \textit{if
$\mathfrak{g}$ is reflexive, the canonical weak symplectic form on $T^\ast G$ is
strong.}

Next we analyze $\Omega_L$ on $G \times \mathfrak{g}_\ast$. As before, we
fix $(g, \rho) \in G \times \mathfrak{g}_\ast$ and study the linear
continuous map $(u_g, \mu) \in T_g G \times \mathfrak{g}_\ast \mapsto
\Omega_L(g, \rho)\left((u_g, \mu), (\cdot, \cdot) \right) \in (T_g G \times
\mathfrak{g}_\ast)^\ast \cong T^\ast_g G \times
\mathfrak{g}$. To prove its surjectivity one needs to find for given 
$\alpha_g \in T_g^\ast G$ and $\xi \in \mathfrak{g}$ 
a vector $u_g \in T_g G$ and a form $\mu \in \mathfrak{g}^\ast$ such that 
\[
\langle \nu, T_gL_{g^{-1}}u_g \rangle - \langle \mu +
\operatorname{ad}^\ast_{T_gL_{g^{-1}}u_g} \rho, T_gL_{g^{-1}}v_g \rangle  =
\langle \alpha_g, v_g \rangle + \langle \nu,
\xi \rangle
\]
for all $v_g \in T_g G$, $\nu \in \mathfrak{g}_\ast$. This identity
implies that $u_g = T_eL_g \xi$, which, unlike the previous case, is
possible. However, this identity also requires that 
$\mu = -\operatorname{ad}^\ast _\xi \rho -  T_e^\ast L_g \alpha_g$ which
is, in general, impossible to achieve since $T_e^\ast L_g \alpha_g \in
\mathfrak{g}^\ast$ but is not necessarily an element of
$\mathfrak{g}_\ast$.

In spite of this obstruction, there is a Poisson bracket on $G \times
\mathfrak{g}_\ast$. Given $f, h \in C^\infty(G \times \mathfrak{g}_\ast)$
their Poisson bracket is given by
\begin{multline}
\label{Poisson bracket in left trivialization}
\{f, h\}(g, \rho) = \left\langle d_1f(g, \rho), T_eL_g d_2h(g, \rho)
\right\rangle - \left\langle d_1h(g, \rho), T_eL_g d_2f(g, \rho)
\right\rangle \\ 
- \left\langle \rho, \left[d_2 f(g, \rho), d_2 h(g,
\rho)\right] \right\rangle,
\end{multline}
where $d_1f(g, \rho) \in T^\ast_g G$ and $d_2f(g, \rho) \in
(\mathfrak{g}_\ast)^\ast = \mathfrak{g}$ are the first and second
partial derivatives of $f$, the pairing in the first two terms is between $T^\ast_g G$
and $T_gG$, whereas in the third term it is between $\mathfrak{g}_\ast$ and
$\mathfrak{g}$. Conditions (i) and (ii) in Definition \ref{Poisson manifold
definition} are satisfied. The proof of the Jacobi identity is a tedious direct
verification. However, condition (iii) does not hold. Indeed formula \eqref{Poisson
bracket in left trivialization} shows that the Hamiltonian vector field of $h$, if
well defined, must have the expression
\begin{equation}
\label{Hamiltonian vector field left trivialized}
X_h(g, \rho) = \left(T_eL_g d_2h(g, \rho), 
\operatorname{ad}^\ast_{d_2h(g,\rho)} \rho - T_e^\ast L_g d_1h(g, \rho)\right)
\end{equation}
The same obstruction encountered in the attempted proof of the strongness
of $\Omega_L$ on $G\times \mathfrak{g}_\ast$ appears here in the second
summand of the second component: the term $T_e^\ast L_g d_1h(g, \rho)$ is,
in general, not an element of $\mathfrak{g}_\ast$. 

Thus, \textit{unlike $T^\ast G$, the precotangent bundle $T_\ast G$ is
naturally endowed with a Poisson bracket, but is not a Poisson manifold in
the sense of Definition} \ref{Poisson manifold definition}. However, before
Lie-Poisson reduction, the unreduced space $G \times \mathfrak{g}_\ast$ is only
weakly symplectic, admits a Poisson bracket, but has no Hamiltonian vector
fields. For functions admitting Hamiltonian vector fields, the Poisson
bracket is naturally induced by the weak symplectic form which is the pull
back of the canonical symplectic form on the cotangent bundle of the
group. Finally, the projection $G\times \mathfrak{g}_\ast \rightarrow
\mathfrak{g}_\ast$ preserves the Poisson brackets, if one changes the sign
of the Lie-Poisson bracket on $\mathfrak{g}_\ast$. Similar considerations
can be carried out for right translations and one obtains, as in finite
dimensions, a dual pair (Weinstein [1983], Vaisman [1994])
\[
{\mathfrak{g}_\ast}_+ \longleftarrow T_\ast G \longrightarrow
{\mathfrak{g}_\ast}_-
\]
where the signs refer to the Lie-Poisson bracket on $\mathfrak{g}_\ast$;
the two arrows are the momentum maps for left and right translations (see
\S \ref{section: momentum maps and reduction} for a presentation of
momentum maps in our setting and Marsden and Ratiu [1994] for more information in the
finite dimensional  case).
\medskip

\section{Momentum maps and reduction}
\label{section: momentum maps and reduction}

In this section we shall explore the relationship between the
 classical theory of reduction for Poisson manifolds discussed in
 \S \ref{section: classical reduction} and that of quantum
reduction presented in \S \ref{section: quantum reduction}. We
shall show that this link will be crucial for the integration and
quantization of Hamiltonian systems.

We shall introduce a definition of the momentum map which is a
direct generalization of this concept from finite dimensional
Poisson geometry.

\begin{definition}
\label{momentum map definition}
A {\bfi momentum map\/} is a Poisson map $J: P \rightarrow
\mathfrak{b}$ from a Banach Poisson manifold $P$ to a Banach
Lie-Poisson space $\mathfrak{b}$.
\end{definition}

Recall that $\mathfrak{b}$ is a Banach space, that $C^\infty
(\mathfrak{b})$ is endowed with a Poisson bracket $\{\cdot,\cdot \}$, that
$\mathfrak{b}^\ast$ is closed under this bracket, and that
$\operatorname{ad}^\ast_{\mathfrak{b}^\ast} \mathfrak{b} \subset
\mathfrak{b}$. Thus $\mathfrak{b}^\ast$ is a Lie algebra and the
prescription $\xi \in \mathfrak{b}^\ast \mapsto \xi_P : = X_{\xi
\circ J}$ defines a (left) Lie algebra action on $P$, that is,
$[\xi_P, \eta_P] = -[\xi, \eta]_P$ for any $\xi, \eta \in
\mathfrak{b}^\ast$.  Indeed, recalling that the Hamiltonian vector
field defined by $h \in C^\infty(P)$ is defined by $df (X_h) =
\{f, h\}$, the Jacobi identity for the Poisson bracket is
equivalent to $[X_f, X_g ] = -X_{\{f,g\}}$. Using this relation
and the fact that $J$ is a Poisson map, we conclude
\[
[\xi_P, \eta_P] = [X_{\xi \circ J}, X_{\eta \circ J}] = -X_{\{\xi
\circ J, \eta \circ J\}} = -X_{\{\xi, \eta\} \circ J} = -[\xi,
\eta]_P
\]
which proves that $\xi \mapsto \xi_P$ is indeed a (left) Lie
algebra action.

\begin{theorem}
\label{noether}
(Noether's Theorem) If $h \in C^\infty(P)$ is
$\mathfrak{b}^\ast$-invariant, that is, $\langle dh, \xi_P \rangle
= 0$ for all $\xi \in \mathfrak{b}^\ast$, then $J$ is conserved
along the flow of the Hamiltonian vector field $X_h$.
\end{theorem}

\noindent\textbf{Proof}. The condition of invariance states that
\[
0 = \langle dh, X_{\xi \circ J} \rangle = \{h, \xi \circ J\} =
-\langle d(\xi \circ J), X_h \rangle
\]
for every $\xi \in \mathfrak{b}^\ast$, which is equivalent to
\[
\frac{d}{dt} \sigma_h(t)^\ast (\xi \circ J) = 0
\]
for every $\xi \in \mathfrak{b}^\ast$, where $\sigma_h(t)$ is the
flow of $X_h$. This in turn means that $\xi \circ J \circ
\sigma_h(t) = \xi \circ J$, which says that $\langle \xi, (J \circ
\sigma_h(t) - J)(b) \rangle = 0$ for all $b \in \mathfrak{b}$.
Since $\xi$ is arbitrary, one concludes that $(J \circ \sigma_h(t)
- J)(b) = 0$ for every $b \in \mathfrak{b}$, that is, $J \circ
\sigma_h(t) = J$ for all $t$. $\quad \blacksquare$

\medskip

A Hamiltonian system $(P, \{\cdot,\cdot\}_P, h)$ is called {\bfi
collective\/} if there is a momentum map $J :P \rightarrow
\mathfrak{b}$ and a function $H \in C^\infty(\mathfrak{b})$ such
that $h = H \circ J$. Therefore, the corresponding Hamiltonian
vector fields $X^P_h$ and $X^\mathfrak{b}_H$ on $P$ and
$\mathfrak{b}$ respectively are $J$-related, that is,
\[
TJ \circ X^P_h = X^\mathfrak{b}_H \circ J
\]
which is equivalent to the commutation of the respective flows
$\sigma_h(t)$ and $\sigma_H(t)$ of $X^P_h$ and $X^\mathfrak{b}_H$
respectively, that is,
\[
\sigma_H (t)\circ J = J \circ \sigma_h(t).
\]
Thus the integration of Hamilton's equations on $\mathfrak{b}$ given by
\begin{equation}
\label{equation on b} \dot{b} = - \operatorname{ad}^\ast_{DH(b)}b,
\end{equation}
where $b \in \mathfrak{b}$, leads to the partial integration of
Hamilton's equations
\begin{equation}
\label{equation on P} \dot{f} = \{f, h\}_P
\end{equation}
on $P$, where $f \in C^\infty(P)$. If $J: P \hookrightarrow
\mathfrak{b}$ is an embedding, then solving \eqref{equation on b}
is equivalent to solving \eqref{equation on P}. In the other
extreme case, namely when the trajectory $\{\sigma_H(t)(b) \mid t
\in \mathbb{R} \}$ is a point $b \in \mathfrak{b}$, any trajectory $\sigma_H(t)(p)$
with $p \in J^{-1}(b)$ remains in the level set $J^{-1}(b)$. If there is a
distribution $E$ covering this level set satisfying the hypotheses of the classical
reduction theorem, then the trajectories above drop to the quotient and one is led to
the problem of solving a reduced system. Provided this can be
integrated, for example, if the reduced system is integrable, then the standard
reconstruction method (Abraham and Marsden [1978], \S4.3) gives the solution of the
original system on the level set of the momentum map.

In the special case when $\mathfrak{b} = \one$, equation
\eqref{equation on b} assumes the form of the nonlinear Liouville-von
Neumann equation
\begin{equation}
\label{Liouville-von Neumann eq}
\dot{\rho} = [DH(\rho), \rho]
\end{equation}
for $\rho \in \one$. The search for a collective Hamiltonian
system on $P$ is equivalent to finding a ``Lax representation" on
$\one$. The functions $T_k : = (\tr \rho^k)/k$, $k \in \mathbb{N}$
are Casimir functions on $\one$. Therefore the functions $t_k : =
T_k \circ J$, $k \in \mathbb{N}$, are, in general, integrals of
motion in involution for the system \eqref{equation on P}. If $J: P \hookrightarrow
\mathfrak{b}$ is an embedding, then they are Casimirs of the
system given by the phases space $P$. So, the problem of integration of
an equation having Lax representation in
$\one$ reduces to a large extent to the integration of the equation 
\eqref{Liouville-von Neumann eq}. 

For direct investigations of the nonlinear Liouville-von Neumann equation
in the physics literature, see for example Leble and Czachor [1998] and references therein. 

We shall illustrate the above considerations with the example of
the infinite Toda lattice system associated to the Banach
Lie-Poisson space $\one$.

\begin{example} \textbf{The infinite Toda lattice system}. 
\normalfont
The details for this example can be found in Odzijewicz and Ratiu [2003].
On the weak symplectic Banach space $\ell^\infty \times \ell^1 = (\ell^1)^\ast \times \ell^1$ define the Hamiltonian
\begin{equation}
\label{Toda lattice Hamiltonian}
H(\mathbf{q}, \mathbf{p}) := \frac{1}{2} \sum_{k=1}^\infty p_k^2 +
\sum_{k=1}^\infty \alpha_k\lambda_k \exp(q_k - q_{k+1}),
\end{equation}
where $\lambda_k \neq 0$, $\{\alpha_k\}, \{\lambda_k\} \in \ell^1$. Since the Banach space on which
this Hamiltonian is defined is only weak symplectic, the existence of the Hamiltonian
vector field associated to $H$ is not guaranteed. Formally, this Hamiltonian is that
of the Toda lattice. We will also assume that $\sum_{k=1}^\infty p_k = 0$, which means
that the velocity of the center of mass is zero. We also observe that $H$ is invariant
relative to the action of $\mathbb{R}$ on the $\mathbf{q}$-space by translation. Thus we
shall consider $H$ defined on the weak symplectic Banach space
$(\ell^\infty /\mathbb{R}\mathbf{q}_0)
\times \ell^1_0$, where $\ell^\infty /\mathbb{R}\mathbf{q}_0$ is the quotient Banach space
by the closed subspace $\mathbb{R}\mathbf{q}_0$, where $q_{0k} = 1$ for all $k \in \mathbb{N}$,
and $\ell^1_0 : = \ker \mathbf{q}_0$. Relative to the canonical coordinates 
$x_k : = q_k - q_{k+1}$ on $\ell^\infty /\mathbb{R}\mathbf{q}_0$ and $p_k$ on $\ell^1_0$
the weak symplectic form has the expression 
\begin{equation}
\label{Toda symplectic from}
\omega = - d\left(\sum_{k=1}^\infty p_k dx_k \right)= \sum_{k=1}^\infty dx_k \wedge dp_k.
\end{equation}

Let $\mathcal{M}$ be a separable Hilbert space. Let 
$P_n : =|n\rangle \langle n| : \mathcal{M} \rightarrow
\mathcal{M}$ be the rank one projection onto the span of
$|n\rangle$. If $\rho \in \one$ and $X \in \bounded$ write 
\[
\rho = \sum_{n, m = 1}^\infty \rho_{nm} |n \rangle \langle m|
\qquad \text{and} \qquad
X = \sum_{n, m = 1}^\infty X_{nm} |n \rangle \langle m|
\]
where $\rho_{nm} : =\langle n|\rho|m\rangle$ and $X_{nm} := \langle
n|X|m\rangle$. From Example \ref{section 7, third example} we know that
$\one_-$ is a Banach Lie-Poisson space whose dual is the Banach Lie
algebra $\bounded_+$. The paring between these two spaces is given by
\[
\langle \rho_-, X_+ \rangle : = \operatorname{tr}(\rho_- X_+) = \sum_{n
\geq m} \rho_{nm} X_{mn},
\]
where $\rho_- \in \one_-$ and $X_+ \in \bounded_+$. Using this pairing,
a direct verification shows that the coinduced Poisson bracket of
$\one_-$ has the expression
\begin{align}
\label{Toda bracket}
\{f, g\}_{\one_-} (\rho_-) &= \operatorname{tr}([(\pi^1_-)^\ast
(df(\rho_-)), 
(\pi^1_-)^\ast (dg(\rho_-))]\rho_-) \nonumber \\
&= \sum_{n \geq m \geq \ell}\left(\frac{\partial f}{\partial \rho_{nm}}
\frac{\partial g}{\partial \rho_{m\ell}} -
\frac{\partial g}{\partial \rho_{nm}}
\frac{\partial f}{\partial \rho_{m\ell}}
\right) \rho_{n\ell},
\end{align}
where $\pi^1_-: \one \rightarrow \one_-$ is the projector given by 

\[
\pi^1_-(\rho): = \sum_{n\geq m} \rho_{nm} |n \rangle \langle m|,
\] 
$(\pi^1_-)^\ast: (\one_-)^\ast \cong \bounded_+ \rightarrow \bounded$ is its 
dual, and the formula 
\[
(\pi^1_-)^\ast (df(\rho_-))_{mn} = \frac{\partial f}{\partial
\rho_{nm}}(\rho_-)
\]
was used in the proof of the second equality.

Define the Flaschka transformation
$J: (\ell^\infty /\mathbb{R}\mathbf{q}_0)
\times \ell^1_0 \rightarrow \one_-$ by
\begin{equation}
\label{flaschka transformation} J(\mathbf{z}) =
\sum_{k=1}^\infty\left(p_k |k\rangle \langle k| + \lambda_k e^{x_k} | k+1 \rangle
\langle k | \right).
\end{equation}
One verifies  that $J$ is a smooth injective map whose tangent map 
at every point is a continuous injection.

According to \S \ref{section: symplectic leaves}, the $\operatorname{Ad}^\ast (GL^\infty_+
(\mathcal{M}))$-orbit
\[
\mathcal{O}_\Lambda := \{\pi^1_-(g \Lambda g^{-1}) \mid g \in GL^\infty_+(\mathcal{M}) \} 
\]
where $\Lambda = J(\mathbf{0}, \mathbf{0}) \in L^1_-$, is a symplectic leaf in the Banach Lie-Poisson
space $L^1_-(\mathcal{M})$ with the symplectic form given by (\ref{KKS}). One can show that:
\begin{enumerate}
\item[(i)] $J^{-1}\left( J \left((\ell^\infty /\mathbb{R}\mathbf{q}_0)
\times \ell^1_0\right) \cap \mathcal{O}_\Lambda \right) = (\ell^\infty
/\mathbb{R}\mathbf{q}_0)
\times \ell^1_{0f} $, where  $\ell^1_{0f}$ is the linear dense subspace of $\ell^1_0$ which consists
of elements $\mathbf{p}$ with only finitely many non vanishing components;
\item[(ii)] if $(\mathbf{x}, \mathbf{p}) \in (\ell^\infty
/\mathbb{R}\mathbf{q}_0)
\times \ell^1_{0f}$ the tangent map $T_{(\mathbf{x}, \mathbf{p})} J :
(\ell^\infty /\mathbb{R}\mathbf{q}_0) \times \ell^1_{0f} \rightarrow
T_{J(\mathbf{x}, \mathbf{p})}
\mathcal{O}_\Lambda$ is a continuous linear bijection;
\item[(iii)] $J^\ast \omega_{\mathcal{O}_\Lambda} = \omega$ on $(\ell^\infty
/\mathbb{R}\mathbf{q}_0)
\times \ell^1_{0f}$.
\end{enumerate}
The map \eqref{flaschka transformation} is not a momentum map in the
sense of Definition \ref{momentum map definition} since $(\ell^\infty /\mathbb{R}\mathbf{q}_0)
\times \ell^1_0$ is a weak symplectic manifold and thus not Poisson
according to our Defintion \ref{Poisson manifold definition}. However,
the above points show that \eqref{flaschka transformation} preserves the
presymplectic forms on $(\ell^\infty
/\mathbb{R}\mathbf{q}_0)
\times \ell^1_{0f}$  and therefore $J$ can be regarded as a momentum map
in a more general setting.

Additionally, let us mention that the Toda lattice Hamiltonian (\ref{Toda lattice Hamiltonian})
is of the form $H = h \circ J$, for 
\begin{equation}
\label{toda hamiltonian on orbits}
h(\rho_-) = \operatorname{tr}(\rho_- + a)^2,
\end{equation}
where $a: = \sum_{n=1}^\infty \alpha_k |k \rangle \langle k+1|
\in \one$, that is, $\sum_{n=1}^\infty |\alpha_k| <\infty$.
Thus, the Toda lattice is a Hamiltonian system on a Poisson submanifold of
$\one_-$ endowed with the bracket \eqref{Toda bracket} and relative to the
Hamiltonian function \eqref{toda hamiltonian on orbits}.

We shall prove below a version of an involution
theorem combining the Kostant-Symes and the Mischchenko-Fomenko
involution theorems for $\one_-$ with the goal to show that the functions
$h_k(\rho_-):= \operatorname{tr}(\rho_- +a)^k/k$, $k \in \mathbb{N}$, are
in involution relative to the Poisson bracket $\{\cdot,
\cdot\}_{\one_-}$. The proof turns out to follow the finite dimensional
one (see Kostant [1979] or Ratiu [1978]). Note that $f$ is a
Casimir function on $\one$ if and only if $X_f = 0$, which by \eqref{eq:ham} is
equivalent to $[df(\rho), \rho] = 0$ for every $\rho \in \one$. Decompose $\one =
\one_- \oplus \one^+$ and $\bounded = \bounded_+ \oplus \bounded^-$ 
where $\one^+ := \{\rho \in \one \mid \rho_{nm} = 0 \text{~for~} n
\geq m\}$ are the strictly upper triangular linear trace class operators
(no diagonal) and $\bounded^- := \{ X \in \bounded \mid X_{nm}=0 
\text{~for~} m\geq n\}$ are the strictly lower triangular bounded linear
operators (no diagonal). Let $\pi^1_-:\one \rightarrow \one_-$,
${\pi^1}^{+}:\one \rightarrow \one^+$, $\pi^\infty_+: \bounded
\rightarrow \bounded_+$, and ${\pi^\infty}^-: \bounded \rightarrow
\bounded^-$ be the projections associated to the Banach space direct
sums $\one =\one_- \oplus \one^+$ and $\bounded = \bounded_+ \oplus
\bounded^{-}$ respectively. With this notation the Poisson bracket
\eqref{Toda bracket} becomes
\begin{equation}
\label{Toda bracket with projections}
\{f, g\}_{\one_-} (\rho_-) =
\operatorname{tr}\left(\left[\pi^\infty_+(d\tilde{f}(\rho_-)),
\pi^\infty_+(d\tilde{g}(\rho_-))\right]\rho_-\right),
\end{equation}
where on the right hand side, $\tilde{f}$ and $\tilde{g}$ are arbitrary extensions
of $f,g : \one_- \rightarrow \mathbb{R}$ to $\one$
respectively. Thus $d\tilde{f}(\rho_-) \in \bounded$ and 
$\pi^\infty_+(d\tilde{f}(\rho_-)) \in \bounded_+$ and similarly for $g$.

\begin{proposition}
\label{involution theorem}
Let $a \in \one$ be a  given element satisfying 
\[
\operatorname{tr}(a[\bounded_+, \bounded_+]) = 0 \quad \text{and} \quad
\operatorname{tr}(a[\bounded^-, \bounded^-]) = 0.
\]
 For any two Casimir 
functions $f, g$ on $\one$ the functions
$f_a(\rho_-):=f(\rho_- + a)$,
$g_a(\rho_-):= g(\rho_- + a)$ are in involution on $\one_-$.
\end{proposition}

\noindent\textbf{Proof}.  Note that in \eqref{Toda bracket with projections} one
can take $\tilde{f}_a(\rho) = f(\rho +  a)$
for any $\rho \in \one$. Since $\operatorname{tr}(a[\bounded_+, \bounded_+]) = 0$, it
follows
\begin{align*}
\{f_a, g_a\}_{\one_-} (\rho_-) &=
\operatorname{tr}\left(\left[\pi^\infty_+(d\tilde{f}_a(\rho_-)),
\pi^\infty_+(d\tilde{g}_a(\rho_-))\right]\rho_- \right) \\
& = \operatorname{tr}\left(\left[\pi^\infty_+(df(\rho_- + a)),
\pi^\infty_+(dg(\rho_- + a))\right](\rho_- + a)\right) \\
& = \operatorname{tr}\left(\pi^\infty_+(df(\rho_- + a))
\left[dg(\rho_- + a) - {\pi^\infty}^{-}(dg(\rho_- + a)),
\rho_- + a \right]\right) \\
&= - \operatorname{tr}\left(\pi^\infty_+(df(\rho_- + a))
\left[{\pi^\infty}^{-}(dg(\rho_- + a)),
\rho_- + a \right]\right) \\
&= - \operatorname{tr}\left(\left[\pi^\infty_+(df(\rho_- + a)),
{\pi^\infty}^{-}(dg(\rho_- + a))\right](\rho_- + a) \right) \\
&= \operatorname{tr}\left(\left[{\pi^\infty}^{-}(dg(\rho_- + a)), 
\pi^\infty_+(df(\rho_- + a))\right](\rho_- + a) \right) \\
&= \operatorname{tr}\left({\pi^\infty}^{-}(dg(\rho_- + a))\left[ 
\pi^\infty_+(df(\rho_- + a)), \rho_- + a \right]\right) \\
&= \operatorname{tr}\left({\pi^\infty}^{-}(dg(\rho_- + a))\left[ 
df(\rho_- + a) - {\pi^\infty}^-(df(\rho_- + a)), 
\rho_- + a \right]\right) \\
&= -\operatorname{tr}\left({\pi^\infty}^{-}(dg(\rho_- + a))
\left[ {\pi^\infty}^-(df(\rho_- + a)), 
\rho_- + a \right]\right) \\
&= -\operatorname{tr}\left(\left[{\pi^\infty}^{-}(dg(\rho_- + a)),
{\pi^\infty}^-(df(\rho_- + a))\right] (\rho_- + a) \right) =0
\end{align*}
because $\operatorname{tr}(\rho_-[\bounded^-, \bounded^-]) =0$ and
$\operatorname{tr}(a[\bounded^-, \bounded^-]) =0$. \quad $\blacksquare$
\medskip

In the case of the Toda lattice one takes $a = \sum_{n=1}^\infty \alpha_k |k+1 
\rangle \langle k| \in \one$ with $\sum_{n=1}^\infty |\alpha_k| <\infty$ and
then it immediately follows that $\operatorname{tr}(a[\bounded^-, \bounded^-]) 
=0$ and $\operatorname{tr}(a[\bounded_+, \bounded_+]) = 0$. Thus, the 
hypotheses of Proposition \ref{involution theorem} are satisfied and we 
conclude that all the functions $h_k(\rho_-):= \operatorname{tr}(\rho_- +a)^k/k$, $k
\in \mathbb{N}$, are in involution relative to the Poisson bracket $\{\cdot,
\cdot\}_{\one_-}$ and hence the relation $J^\ast \omega_{\mathcal{O}_\Lambda} =
\omega$ shows that $h_k \circ J$ are commuting conserved quantities for the Toda
Hamiltonian $H : = h_2 \circ J$.
\quad $\blacklozenge$
\end{example} 

Next we discuss the Poisson reduction in a Banach Lie-Poisson
space $\mathfrak{b}$ associated to a quantum reduction operator $R
: \mathfrak{b} \rightarrow \mathfrak{b}$. Assume that $i: N
\hookrightarrow \mathfrak{b}$ is a (locally closed) submanifold.
Since $T\mathfrak{b} |_N = N \times \mathfrak{b}$, define the
subbundle $E \subset T\mathfrak{b} |_N$ by $E_b := \{b\} \times
\ker R$. Next, make the topological assumption that $E \cap TN$ is
the tangent bundle to a regular foliation $\mathcal{F}$ and that 
the space of leaves $M: = N/\mathcal{F}$ is a smooth manifold with
the projection $\pi :N \rightarrow M$ a submersion.

\begin{lemma}
The subbundle $E$ is compatible with the Poisson structure of $\mathfrak{b}$.
\end{lemma}

\noindent\textbf{Proof}. Let $U$ be an open subset of
$\mathfrak{b}$ and $f, g \in C^\infty(U, \mathbb{C})$ have the
property that $df, dg$ vanish on $E$, that is, $\langle df(b),
\ker R \rangle = 0$ and $\langle dg(b), \ker R \rangle = 0$ for $b
\in U \cap N$. Therefore, there exist functions $\tilde{f},
\tilde{g} \in C^\infty(R(U), \mathbb{C})$ such that $f = \tilde{f}
\circ R$ and $g = \tilde{g} \circ R$. Recall from 
\S \ref{section: symplectic leaves}
that the quantum reduction $R:
(\mathfrak{b}, \{\cdot,\cdot \} ) \rightarrow (\operatorname{im}R,
\{\cdot,\cdot\}_R)$ is a Poisson map. Thus, $\{f, g \} = \{\tilde{f} \circ
R, \tilde{g} \circ R\} = \{\tilde{f}, \tilde{g}\}_R \circ R$
whence $d\{f, g \}(b) = d\{\tilde{f}, \tilde{g}\}_R(R(b)) \circ R$
which implies that $d\{f, g \}(b)$ vanishes on $E_b = \ker R$.
\quad $\blacktriangledown$

\medskip

The following commutative diagram 

\unitlength=5mm
\begin{center}
\begin{picture}(9,8)
\put(1,7){\makebox(0,0){$N$}}
\put(9,7){\makebox(0,0){$\mathfrak{b}$}}
\put(1,2){\makebox(0,0){$N/\mathcal{F}$}}
\put(9,2){\makebox(0,0){$\operatorname{im}R$}}
\put(1,6){\vector(0,-1){3}} \put(9,6){\vector(0,-1){3}}
\put(2,7){\vector(1,0){5.5}}
\put(2,2){\vector(1,0){5.7}}
\put(0.3,4.3){\makebox(0,0){$\pi$}}
\put(9.4,4.3){\makebox(0,0){$R$}} \put(5,7.5){\makebox(0,0){$i$}}
\put(5,2.5){\makebox(0,0){$J$}}
\end{picture}
\end{center}
summarizes the maps involved in the theorem below.

\begin{theorem}
\label{quantum reduction theorem}
Let $i: N \hookrightarrow
\mathfrak{b}$ be a submanifold, $R: \mathfrak{b} \rightarrow
\mathfrak{b}$ be a quantum reduction, and $E$ be the distribution
on $\mathfrak{b}$ given at every point by $ker R$. Assume that:
\begin{itemize}
\item[{\rm (i)}] $E \cap TN$ is the tangent bundle of a regularfoliation
$\mathcal{F}$ on $N$ and the projection $\pi :N
\rightarrow M:= N/\mathcal{F}$ is a submersion;
\item[{\rm (ii)}] $\sharp(\ker R)^\circ \subset \overline{\ker R + T_n N}$
for every $n \in N$.
\end{itemize}
Then $M$ is the reduction of $\mathfrak{b}$ by $(N, E)$ and is
thus a Banach Poisson manifold. The map $J: M \rightarrow
\operatorname{im}R$ defined by $J([n]):= (R\circ i)(n)$ is
Poisson, that is, $J$ is a momentum map.
\end{theorem}

\noindent\textbf{Proof}. In view of the previous lemma, by Theorem
\ref{reduction theorem}, the two hypotheses guarantee that the
triple $(\mathfrak{b}, N, E)$ is reducible. Thus $M$ is a Banach
Poisson manifold.

Since $\operatorname{im} R$ can be regarded as the quotient
manifold obtained by collapsing the fibers of $R$, that is, by
dividing with $\ker R$, the inclusion map $i$ is obviously
compatible with the equivalence relations on $N$ and on
$\mathfrak{b}$. Therefore, $i$ induces a smooth map $J:M \rightarrow \operatorname{im}
R$ on the quotients (see, for example, Abraham, Marsden, Ratiu [1988] or Bourbaki
[1967]) given by $J([b]):= R(i(b))$. The diagram above commutes by construction.
It remains to be shown that $J$ is a Poisson map.

Let $f,g \in C^\infty(\operatorname{im}R, \mathbb{C})$. Then $f
\circ R \in C^\infty (\mathfrak{b}, \mathbb{C})$ is an extension
of $f\circ J \circ \pi \in C^\infty (N, \mathbb{C})$ and similarly
for $g$. Therefore, by the definition of the reduced bracket on
$M$, since $R$ is a Poisson map, we get
\[
\{f\circ J, g \circ J \}_M \circ \pi = \{f \circ R, g \circ R\}
\circ i = \{f, g \}_R \circ R \circ i= \{f, g\}_R \circ J \circ \pi.
\]
Since $\pi$ is a surjective map, this implies that $J: M
\rightarrow \operatorname{im}R$ is a Poisson map. \quad $\blacksquare$

\medskip

\begin{example} \textbf{Averaging}. \normalfont
Let $i :N \hookrightarrow \one$ be the inclusion map of a smooth
(regular) Banach submanifold in $\one$. Let $G$ be a compact Lie
group and denote by $\mu(g)$ the normalized Haar measure on $G$.
Given are:
\begin{itemize}
\item a smooth left action $\sigma: G \rightarrow
\operatorname{Diff}(N)$ of $G$ on $N$,
\item  a smooth Lie group homomorphism $U: G
\rightarrow \invertible$ such that $U(g)$ is unitary for each $g \in G$.
\end{itemize}
    
Assume also that the inclusion $i:N \hookrightarrow \one$ is
equivariant, that is, $i(\sigma(g)(n)) = U(g) i(n)
U(g)^{-1}$, for all $n \in N$ and all $g \in G$.

The homomorphism $U$ defines the operator $R: \one
\rightarrow \one$ by
\begin{equation}\label{averaging}
R(\rho) := \int_G U(g)\rho U(g)^\ast d\mu(g).
\end{equation}
We shall prove below that this $R$ is a quantum reduction operator.

We begin by showing that $R$ is a projector. By invariance of the
Haar measure under translations, we have for $\rho \in \one$,
\begin{align*}
R^2(\rho) &= \int_G U(h) \left(\int_G U(g)\rho U(g)^\ast d\mu(g)\right) U(h)^\ast
d\mu(h)\\ &= \int_G \left(\int_G
U(hg)\rho U(hg)^\ast d\mu(g)\right) d\mu(h)\\ &= \int_G R(\rho)
d\mu(g) = R(\rho).
\end{align*}

Next, we show that $\|R\|= 1$. Since $R$ is a projector we have
$\|R\| \leq 1$.  To prove equality, note first that if $\rho \geq
0$, that is, $\langle \psi | \rho |\psi \rangle \geq 0$ for all
$|\psi \rangle \in \mathcal{M}$, then we also have $ \langle \psi| U(g)\rho U(g)^\ast
|\psi \rangle \geq 0$ for all $|\psi \rangle
\in \mathcal{M}$ and integration over $G$ yields $R(\rho) \geq 0$.
Thus we showed that $\rho \geq 0$ implies $R(\rho) \geq 0$.  Continuity of the trace
in the $\|\cdot \|_1$--norm gives then for $\rho \geq 0$,
\[
\|R(\rho)\|_1 = \tr R(\rho) = \int_G \tr \left( U(g)\rho U(g)^\ast
\right) d\mu(g) =\int_G \tr \rho \,d\mu(g) = \tr\rho = \|\rho
\|_1,
\]
which proves that $\|R\|= 1$.

Finally, we need to show that $\operatorname{im}R^\ast$ is a
Banach Lie subalgebra of $\bounded$. It is easy to see that for
any $X \in \bounded$
\begin{equation*}
R^\ast(X) = \int_G U(g)^\ast X U(g) d\mu(g).
\end{equation*}
Using this formula we find
\[
R^\ast(X) R^\ast(Y) = R^\ast (X R^\ast(Y)) =  R^\ast(R^\ast(X) R^\ast(Y))
\]
and the condition that $\operatorname{im}R^\ast$ is a Lie
subalgebra of $\bounded$ follows immediately.

Thus all conditions in the definition of a quantum reduction
map are satisfied and hence $R$ given by \eqref{averaging} is a
quantum reduction operator. One can regard $R: \one \rightarrow \im R$ as a
momentum map. Consider the distribution $E$ on $N$
given at every point $n \in N$ by $E_n = \{n\} \times \ker R$.
Assume that $E\cap TN$ is the tangent bundle of a regular
foliation $\mathcal{F}$ on $N$ and that the projection $\pi : N
\rightarrow N/\mathcal{F}$ is a submersion. If $\sharp(\ker R)^\circ
\subset \overline{\ker R + T_n N}$ for every $n \in N$, the
conditions of Theorem \ref{quantum reduction theorem} are
satisfied and we obtain a momentum map $J: N/\mathcal{F}
\rightarrow \operatorname{im}R$. 

Note that the tanget spaces to the $G$-orbits determine a distribution on $N$ that is
included in the distribution $E \cap TN$. Supposing that the quotient by the group
action is regular, that is, that $N/G$ is a smooth manifold with the canonical
projection $N
\rightarrow N/G$ a surjective submersion, it follows that there is smooth surjective
map between quotient spaces $\Sigma : N/G \rightarrow N/\mathcal{F}$ and hence a map
$J \circ \Sigma : N/G\rightarrow
\im R$. If, in addition, $N/G$ is a Poisson manifold and $\Sigma$ is a Poisson map,
then $J \circ \Sigma : N/G\rightarrow \im R$ is a momentum map. This is satisfied,
for example, if $\Sigma = identity$.
\quad $\blacklozenge$
\end{example} 

Certain momentum maps play a special role in
the physical description of various systems. An important class of
such momentum maps are the coherent states maps. We shall
introduce this notion in the context of Banach Poisson manifolds,
modeling it on the definition introduced by Odzijewicz [1992] for
the case of a canonical map between a finite dimensional
symplectic manifold and the projectivization of a complex Hilbert
space.

\begin{definition}
\label{coherent states map definition} Let $P$ be a Banach Poisson
manifold and $\mathfrak{b}$ be a Banach Lie-Poisson space. A {\bfi
coherent states map\/} of $P$ into $\mathfrak{b}$ is a Poisson
embedding $\mathcal{K}: P \rightarrow \mathfrak{b}$ with linearly
dense range, that is, the span of $\operatorname{im} \mathcal{K}$
equals $\mathfrak{b}$.
\end{definition}

The situation investigated by Odzijewicz [1992] is the case when
$P$ is a finite dimensional Poisson manifold, $\mathfrak{b} =
\hermone$ is the Banach space of Hermitian trace class operators
on a separable complex Hilbert space $\mathcal{M}$,  and
$\mathcal{K}(p)$ is a rank one orthogonal projector for every $p
\in P$. In this case, the range of $\mathcal{K}$ lies in the
projectivization $\mathbb{CP}(\mathcal{M})$ of $\mathcal{M}$ by
identifying a rank one projector with the point in projective
space determined by its image. To illustrate this situation, let
us recall how $\mathcal{K}$ is used in the quantization of a
physical system. For the Poisson diffeomorphism $\sigma: P
\rightarrow P$, we assume that there is a linear Poisson
automorphism $\Sigma : \hermone \rightarrow \hermone$, such that
the diagram of canonical maps

\unitlength=5mm
\begin{center}
\begin{picture}(9,8)
\put(1,7){\makebox(0,0){$P$}}
\put(9,7){\makebox(0,0){$\hermone$}}
\put(1,2){\makebox(0,0){$P$}}
\put(9,2){\makebox(0,0){$\hermone$}}
\put(1,6){\vector(0,-1){3}}
\put(9,6){\vector(0,-1){3}}
\put(2,7){\vector(1,0){5.5}}
\put(2,2){\vector(1,0){5.7}}
\put(0.3,4.3){\makebox(0,0){$\sigma$}}
\put(9.4,4.3){\makebox(0,0){$\Sigma$}}
\put(5,7.5){\makebox(0,0){$\mathcal{K}$}}
\put(5,2.5){\makebox(0,0){$\mathcal{K}$}}
\end{picture}
\end{center}
commutes. By a theorem of Wigner, the automorphism $\Sigma$ is of
the form $\Sigma(\rho) = U \rho U^\ast$, where $U$ is a unitary or
anti-unitary operator on $\mathcal{M}$. Due to the hypothesis that
$\operatorname{im}\mathcal{K}$ is linearly dense in
$\mathfrak{b}$, if such an automorphism $\Sigma$ exists, it is
necessarily unique. It is natural to interpret $\Sigma$ as the
quantization of $\sigma$. Denote by $\operatorname{Aut}(\hermone)$
the linear Poisson isomorphisms of $\hermone$. The set of all
Poisson diffeomorphisms $\sigma$ for which a $\Sigma$ as above
exists, form a subgroup $\operatorname{Diff}_{\mathcal{K}}(P,
\{\cdot,\cdot\})$ of the Poisson diffeomorphism group
$\operatorname{Diff}(P, \{\cdot,\cdot\})$ of $P$. The map
\[
\mathcal{E}: \operatorname{Diff}_{\mathcal{K}}(P, \{\cdot, \cdot\})
\rightarrow \operatorname{Aut}(\hermone)
\]
so defined, is a group homomorphism which will be called,
according to Odzijewicz [1992], {\bfi Ehrenfest quantization}.

Consider now the flow $\sigma_t$ of the Hamiltonian vector field
$X_h$ on $M$ and assume that $\sigma_t \in
\operatorname{Diff}_{\mathcal{K}}(P, \{\cdot, \cdot\})$ for all $t$. It is
known that the set of all Hamiltonian functions satisfying this
condition form a Poisson subalgebra in the Poisson algebra of all
smooth functions on $P$. Then
\begin{equation}
\label{ehrenfest}
\mathcal{E}(\sigma_t)(\rho) = \exp (itH) \rho \exp(-itH)
\end{equation}
for $H$ a self-adjoint operator (unbounded, in general) whose
domain includes the linear span of the set
$\mathcal{K}(P)(\mathcal{M})$. The correspondence
\[
\mathcal{Q}: h \mapsto H
\]
defined above is linear and satisfies the relation
\[
\mathcal{Q}(\{h_1, h_2 \}) = i [\mathcal{Q}(h_1),
\mathcal{Q}(h_2)],
\]
that is, $\mathcal{Q}$ is the Lie algebra homomorphism induced by
$\mathcal{E}$. For more details on the precise relationship
between the coherent states map quantization and the
Kostant-Souriau geometric quantization as well as $*$-product
quantization, we refer to Odzijewicz [1988], [1992], [1994].

The traditional coherent states map sends classical states to pure
quantum states. Definition \ref{coherent states map definition}
generalizes this idea by letting the coherent states map to send
classical states to mixed quantum states. Mathematically, mixed
quantum states are in $\one$, or, more generally, in a Banach
Lie-Poisson space. Definition \ref{coherent states map definition}
further generalizes the usual approach by also allowing in this
scheme infinite dimensional classical systems.

\begin{example} 
\normalfont
(See Odzijewicz [1993])
 Consider the coherent states map $\mathcal{K}: P \rightarrow
\hermone$ from a finite dimensional Poisson manifold $(P, \{\cdot,
\cdot\})$ into the Banach Lie-Poisson space $\hermone$ of Hermitian trace
class operators on the separable complex Hilbert space
$\mathcal{M}$. Assume that the functions $f_1, \dots, f_k \in C^\infty(P)$
are in involution, that is, $\{f_i, f_j \} = 0$, for all $i,j = 1, \dots,
k$, and that their differentials $df_1(p), \dots, df_k(p)$ are linearly
independent for all $p \in N$, where $\iota: N \hookrightarrow P$ is a
given submanifold of $P$, invariant under the Hamiltonian flows
$\sigma_1(t), \dots, \sigma_k(t)$, $t \in \mathbb{R}$, generated by $f_1,
\dots , f_k$ respectively.

Quantize the flows $\sigma_1(t), \dots, \sigma_k(t)$ using the Ehrenfest 
quantization procedure \eqref{ehrenfest}, $\mathcal{E}: \sigma_i(t)
\mapsto \Sigma_i(t)$, where it is assumed that the generators $F_i$, $i =
1, \dots, k$ of the quantum flows $\Sigma_i(t)$ are all self adjoint
operators with discrete spectrum, i.e., 
\begin{equation}
\label{spectrum of f i}
F_i = \sum_{n=1}^\infty \lambda_n^i P_n^i \in \herm
\end{equation}
where $\{P_n^i\}_{n=1}^\infty$ is an orthonormal decomposition of the unit
related to $F_i$. 

Consider the orthogonal projector 
\[
P_{\lambda_{i_1}^1 \dots
\lambda_{i_k}^k}:  = P^1_{\lambda_{i_1}^1}  \dots
P^k_{\lambda_{i_k}^k}
\]
of $\mathcal{M}$ onto the common eigensubspace
\[
\mathcal{M}_{\lambda_{i_1}^1
\dots \lambda_{i_k}^k} := \{v \in \mathcal{M} \mid F_1 v= \lambda_{i_1}^1
v, \dots F_k v= \lambda_{i_k}^k v\}
\]
of the generators $F_1, \dots, F_k$. According to Example \ref{section 7,
first example}, the map $R_{\lambda_{i_1}^1 \dots
\lambda_{i_k}^k}: \hermone \rightarrow \hermone$ defined by
\begin{equation}
\label{r in last example}
R_{\lambda_{i_1}^1 \dots \lambda_{i_k}^k}(\rho): = P_{\lambda_{i_1}^1
\dots \lambda_{i_k}^k} \rho P_{\lambda_{i_1}^1 \dots \lambda_{i_k}^k}
\end{equation}
is a quantum reduction.

Assume now that conditions (i) and (ii) of Theorem \ref{quantum reduction
theorem} hold. In addition, assume that

\begin{enumerate}
\item[(iii)] the foliation $\mathcal{F}$ is given by the Hamiltonian
vector fields $X_{f_1}, \dots, X_{f_k}$.
\end{enumerate}

Then the quotient manifold $N/\mathcal{F}=:M$ is a Poisson manifold and
the map
\begin{equation}
\label{last example momentum map}
\mathcal{I}_{\lambda_{i_1}^1 \dots \lambda_{i_k}^k} : = 
R_{\lambda_{i_1}^1 \dots \lambda_{i_k}^k} \circ \iota
\end{equation}
is a momentum map of $M$ with values in $\operatorname{im}
R_{\lambda_{i_1}^1 \dots \lambda_{i_k}^k}$.

In the special case when $N = \mathbf{f}^{-1}(\mu)$, $\mu \in
\mathbb{R}^k$, is the level set of the map $\mathbf{f}: = (f_1, \dots,
f_k)$, conditions (i), (ii), and (iii) imply certain restrictions on
$\mu$. For example, it can happen that $\mu_1 = \lambda^1_{i_1}, \dots,
\mu_k =
\lambda^k_{i_k}$. In this case, the existence of
the momentum map
\eqref{last example momentum map}, i.e., the existence  of the Ehrenfest
quantization for the quotient system $N/\mathcal{F} = M$ leads to a
discretization (quantization) of $\mathbf{f}:P \rightarrow \mathbb{R}^k$.
The above method has been applied to the quantization of the 
MIC-Kepler system in Odzijewicz
and Swietochowski [1997]. \quad $\blacklozenge$
\end{example} 

The above examples show the importance of the relation between the
classical and quantum reduction procedures for  the
quantization and the integration of Hamiltonian systems. The present paper
raises several important questions regarding these connections to which we
shall return in future publications.

\bigskip

\addcontentsline{toc}{section}{Acknowledgments}
\noindent\textbf{Acknowledgments.} We thank A.B Antonewich, H. Flaschka,
J. Huebschmann, A.V. Lebedev J.-P. Ortega, and J. Mars\-den for several useful
discussions that improved the exposition. Special thanks to P. Bona for
his interest in our work and his inspired remarks. The first author was partially
supported by KBN under grant 2 PO3 A 012 19. The second author was
partially supported by the European Commission and the Swiss Federal
Government through funding for the Research Training Network
\emph{Mechanics and Symmetry in Europe} (MASIE) as well as the Swiss
National Science Foundation.

\bigskip

\section*{References}

\begin{description}

\item Abraham, R. and J.E. Marsden [1978]
{\it Foundations of Mechanics.\/} Second Edition, Addison-Wesley.

\item Abraham, R., Marsden, J.E., and  Ratiu, T.S. [1988]
{\it Manifolds, Tensor Analysis, and Applications.\/}  Second Edition, Applied 
Mathematical Sciences {\bf 75},
Springer-Verlag.

\item Accardi, L., Frigerio, A. and Gorini, V. [1984] \textit{Quantum
probability and Applications\/}, Lecture Notes in Math. {\bf 1136},
Springer-Verlag.

\item Accardi, L. and von Waldenfels [1988] \textit{Quantum
probability and Applications III\/}, Lecture Notes in Math. {\bf 1396},
Springer-Verlag.

\item Arnold, V.I. [1989]
{\it Mathematical Methods of Classical Mechanics.\/} Second Edition, Graduate 
Texts in Mathematics {\bf 60},
Springer-Verlag.

\item Bona, P. [2000] Extended quantum mechanics. {\it Acta
Physica Slovaca} {\bf 50} (1), 1 -- 198.

\item Boos, B. and Bleecker, D. [1983] {\it Topology and Analysis.
The Atiyah-Singer Index Formula and Gauge-Theoretic Physics}, Universitext, 
Springer-Verlag, New York.

\item Bourbaki, N. [1959] {\it Int\'egration\/}, Chapitre 6,
Hermann, Paris.

\item Bourbaki, N. [1967] {\it Vari\'et\'es diff\'erentielles et
analytiques. Fascicule de r\'esultats. Paragraphes 1 \`a 7},
Hermann, Paris.

\item Bourbaki, N. [1972] {\it Groupes et alg\`ebres de Lie},
Chapitre 3, Hermann, Paris.

\item Bratteli, O. and Robinson, D.W. [1979] \textit{ Operator Algebras
and Quantum Statistical Mechanics I\/}, Springer-Verlag.

\item Bratteli, O. and Robinson, D.W. [1981] \textit{ Operator Algebras
and Quantum Statistical Mechanics II\/}, Springer-Verlag.

\item Chernoff, P. R. and Marsden, J. E. [1974]
\textit{Properties of Infinite Dimensional Hamiltonian Systems}.
Lecture Notes in Mathematics, {\bf 425}. Springer Verlag.

\item Ehrenfest, P. [1927] {\it Z. Physik} {\bf 45}, 455--

\item Emch, G. [1972] \textit{Algebraic Methods in Statistical
Mechanics\/}, Wiley Interscience.

\item Holevo, A. [2001]  {\it Statistical Structure of Quantum
theory\/.} Lecture Notes in Physics. Monographs. Springer-Verlag.

\item Kirillov, A. A. [1993] The orbit method, II: Infinite-dimensional
Lie groups and Lie algebras, {\it Cont. Math.\/}, {\bf 145}, {\it
Representation Theory of Groups and Algebras\/}, 33--63.

\item Klauder, J.R. and  Skagerstem, Bo-Sure [1985] {\it Coherent
States -- Applications in Physics and Mathematical Physics\/.} World 
Scientific, Singapore.

\item Kobayashi, S. and Nomizu, K. [1963] {\it Foundations of
Differential Geometry.\/} Wiley

\item Kostant [1966] Orbits, symplectic structures and representation
theory, \textit{Proc, US-Japan Seminar on Diff. Geom., Kyoto. \/} Nippon
Hyronsha, Tokyo {\bf 77}.

\item Kostant, B. [1970] Quantization and unitary
representations, {\it Lecture Notes in Math.\/} {\bf 170}, 87--208.

\item Kostant, B. [1979] The solution to a generalized Toda lattice and 
representation theory, \textit{Advances in Math.\/} \textbf{34}, 195--338.

\item Landsman, N.P. [1998] {\it Mathematical Topics Between
Classical and Quantum Mechanics.\/} Springer Monographs in Mathematics. 
Springer-Verlag.

\item Leble S.B. and Czachor M. [1998] Darboux-integrable Liouville-von Neumann
equations, \textit{Phys. Rev. E} {\bf 58}, 7091.

\item Libermann, P. and  Marle, C.-M. [1987]
{\it Symplectic Geometry and Analytical Mechanics.\/} Kluwer Academic 
Publishers.

\item Lichnerowicz, A. [1988] Vari\'et\'es de Jacobi et espaces 
homog\`enes de contact complexes, \textit{Journ. Math. Pures et Appl.\/}
{\bf 67}, 131--173.

\item Lie, S. [1890] \textit{Theorie der Transformationsgruppen, Zweiter
Abschnitt.} Teubner.

\item Marsden, J.E. and  Ratiu, T.S. [1986]
Reduction of Poisson manifolds, {\it Lett. in Math. Phys.\/} {\bf 11},
161--170.

\item Marsden, J.E. and  Ratiu, T.S. [1994]  {\it
Introduction to Mechanics and Symmetry.\/} Texts in Applied Mathematics, {\bf  
17}, Second Edition, second printing 2003,
Springer-Verlag.

\item Murphy, G.J. [1990] $C^\ast$-{\it algebras and Operator
Theory}, Academic Press, San Diego.

\item Nunes da Costa, J.M. [1997] Reduction of complex Poisson manifolds,
Portugaliae Math. {\bf 54}, 467--476.

\item Odzijewicz, A. [1988] On reproducing
kernels and quantization of states, {\it Comm. Math. Phys.} {\bf 114},
577--597.

\item Odzijewicz, A. [1992] Coherent states
and geometric quantization, {\it Comm. Math. Phys.} {\bf 150}, 385--413.

\item Odzijewicz, A. [1993] Coherent states
for reduced phase spaces, {\it Quantization and Coherent States Methods}, S.T. 
Ali, I.M. Mladenov, A.
Odzijewicz, eds., 161--169. World Sci. Press, Singapore.

\item Odzijewicz, A. [1994] Covariant and contravariant
Berezin symbols of bounded operators, {\it Quantization and Infinite-
Dimensional Systems}, J.-P. Antoine, S.T.
Ali, W. Lisiecki, I.M. Mladenov, A. Odzijewicz, eds., 99--108. World Sci. 
Press, Singapore.

\item Odzijewicz, A. and Ratiu, T.S. [2003] The Banach Poisson geometry of the
infinite Toda lattice, preprint.

\item Odzijewicz, A. and Swietochowski,  M. [1997]  Coherent states
map for MIC-Kepler system, {\it J. Math. Phys.\/} {\bf 38}, 5010--5030.

\item Ratiu, T.S. [1978] Involution theorems, in \textit{Geometric Methods
in Mathematical Physics\/}, G. Kaiser and J.E. Marsden, eds., Springer
Lecture Notes \textbf{775}, 219--257.

\item Sakai, S. [1971] {\it $C^\ast$-Algebras and
$W^\ast$-Algebras.\/} Ergebnisse der Mathematik und ihrer
Grenzgebiete, {\bf  60}, 1998 reprint, Springer-Verlag.

\item Schr\"odinger, E. [1926] {\it Naturwissenschaften} {\bf
14}, 664

\item Stefan, P. [1974] Accessible sets, orbits and foliations with
singularities, {\it Proc. London Math. Soc.} {\bf 29}, 699--713.

\item Souriau, J.-M. [1966] Quantification g\'eom\'etrique, \textit{Comm,
Math. Phys.\/} {\bf 1}, 374--398.

\item Souriau, J.-M. [1967] Quantification g\'eom\'etrique. Applications.
\textit{Ann. Inst. H. Poincar\'e.\/} {\bf 6}, 311--341.

\item Takesaki, M. [1972] Conditional expectations in von Neumann
algebra, \textit{J. Funct. Anal.\/} {\bf 9}, 306--321.

\item Takesaki, M. [1979] {\it Theory of Operator Algebras I}.
Springer-Verlag.

\item Vaisman, I. [1994] {\it Lectures on the Geometry of
Poisson Manifolds}. Progress in Mathematics, {\bf 118}. Birkh\"auser
Verlag, 
Basel, 1994.

\item von Neumann, J. [1955] {\it Mathematical Foundations of Quantum 
Mechanics\/}, Princeton Univ. Press,
Princeton, N.J.

\item Weinstein, A. [1983] The local structure of Poisson manifolds, 
{\it Journ. Diff. Geom\/} {\bf 18}, 523--557.

\item Weinstein, A. [1998] Poisson geometry.
{\it Differential Geom. Appl.} {\bf 9}, 213--238.

\end{description}

\end{document}